\documentclass[11pt]{article}

\usepackage{amscd,amsmath,amssymb,fancyhdr}

\newcommand{\version}{version 2.0,\ \   Oct 30, 2012}

\renewcommand{\index}[2]{{}}

\def\eqref#1{(\ref{#1})}
\newcommand{\goth}{\frak}
\newcommand{\g}{{\frak g}}
\newcommand{\arrow}{{\:\longrightarrow\:}}
\newcommand{\Z}{{\Bbb Z}}
\newcommand{\C}{{\Bbb C}}
\newcommand{\R}{{\Bbb R}}
\newcommand{\Q}{{\Bbb Q}}
\newcommand{\6}{\partial}
\newcommand{\1}{\sqrt{-1}\:}
\newcommand{\restrict}[1]{{\left|_{{\phantom{|}\!\!}_{#1}}\right.}}

\renewcommand{\c}[1]{{\cal #1}}
\newcommand{\calo}{{\cal O}}

\let\oldtilde=\tilde

\renewcommand{\tilde}{\widetilde}
\renewcommand{\bar}{\overline}
\renewcommand{\phi}{\varphi}
\renewcommand{\epsilon}{\varepsilon}
\renewcommand{\geq}{\geqslant}
\renewcommand{\leq}{\leqslant}

\newcommand{\fl}{{\rm fl}}
\newcommand{\im}{\operatorname{im}}
\newcommand{\End}{\operatorname{End}}
\newcommand{\Tot}{\operatorname{Tot}}
\newcommand{\Id}{\operatorname{Id}}
\newcommand{\id}{\operatorname{\text{\sf id}}}
\newcommand{\Vol}{\operatorname{Vol}}

\newcommand{\Aut}{\operatorname{Aut}}
\newcommand{\codim}{\operatorname{codim}}
\newcommand{\coim}{\operatorname{coim}}
\newcommand{\coker}{\operatorname{coker}}
\newcommand{\slope}{\operatorname{slope}}
\newcommand{\rk}{\operatorname{rk}}
\newcommand{\Lin}{\operatorname{Lin}}
\newcommand{\Hor}{\operatorname{Hor}}
\newcommand{\Def}{\operatorname{Def}}
\newcommand{\Tw}{\operatorname{Tw}}

\newcommand{\Spec}{\operatorname{Spec}}

\newcommand{\comment}[1]{{}}

\def\blacksquare{\hbox{\vrule width 4pt height 4pt depth 0pt}}
\def\endproof{\blacksquare}

\makeatletter

\@ifundefined{Bbb}
     {\newcommand{\Bbb}[1]{{\mathbb #1}}}%
{}

\newcounter{Mycounter}[section]
\newcounter{lemma}[section]
\renewcommand{\thelemma}{{Lemma \thesection.\arabic{lemma}}}
\newcommand{\lemma}{%
     \setcounter{lemma}{\value{Mycounter}}
     \refstepcounter{lemma}
     \stepcounter{Mycounter}
     {\bf \thelemma:\ }}

\newcounter{claim}[section]
\renewcommand{\theclaim}{{Claim \thesection.\arabic{claim}}}
\newcommand{\claim}{%
     \setcounter{claim}{\value{Mycounter}}
     \refstepcounter{claim}
     \stepcounter{Mycounter}
     {\bf \theclaim:\ }}

\newcounter{sublemma}[section]
\renewcommand{\thesublemma}{{Sublemma \thesection.\arabic{sublemma}}}
\newcommand{\sublemma}{%
     \setcounter{sublemma}{\value{Mycounter}}
     \refstepcounter{sublemma}
     \stepcounter{Mycounter}
     {\bf \thesublemma:\ }}

\newcounter{corollary}[section]
\renewcommand{\thecorollary}{{Corollary \thesection.\arabic{corollary}}}
\newcommand{\corollary}{%
     \setcounter{corollary}{\value{Mycounter}}
     \refstepcounter{corollary}
     \stepcounter{Mycounter}
     {\bf \thecorollary:\ }}

\newcounter{theorem}[section]
\renewcommand{\thetheorem}{{Theorem \thesection.\arabic{theorem}}}
\newcommand{\theorem}{%
     \setcounter{theorem}{\value{Mycounter}}
     \refstepcounter{theorem}
     \stepcounter{Mycounter}
     {\bf \thetheorem:\ }}

\newcounter{conjecture}[section]
\renewcommand{\theconjecture}{{Conjecture \thesection.\arabic{conjecture}}}
\newcommand{\conjecture}{%
     \setcounter{conjecture}{\value{Mycounter}}
     \refstepcounter{conjecture}
     \stepcounter{Mycounter}
     {\bf \theconjecture:\ }}

\newcounter{proposition}[section]
\renewcommand{\theproposition}
       {{Proposition \thesection.\arabic{proposition}}}
\newcommand{\proposition}{%
     \setcounter{proposition}{\value{Mycounter}}
     \refstepcounter{proposition}
     \stepcounter{Mycounter}
     {\bf \theproposition:\ }}

\newcounter{definition}[section]
\renewcommand{\thedefinition}
       {{Definition \thesection.\arabic{definition}}}
\newcommand{\definition}{%
     \setcounter{definition}{\value{Mycounter}}
     \refstepcounter{definition}
     \stepcounter{Mycounter}
     {\bf \thedefinition:\ }}

\newcounter{example}[section]
\renewcommand{\theexample}{{Example \thesection.\arabic{example}}}
\newcommand{\example}{%
     \setcounter{example}{\value{Mycounter}}
     \refstepcounter{example}
     \stepcounter{Mycounter}
     {\bf \theexample:\ }}

\newcounter{remark}[section]
\renewcommand{\theremark}{{Remark \thesection.\arabic{remark}}}
\newcommand{\remark}{%
     \setcounter{remark}{\value{Mycounter}}
     \refstepcounter{remark}
     \stepcounter{Mycounter}
     {\bf \theremark:\ }}

\newcounter{problem}[section]
\newcounter{question}[section]
\renewcommand{\thequestion}{{Question \thesection.\arabic{question}}}
\newcommand{\question}{%
     \setcounter{question}{\value{Mycounter}}
     \refstepcounter{question}
     \stepcounter{Mycounter}
     {\bf \thequestion:\ }}

\@addtoreset{equation}{section}
\@addtoreset{footnote}{section}
\makeatother

\begin{document}

\begin{center}
{\Large\bf
Hyperholomorphic sheaves \\ and
new examples of hyperk\"ahler manifolds}\\[4mm]
Misha Verbitsky,\\[4mm]
{\tt verbit@thelema.dnttm.rssi.ru, verbit@ihes.fr}
\end{center}

\hfill

{\small 
\hspace{0.2\linewidth}
\begin{minipage}[t]{0.7\linewidth}
Given a compact hyperk\"ahler manifold $M$ and a holomorphic bundle
$B$ over $M$, we consider a Hermitian connection $\nabla$ on $B$
which is compatible with all complex structures on $M$ induced by
the hyperk\"ahler structure.  Such a connection is unique, because it is 
Yang-Mills.  We call the bundles admitting such connections
{\bf hyperholomorphic bundles}. A bundle is hyperholomorphic if and
only if its Chern classes $c_1$, $c_2$ are $SU(2)$-invariant, with
respect to the natural $SU(2)$-action on the cohomology. For several
years, it was known that the moduli space of stable hyperholomorphic
bundles is singular hyperk\"ahler.  More recently, it was proven that
singular hyperk\"ahler varieties admit a canonical hyperk\"ahler
desingularization.  In the present paper, we show that a moduli
space of stable hyperholomorphic bundles is compact, given some
assumptions on Chern classes of $B$ and hyperk\"ahler geometry of $M$
(we also require $\dim_{\C}M>2$).  Conjecturally, this leads to
new examples of hyperk\"ahler manifolds. We develop the theory of
hyperholomorphic sheaves, which are (intuitively speaking)
coherent sheaves compatible with hyperk\"ahler structure. We show
that \index{terms}{coherent sheaves!hyperholomorphic} hyperholomorphic 
sheaves with isolated singularities can be
canonically desingularized by a blow-up. This theory is used to
study degenerations of hyperholomorphic bundles.
\end{minipage}
}
\index{terms}{moduli space!of stable hyperholomorphic bundles} 
\index{terms}{vector bundles!stable} 
\index{terms}{connections!Yang--Mills} 

\tableofcontents


\section{Introduction}
\label{_intro_Section_}


\hfill

For an introduction to basic results and the history of
hyperk\"ahler geometry, see \cite{_Besse:Einst_Manifo_}. 
\index{names}{Besse, A.}

\hfill

This Introduction is independent from the rest of this paper.

\subsection{An overview}

\subsubsection{Examples of hyperk\"ahler manifolds}

 A Riemannian manifold $M$ is called {\bf hyperk\"ahler}
if the tangent bundle of $M$ is equipped with an action of
quaternian algebra, and its metric is K\"ahler with respect to the
complex structures $I_\iota$, for all 
embeddings $\C \stackrel{\iota}\hookrightarrow \Bbb H$.
The complex structures $I_\iota$ are called 
{\bf induced complex structures}; the corresponding
\index{terms}{induced complex structures}
K\"ahler manifold is denoted by $(M, I_\iota)$.

For a more formal definition of a hyperk\"ahler manifold, see 
\ref{_hyperkahler_manifold_Definition_}.
The notion of a hyperk\"ahler manifold was introduced by E. \index{names}{Calabi, E.} Calabi
(\cite{_Calabi_}).

Clearly, the real dimension of $M$ is divisible by 4.
For $\dim_\R M= 4$, there are only two classes of compact 
hyperk\"ahler manifolds: compact tori and K3 surfaces. 

Let $M$ be a complex surface and $M^{(n)}$ be its $n$-th 
symmetric power, $M^{(n)} = M^n/S_n$. The variety $M^{(n)}$
admits a natural desingularization  $M^{[n]}$, called
{\bf the Hilbert scheme of points}. 

The manifold $M^{[n]}$ admits a hyperk\"ahler metrics
whenever the surface $M$ is compact and hyperk\"ahler 
(\cite{_Beauville_}). This way, \index{names}{Beauville, A.} Beauville constructed 
two series of examples of hyperk\"ahler manifolds,
associated with a torus (so-called 
``higher Kummer variety'')
and a K3 surface. It was conjectured that
all compact hyperk\"ahler manifolds $M$ with $H^1(M) =0$, $H^{2,0}(M)=\C$
are deformationally equivalent to one of these examples. 
In this paper, we study the deformations of coherent sheaves 
over higher-dimensional hyperk\"ahler manifolds 
in order to construct counterexamples to 
this conjecture. A different approach to the 
construction of new examples of hyperk\"ahler manifolds
is found in the recent paper of K. O'Grady, who studies the moduli of 
\index{names}{O'Grady, K.} 
\index{terms}{vector bundles!semistable!moduli of} 
semistable bundles over a K3 surface and resolves the 
singularities using methods of symplectic geometry (\cite{_O'Grady_}).

\subsubsection{Hyperholomorphic bundles}

Let $M$ be a compact hyperk\"ahler manifold, and $I$ an induced
complex structure. It is well known that the differential
forms and cohomology
of $M$ are equipped with a natural $SU(2)$-action
(\ref{_SU(2)_commu_Laplace_Lemma_}). 
In \cite{_Verbitsky:Hyperholo_bundles_}, we studied the
holomorphic vector bundles $F$ on $(M,I)$
which are compatible with a hyperk\"ahler
\index{terms}{hyperk\"ahler structures} 
structure, in the sense that any of the following
conditions hold:
\begin{equation}\label{_hyperho_condi_Equation_}
\begin{minipage}[m]{0.8\linewidth}
\begin{description}
\item[(i)] The bundle $F$ admits a Hermitian connection $\nabla$ 
with a
curvature $\Theta\in \Lambda^2(M, \End(F))$ which is
of Hodge type (1,1) with respect to any of induced complex
structures.
\item[(ii)] The bundle $F$ is a direct sum of \index{terms}{vector bundles!stable} stable bundles,
and its Chern classes $c_1(F)$, $c_2(F)$ are $SU(2)$-invariant.
\end{description}
\end{minipage}
\end{equation}

These conditions are equivalent
(\ref{_inva_then_hyperho_Theorem_}). Moreover, the
connection $\nabla$ of \eqref{_hyperho_condi_Equation_} (i)
is \index{terms}{vector bundles!Yang--Mills} 
Yang-Mills (\ref{_hyperholo_Yang--Mills_Proposition_}), 
and by \index{terms}{Uhlenbeck--Yau Theorem} Uhlenbeck--\index{names}{Yau, S.-T.}Yau theorem (\ref{_UY_Theorem_}), 
it is unique.

A holomorphic vector bundle satisfying any of the
conditions of \eqref{_hyperho_condi_Equation_}
is called {\bf \index{terms}{vector bundles!hyperholomorphic!definition of} 
hyperholomorphic}
(\cite{_Verbitsky:Hyperholo_bundles_}).

Clearly, a stable deformation of a 
\index{terms}{vector bundles!hyperholomorphic!deformations of} 
hyperholomorphic bundle
is again a hyperholomorphic bundle. In  
\cite{_Verbitsky:Hyperholo_bundles_}, we proved that a
deformation space of hyperholomorphic bundles is a singular
hyperk\"ahler variety. A recent development in the theory
of singular hyperk\"ahler varieties 
(\cite{_Verbitsky:Desingu_}, 
\cite{_Verbitsky:DesinguII_}, \cite{_Verbitsky:hypercomple_}) 
gave a way to \index{terms}{hyperk\"ahler desingularization} desingularize singular hyperk\"ahler manifolds,
in a canonical way. It was proven (\ref{_desingu_Theorem_}) 
that a normalization of a singular hyperk\"ahler variety
(taken with respect to any induced complex structure $I$)
is a smooth hyperk\"ahler manifold. 

This suggested a possibility of constructing new examples of
compact hyperk\"ahler manifolds, obtained as deformations of
hyperholomorphic bundles. Two problems arise.

\hfill

{\bf Problem 1.} The deformation space 
of hyperholomorphic bundles is {\it a priori} non-compact
and must be compactified.

\hfill

{\bf Problem 2.} The geometry of deformation spaces is notoriously
hard to study. Even the dimension of a deformation space
is difficult to compute, in simplest examples. How to 
find, for example, the dimension of the deformation
space of a tangent bundle, on a Hilbert scheme of points on a K3
surface? The Betti numbers are even more difficult to compute.
Therefore, there is no easy way to distinguish 
a deformation space of 
\index{terms}{vector bundles!hyperholomorphic!moduli of} 
hyperholomorphic bundles from already
known examples of hyperk\"ahler manifolds.

\hfill

In this paper, we address Problem 1. Problem 2 can be solved 
by studying the algebraic geometry of moduli spaces.
It turns out that, for a generic deformation of a complex
structure, the Hilbert scheme of points on a K3 surface has no 
closed complex subvarieties (\cite{_Verbitsky:Hilbert_};
see also \ref{_no_triana_subva_of_Hilb_Theorem_}).
It is possible to find a 21-dimensional family
of deformations of the moduli space $\Def(B)$ of hyperholomorphic
bundles, with all fibers having complex subvarieties
(\ref{_double_Fou_embedding_Lemma_}).
Using this observation, it is possible to show that
 $\Def(B)$ is a new example of a hyperk\"ahler manifold.
Details of this approach are given in
Subsection \ref{_new_exa_F-M_checking_Subsection_},
and the complete proofs will be given in a forthcoming paper.

It was proven
that a Hilbert scheme of a generic K3 surface
has no \index{terms}{trianalytic subvarieties} trianalytic subvarieties.\footnote{Trianalytic subvariety
(\ref{_trianalytic_Definition_}) is a closed subset which is
complex analytic with respect to any of induced complex structures.}
Given a hyperk\"ahler manifold $M$ and an appropriate 
hyperholomorphic bundle $B$, denote the
deformation space of \index{terms}{connections!hyperholomorphic} hyperholomorphic connections
on $B$ by $\Def(B)$. Then the moduli of complex structures
on $M$ are locally embedded to a moduli of complex structures
on $\Def(B)$ (\ref{_maps_pre_tw_curves_Claim_}).
Since the dimension of the moduli of complex
structures on $\Def(B)$ is equal to 
its second Betti number minus 2 
(\ref{_Bogomo_etale_Theorem_}), the second 
Betti number of $\Def(B)$ is no less than the second Betti number
of $M$. The Betti numbers of \index{names}{Beauville, A.} Beauville's examples of 
simple hyperk\"ahler manifolds are 23 
(Hilbert scheme of points on a K3 surface) 
and 7 (generalized Kummer variety). Therefore, 
for $M$ a generic deformation of a 
Hilbert scheme of points on K3, $\Def(B)$ is either a new
manifold or a generic deformation of a Hilbert scheme of points on K3.
It is easy to construct \index{terms}{trianalytic subvarieties} trianalytic subvarieties of
the varieties $\Def(B)$, for hyperholomorphic $B$
(see \cite{_Verbitsky:Symplectic_I_}, 
Appendix for details). This was the motivation of our work on trianalytic
subvarieties of the Hilbert scheme of points on a K3 surface
(\cite{_Verbitsky:Hilbert_}). 

For a generic complex structure on a hyperk\"ahler
manifold, all \index{terms}{vector bundles!stable} 
stable bundles are hyperholomorphic 
(\cite{_Verbitsky:Symplectic_I_}).
Nevertheless, 
\index{terms}{vector bundles!hyperholomorphic!over higher-dimensional hyperk\"ahler manifolds} 
hyperholomorphic bundles over 
higher-dimensional hyperk\"ahler manifolds are in short supply.
In fact, the only example to work with is the
tangent bundle and its tensor powers, and their
Chern classes are not prime. 
Therefore, there is no way to insure that
a deformation of a stable bundle will remain 
stable (like it happens, for instance, in the case of 
deformations of stable bundles
of rank 2 with odd first Chern class
over a K3 surface). Even worse, a new
kind of singularities may appear which never appears
for 2-dimensional base manifolds: a deformation
of a \index{terms}{vector bundles!stable!deformations of} 
stable bundle can have
a singular \index{terms}{reflexization}
reflexization. We study the singularities
of stable coherent sheaves over hyperk\"ahler manifolds, 
using Yang-Mills theory for 
\index{terms}{coherent sheaves!reflexive}
reflexive
sheaves developed by S. Bando and Y.-T. \index{names}{Siu, Y.-T.} Siu 
(\cite{_Bando_Siu_}).

\subsubsection{Hyperholomorphic sheaves}

A compactification of the moduli of 
\index{terms}{vector bundles!hyperholomorphic!moduli of} 
hyperholomorphic bundles
is the main purpose of this paper. We require the compactification
to be singular hyperk\"ahler. A natural approach to this problem
requires one to study the coherent sheaves which are
compatible with a hyperk\"ahler structure, in the same
sense as \index{terms}{vector bundles!hyperholomorphic} hyperholomorphic bundles are holomorphic
bundles compatible with a hyperk\"ahler structure. Such
sheaves are called {\bf hyperholomorphic sheaves}
(\ref{_hyperho_shea_Definition_}). Our approach to the
theory of hyperholomorphic sheaves uses the notion of
admissible \index{terms}{connections!Yang--Mills!in coherent sheaves} 
Yang-Mills connection on a coherent sheaf
(\cite{_Bando_Siu_}). 

The equivalence of conditions \eqref{_hyperho_condi_Equation_} (i)
and \eqref{_hyperho_condi_Equation_} (ii) is based on \index{terms}{Uhlenbeck--Yau Theorem} Uhlen\-beck--\-\index{names}{Yau, S.-T.}Yau
theorem (\ref{_UY_Theorem_}), which states that every
\index{terms}{coherent sheaves!stable!singularities of} 
stable bundle
$F$ with $\deg c_1(F) =0$ admits a unique \index{terms}{connections!Yang--Mills} Yang--Mills connection,
that is, a connection $\nabla$ satisfying $\Lambda\nabla^2=0$
(see Subsection \ref{_sta_bu_and_YM_Subsection_} for details). 
S. Bando and Y.-T. \index{names}{Siu, Y.-T.} Siu developed a similar approach to the
Yang--Mills theory on (possibly singular) coherent sheaves.
Consider a coherent sheaf $F$ and a Hermitian metric 
$h$ on a locally trivial part of $F\restrict{U}$. Then $h$ is called
admissible (\ref{_admi_metri_Definition_}) if the curvature $\nabla^2$
of the Hermitian connection on $F\restrict{U}$ is square-integrable,
and the section $\Lambda\nabla^2\in \End(F\restrict{U})$ is
uniformly bounded. The admissible metric is called {\bf 
\index{terms}{metric!Yang--Mills!admissible} 
Yang-Mills}
if $\Lambda\nabla^2=0$ (see \ref{_Yang-Mills_sheaves_Definition_} 
for details). There exists an analogue of \index{terms}{Uhlenbeck--Yau Theorem} Uhlenbeck--\index{names}{Yau, S.-T.}Yau theorem
for coherent sheaves (\ref{_UY_for_shea_Theorem_}): a
\index{terms}{coherent sheaves!stable!connections on} 
stable sheaf
admits a unique admissible Yang--Mills metric, and conversely,
a sheaf admitting a \index{terms}{metric!Yang--Mills} Yang--Mills metric is a direct sum of stable
sheaves with the first Chern class of zero degree.

A coherent sheaf $F$ is called {\bf 
\index{terms}{coherent sheaves!reflexive}
reflexive} if it is isomorphic
to its second dual sheaf $F^{**}$. The sheaf $F^{**}$ is always
reflexive, and it is called {\bf a \index{terms}{reflexization}
reflexization} of $F$
(\ref{_refle_Definition_}). 

Applying the arguments of Bando and \index{names}{Siu, Y.-T.} Siu 
to a reflexive coherent sheaf $F$ over a hyperk\"ahler 
manifold $(M, I)$, we show that the following 
conditions are equivalent (\ref{_hyperho_conne_exi_Theorem_}).

\begin{description}
\item[(i)] The sheaf $F$ is stable and its Chern classes
$c_1(F)$, $c_2(F)$ are $SU(2)$-invariant
\item[(ii)] $F$ admits an admissible 
\index{terms}{connections!Yang--Mills!admissible} Yang--Mills connection,
and its curvature is of type (1,1) with respect to all
induced complex structures.
\end{description}

A reflexive sheaf satisfying any of the these conditions is called
{\bf reflexive 
\index{terms}{coherent sheaves!reflexive}
\index{terms}{coherent sheaves!hyperholomorphic!stable} 
\index{terms}{coherent sheaves!hyperholomorphic!semistable} 
stable 
hyperholomorphic}. An arbitrary torsion-free coherent sheaf is
called {\bf stable hyperholomorphic} if its \index{terms}{reflexization}
reflexization is
hyperholomorphic, and its second Chern class is $SU(2)$-invariant,
and {\bf semistable hyperholomorphic} if it is a successive extension of
stable hyperholomorphic sheaves (see 
\ref{_hyperho_shea_Definition_} for details).

This paper is dedicated to the study of hyperholomorphic sheaves.

\subsubsection{Deformations of hyperholomorphic sheaves}

By \ref{_generic_are_dense_Proposition_}, for an induced 
complex structure $I$ of general type, {\bf all} coherent sheaves
are 
\index{terms}{coherent sheaves!hyperholomorphic!over generic complex structures} hyperholomorphic. 
However, the complex structures of general
type are never algebraic, and in complex analytic situation,
the 
\index{terms}{moduli space!of coherent sheaves} 
moduli of coherent sheaves are, generaly speaking, non-compact.
We study the flat deformations of 
\index{terms}{coherent sheaves!hyperholomorphic!flat deformations of} 
hyperholomorphic sheaves
over $(M,I)$, where $I$ is an algebraic complex structure.

{\it A priori}, a flat deformation of a hyperholomorphic sheaf will
be no longer hyperholomorphic. We show that 
for some algebraic complex structures, called {\bf \index{terms}{$C$-restricted complex structures} $C$-restricted
complex structures}, a flat deformation of a hyperholomorphic sheaf 
remains hyperholomorphic (\ref{_sheaf_on_C_restr_hyperho_Theorem_}). 
This argument is quite convoluted, and takes two sections (Sections
\ref{_cohomo_hype_Section_} and \ref{_C_restri_Section_}).

Further on, we study the local structure of 
\index{terms}{coherent sheaves!stable} 
stable reflexive 
\index{terms}{coherent sheaves!hyperholomorphic!with isolated singularities} 
hyperholomorphic sheaves with isolated singularities.
We prove the Desingularization Theorem for such hyperholomorphic
\index{terms}{coherent sheaves!hyperholomorphic!desingularization of} 
sheaves (\ref{_desingu_hyperho_Theorem_}). It turns out that
such a sheaf can be desingularized by a single blow-up.
The proof of this result is parallel to the proof of
Desingularization Theorem for singular hyperk\"ahler varieties
(\ref{_desingu_Theorem_}). 

The main idea of the \index{terms}{hyperk\"ahler desingularization} desingularization of singular hyperk\"ahler
varieties (\cite{_Verbitsky:DesinguII_}) is the following.
Given a point $x$ on a singular hyperk\"ahler variety $M$
and an induced complex structure $I$, the complex variety
$(M, I)$ admits a local 
\index{terms}{$\C^*$-action on hyperk\"ahler manifolds} 
$\C^*$-action which preserves $x$ and
acts as a dilatation on the Zariski tangent space of $x$.
Here we show that any stable 
\index{terms}{coherent sheaves!hyperholomorphic!$\C^*$-equivariant} 
hyperholomorphic sheaf $F$ is equivariant
with respect to this $\C^*$-action 
(\ref{_Psi_equiv_hyperho_Theorem_},
\ref{_C^*_stru_on_sge_Definition_}). 
Then an elementary algebro-geometric argument
(\ref{_desingu_C^*_equi_Proposition_}) 
implies that $F$ is desingularized by a blow-up.

Using the desingularization of hyperholomorphic sheaves, we prove
that a hyperholomorphic deformation of a hyperholomorphic bundle
is again a bundle (\ref{_reflexi_defo_loca_trivi_Theorem_}),
assuming that it has isolated singularities. 
The proof of this result is conceptual but 
quite difficult, it takes 3 sections
(Sections 
\ref{_twisto_tra_Section_}--\ref{_modu_hyperho_Section_}), 
and uses arguments of
\index{terms}{quaternionic-K\"ahler manifolds} quaternionic-K\"ahler geometry (\cite{_Swann_},
\cite{_Nitta:Y-M_}) and twistor \index{terms}{twistor transform} transform (\cite{_NHYM_}).

In our study of deformations of hyperholomorphic
sheaves, we usually assume
that a deformation of a hyperholomorphic sheaf over
$(M, I)$ is again hyperholomorphic,
i. e. that an induced complex structure $I$ is \index{terms}{$C$-restricted complex structures} $C$-restricted,
for $C$ sufficiently big (\ref{_C_restri_Definition_}). 
Since $C$-restrictness is a tricky condition, it is preferable to
get rid of it. For this purpose, we use the theory of twistor
paths, developed in \cite{_coho_announce_}, to show
that the moduli spaces of hyperholomorphic sheaves are 
real analytic equivalent for different complex structures
$I$ on $M$ (\ref{_iso_Bun_exists_gene_pola_Theorem_}). 
This is done as follows.

A hyperk\"ahler structure on $M$ admits a 2-dimensional sphere of induced
complex structures. This gives a rational curve in the moduli space
$Comp$ of complex structures on $M$, so-called {\bf 
\index{terms}{twistor curves!definition of} 
twistor curve}. 
A sequence of such rational curves connect any two points of $Comp$ 
(\ref{_twistor_connect_Theorem_}).
A sequence of connected twistor curves is called {\bf a twistor
path}. If the intersection points of these curves
are generic, the \index{terms}{twistor paths} 
twistor path is called {\bf admissible}
(\ref{_admi_twi_path_Definition_}). It is known
(\ref{_admi_twi_impli_Theorem_}) that an admissible twistor
path induces a real analytic isomorphism of the moduli spaces
of \index{terms}{vector bundles!hyperholomorphic!moduli of} 
hyperholomorphic bundles. There exist 
admissible twistor paths connecting 
any two complex structures (\ref{_admi_pa_exist_for_gene_pol_Claim_}).
Thus, if we prove that the moduli
\index{terms}{moduli space!of stable hyperholomorphic bundles} 
of deformations of hyperholomorphic bundles
are compact for one generic hyperk\"ahler 
\index{terms}{hyperk\"ahler structures!generic} 
structure, we prove a similar result for all generic
hyperk\"ahler structures (\ref{_iso_Bun_exists_gene_pola_Theorem_}).
Applying this argument to the moduli of deformations
\index{terms}{moduli space!of deformations of tangent bundle} 
of a tangent bundle, we obtain the following theorem.

\hfill

\theorem \label{_defo_tange_Hilb_compact_intro_Theorem_}
Let $M$ be a Hilbert scheme of points on a K3 surface,
$\dim_{\Bbb H}(M)>1$ and $\c H$ a generic
hyperk\"ahler \index{terms}{hyperk\"ahler structures!generic} 
structure on $M$. Assume
that for all induced complex structures $I$, except at most 
a finite many of, all semistable bundle deformations of the tangent bundle
$T(M, I)$ are 
\index{terms}{vector bundles!stable} 
\index{terms}{vector bundles!semistable} 
stable. 
Then, for all complex structures $J$ on $M$ and 
all generic polarizations $\omega$
on $(M, J)$, the deformation space 
$\c M_{J, \omega}(T(M, J))$ is singular 
hyperk\"ahler and compact, and admits
a smooth compact hyperk\"ahler \index{terms}{hyperk\"ahler desingularization} desingularization.

{\bf Proof:} This is \ref{_defo_tange_compact_Theorem_}. 
\endproof

\hfill

In the course of this paper, we develop the theory of 
\index{terms}{$C$-restricted complex structures} $C$-restricted complex structures (Sections
\ref{_cohomo_hype_Section_} and \ref{_C_restri_Section_})
and another theory, which we called {\bf the \index{names}{Swann, A.} 
\index{terms}{Swann's formalism!for vector bundles} Swann's formalism
for vector bundles} (Sections 
\ref{_twisto_tra_Section_} and \ref{_C_equiv_twi_spa_Section_}).
These themes are of independent interest. We give a separate
introduction to $C$-restricted complex structures
(Subsection \ref{_C_restri_intro_Subsection_}) 
and Swann's formalism (Subsection \ref{_Swann's_intro_Subsection_}).

\subsection{$C$-restricted complex structures: an introduction}
\label{_C_restri_intro_Subsection_}

This part of the Introduction is highly non-precise.
Our purpose is to clarify the intuitive meaning
of \index{terms}{$C$-restricted complex structures} $C$-restricted complex structure. 

Consider a compact hyperk\"ahler manifold $M$,
which is {\bf simple}
(\ref{_simple_hyperkahler_mfolds_Definition_}),
that is, satisfies $H^1(M) =0$, $H^{2,0}(M) = \C$.

A {\bf reflexive 
\index{terms}{coherent sheaves!hyperholomorphic!definition of} 
hyperholomorphic sheaf} is by definition a semistable
sheaf which has a filtration of stable sheaves with
$SU(2)$-invariant $c_1$ and $c_2$. A {\bf hyperholomorphic
sheaf} is a torsion-free sheaf which has hyperholomorphic
\index{terms}{reflexization}
reflexization and has $SU(2)$-invariant $c_2$
(\ref{_hyperho_shea_Definition_}). If the 
complex structure $I$ is of general type,
all coherent sheaves are hyperholomorphic
(\ref{_generic_manifolds_Definition_},
\ref{_generic_are_dense_Proposition_}),
because all integer $(p,p)$-classes are 
$SU(2)$-invariant. However, for generic
complex structures $I$, the corresponding
complex manifold $(M, I)$ is never algebraic.
If we wish to compactify the moduli
\index{terms}{moduli space!of stable hyperholomorphic bundles} 
of holomorphic bundles, we need to consider algebraic
complex structures, and if we want to stay in hyperholomorphic
category, the complex structures must be generic. 
This paradox is reconciled by considering the $C$-restricted complex
structures (\ref{_C_restri_Definition_}).

Given a generic hyperk\"ahler structure $\c H$, consider
an algebraic complex \index{terms}{hyperk\"ahler structures!generic} 
structure $I$ with
$Pic(M, I) = \Z$. The group of rational $(p,p)$-cycles
has form 
\begin{equation} \label{_Pic_1_decomposi_Equation_}
\begin{split}
 H^{p,p}_I(M, \Q) = &H^{2p}(M ,\Q)^{SU(2)} \oplus
   a \cdot H^{2p}(M ,\Q)^{SU(2)}  \\
   & a^2 \cdot\oplus H^{2p}(M ,\Q)^{SU(2)} \oplus...
\end{split}
\end{equation}
where $a$ is a generator of $Pic(M, I)\subset H^{p,p}_I(M, \Z)$
and $H^{2p}(M ,\Q)^{SU(2)}$ is the group of rational
$SU(2)$-invariant cycles. This decomposition follows
from an explicit description of the algebra of cohomology
given by \ref{_S^*H^2_is_H^*M_intro-Theorem_}.
Let 
\[ \Pi:\;  H^{p,p}_I(M, \Q)\arrow a \cdot
     H^{2p}(M ,\Q)^{SU(2)}\oplus a^2 \cdot H^{2p}(M ,\Q)^{SU(2)}
   \oplus ...
\]
be the projection onto non-$SU(2)$-invariant part. 
Using Wirtinger's equality, we prove that a fundamental
class $[X]$ of a complex subvariety $X\subset (M, I)$ is $SU(2)$-invariant
unless $\deg \Pi([X])\neq 0$
(\ref{_Wirti_hyperka_Proposition_}). A similar result holds for
the second Chern class of a stable bundle 
(\ref{_stable_shea_degree_Corollary_},). 

A \index{terms}{$C$-restricted complex structures} $C$-restricted complex structure is, heuristically,
a structure for which the decomposition
\eqref{_Pic_1_decomposi_Equation_} folds, and
$\deg a>C$. For a $C$-restricted complex structure $I$,
and a fundamental class $[X]$ of a complex subvariety $X\subset (M, I)$
of complex codimension 2, we have $\deg [X]>C$
or $X$ is \index{terms}{trianalytic subvarieties} trianalytic. A version of Wirtinger's inequality
for vector bundles (\ref{_stable_shea_degree_Corollary_})
implies that a stable vector bundle $B$ over $(M, I)$ is
hyperholomorphic, unless $|\deg c_2(B)| >C$.
Therefore, over a $C$-restricted $(M, I)$,
all torsion-free \index{terms}{coherent sheaves!semistable} 
semistable coherent sheaves 
with bounded degree of the second
Chern class are hyperholomorphic
(\ref{_sheaf_on_C_restr_hyperho_Theorem_}). 

The utility of 
\index{terms}{$C$-restricted complex structures!utility of} 
$C$-restricted induced complex structures is that
they are algebraic, but behave like generic induced complex
structures with respect to the sheaves $F$ with low 
$|\deg c_2(F)|$ and $|\deg c_1(F)|$.

We prove that a generic hyperk\"ahler structure admits
$C$-restricted induced complex structures for all $C$,
and the set of $C$-restricted induced complex structures
is dense in the set of all induced complex structures
(\ref{_C_restri_dense_Theorem_}). 
We prove this by studying the algebro-geometric properties
of the 
\index{terms}{moduli space!of hyperk\"ahler structures} 
moduli of hyperk\"ahler structures on a given
hyperk\"ahler manifold (Subsection 
\ref{_modu_and_C-restri_Subsection_}).

\subsection{Quaternionic-K\"ahler manifolds and Swann's formalism}
\label{_Swann's_intro_Subsection_}

Quaternionic-K\"ahler manifolds (Subsection \ref{_B_2_bundles_Subsection_})
are a beautiful 
subject of Riemannian geometry. We are interested in these manifolds
because they are intimately connected with singularities of
\index{terms}{coherent sheaves!hyperholomorphic!singularities of} 
\index{terms}{coherent sheaves!hyperholomorphic!connections on} 
hyperholomorphic sheaves. A stable hyperholomorphic sheaf
is equipped with a natural connection, which is called 
{\bf \index{terms}{connections!hyperholomorphic!definition of} 
hyperholomorphic connection}. By definition, a 
{hyperholomorphic connection} on a torsion-free coherent sheaf
is a connection $\nabla$ defined outside of singularities of $F$,
with square-integrable curvature $\nabla^2$ which is an $SU(2)$-
invariant 2-form (\ref{_hyperholo_co_Definition_}). 
We have shown that a stable hyperholomorphic sheaf admits
a hyperholomorphic connection, and conversely, a 
\index{terms}{coherent sheaves!reflexive}
reflexive sheaf admitting a 
\index{terms}{connections!hyperholomorphic!existence of} 
hyperholomorphic connection is a direct
sum of stable hyperholomorphic sheaves
(\ref{_hyperho_conne_exi_Theorem_}).

Consider a 
\index{terms}{coherent sheaves!reflexive}
reflexive sheaf $F$ over $(M, I)$ 
with an isolated singularity in $x\in M$. Let $\nabla$ be 
a \index{terms}{connections!hyperholomorphic} hyperholomorphic 
connection on $F$. We prove that $F$
can be desingularized by a blow-up of its 
singular set. In other words, for $\pi:\; \tilde M \arrow (M, I)$
a blow-up of $x\in M$, the pull-back
$\pi^* F$ is a bundle over $\tilde M$.

Consider the restriction $\pi^* F\restrict C$ of 
$\pi^* F$ to the blow-up divisor \[ C = {\Bbb P} T_x M \cong\C P^{2n-1}. \]
To be able to deal with the singularities of $F$ effectively,
we need to prove that the bundle $\pi^* F\restrict C$
is a direct sum of stable bundles.

To study the singularities of coherent sheaves on hyperk\"ahler manifolds,
we consider vector bundles over \index{terms}{quaternionic-K\"ahler manifolds} quaternionic-K\"ahler manifolds,
studied by Berard Bergery \index{names}{Berard Bergery, L.}
and \index{names}{Salamon, S.} Salamon (\cite{_Salamon_}).
\index{terms}{quaternionic-K\"ahler manifolds!definition of} 
A quaternionic-K\"ahler manifold 
(\ref{_q-K_Definition_})
is a Riemannian manifold $Q$ equipped with a
bundle $W$ of algebras acting on its tangent bundle, and satisfying
the following conditions. The fibers of $W$ are 
(non-canonically) isomorphic 
to the quaternion algebra, the map $W\hookrightarrow \End(TQ)$
is compatible with the Levi-Civita connection, and the
unit quaternions $h\in W$ act as orthogonal automorphisms on $TQ$.
For each quaternionic-K\"ahler manifold $Q$,
one has a 
\index{terms}{twistor space!of a quaternionic-K\"ahler manifold} 
twistor space $\Tw(Q)$ (\ref{_twi_q-K_Definition_}), 
which is a total space of a spherical fibration consisting
of all $h\in W$ satisfying $h^2=-1$. 
Further on, we shall use the term ``quaternionic-K\"ahler''
\index{terms}{quaternionic-K\"ahler manifolds!with non-trivial scalar curvature} 
for manifolds with non-trivial $W$.

Consider the 
\index{terms}{twistor space!of a hyperk\"ahler manifold} 
twistor space $\Tw(M)$ of a hyperk\"ahler manifold
$M$, \\ equipped with a natural map \[ \sigma:\;\Tw(M) \arrow M.\]
Let $(B, \nabla)$ be a bundle over $M$ equipped with a
\index{terms}{connections!hyperholomorphic!twistor transform for} 
hyperholomorphic connection. A pullback 
$(\sigma^* B, \sigma^*\nabla)$ 
is a holomorphic bundle on $\Tw(M)$
(\ref{_autodua_(1,1)-on-twi_Lemma_}), that is,
the operator $\sigma^*\nabla^{0,1}$ is a holomorphic structure
operator on $\sigma^* B$. This correspondence
is called {\bf the direct twistor 
\index{terms}{twistor transform!direct} 
\index{terms}{twistor transform!inverse} 
transform}. It is
invertible: from a holomorphic bundle 
$(\sigma^* B, \sigma^*\nabla^{0,1})$ on $\Tw(M)$
it is possible to reconstruct $(B, \nabla)$, which 
is unique (\cite{_NHYM_}; see also \ref{_dire_inve_twisto_Theorem_}). 

A similar construction exists on 
\index{terms}{quaternionic-K\"ahler manifolds!twistor transform for} 
quaternionic-K\"ahler manifolds,
due to T. \index{names}{Nitta, T.} Nitta (\cite{_Nitta:bundles_}, \cite{_Nitta:Y-M_}). 
A bundle $(B, \nabla)$ on a quaternionic-K\"ahler manifold $Q$
is called {\bf a 
\index{terms}{$B_2$-bundles!definition of} 
$B_2$-bundle} if its curvature $\nabla^2$ is
invariant with respect to the adjoint action of ${\Bbb H}^*$
on $\Lambda^2(M, \End(B))$ (\ref{_B_2_bu_Definition_}). 
An analogue of direct and inverse 
\index{terms}{twistor transform!for $B_2$-bundles} 
transform exists
for $B_2$-bundles (\ref{_dire_inve_q-K_Theorem_}). 
Most importantly, T. \index{names}{Nitta, T.} Nitta proved that on a quaternionic-K\"ahler manifold
of positive scalar curvature a twistor transform
of a $B_2$-bundle is a 
\index{terms}{vector bundles!Yang--Mills} 
Yang-Mills bundle on
$\Tw(Q)$ (\ref{_twi_tra_YM_q-K_Theorem_}). This implies
that a twistor transform of a Hermitian $B_2$-bundle is a direct
sum of 
\index{terms}{vector bundles!stable!obtained as twistor transform of $B_2$-bundles} 
stable bundles with $\deg c_1 =0$.

In the situation described in the beginning of this Subsection,
we have a manifold 
$C= {\Bbb P} T_x M \cong \C P ^{2n-1}$ which is a twistor
space of a quaternionic \index{terms}{quaternionic projective space}
projective space 
\[ {\Bbb P}_{\Bbb H}T_x M = 
   \bigg(T_x M\backslash 0\bigg)/{\Bbb H}^* \cong {\Bbb H} P^n.
\] 
To prove that $\pi^* F\restrict C$ is a direct sum of stable 
bundles,
we need to show that $\pi^* F\restrict C$ is obtained
as a twistor transform of some Hermitian $B_2$-bundle 
\index{terms}{twistor transform!for $B_2$-bundles} 
on ${\Bbb P}_{\Bbb H}T_x M$. 

This is done using an equivalence between the category of 
qua\-ter\-ni\-onic-\-K\"ah\-ler manifolds of positive scalar curvature 
and the category of hyperk\"ahler manifolds
equipped with a special type of ${\Bbb H}^*$-action,
constructed by A. \index{names}{Swann, A.} Swann (\cite{_Swann_}). 
Given a 
\index{terms}{quaternionic-K\"ahler manifolds!hyperk\"ahler manifolds associated with} 
quaternionic-K\"ahler manifold $Q$, we consider
a principal bundle $\c U(Q)$ consisting of all quaternion
frames on $Q$ (\ref{_specia_and_q-K-Subsection_}). 
Then $\c U(Q)$ is fibered over $Q$ with a fiber
${\Bbb H}/\{\pm 1\}$. It is easy to show that 
$\c U(Q)$ is equipped with an action of quaternion algebra
in its tangent bundle. A. \index{names}{Swann, A.} Swann proved that if $Q$ has
with positive scalar curvature, then this action of quaternion
algebra comes from a hyperk\"ahler structure on $\c U(M)$
(\ref{_U(Q)_hyperk_Theorem_}). 

The correspondence $Q\arrow \c U(Q)$ induces an equivalence of
appropriately defined categories (\ref{_U(Q)_equiva_cate_Theorem_}).
We call this construction {\bf \index{names}{Swann, A.} \index{terms}{Swann's formalism} Swann's formalism}.

\index{terms}{twistor space!of a hyperk\"ahler manifold} 
The twistor space $\Tw(\c U(Q))$ of the hyperk\"ahler manifold 
$\c U(Q)$ is equ\-ipped with a holomorphic action of $\C^*$.
Every $B_2$-bundle corresponds to a $\C^*$-invariant
holomorphic bundle on $\Tw(\c U(Q))$ and this correspondence
induces an equivalence of appropriately defined categories,
called {\bf \index{names}{Swann, A.} \index{terms}{Swann's formalism!for vector bundles} Swann's formalism for budnles}
(\ref{_B_2_to_C^*_equiva_Theorem_}). Applying this equivalence
to the $\C^*$-equivariant sheaf obtained as an associate graded
sheaf of a hyperholomorphic sheaf, we obtain a $B_2$ bundle on
${\Bbb P}_{\Bbb H}T_x M$, and $\pi^* F\restrict C$ is obtained
from this $B_2$-bundle by a twistor transform. 
\index{terms}{twistor transform!for $B_2$-bundles} 

The correspondence between $B_2$-bundles on $Q$ and 
$\C^*$-invariant holomorphic bundles on $\Tw(\c U(Q))$
is a geometric phenomenon which is of independent 
interest. We construct it by reduction to $\dim Q=0$,
where it follows from an explicit calculation involving
2-forms over a flat manifold of real dimension 4.

\subsection{Contents}

The paper is organized as follows.

\begin{itemize} 

\item Section \ref{_intro_Section_} is an introduction.
It is independent from the rest of this apper.

\item Section \ref{_basics_hyperka_Section_} is an introduction 
to the theory of hyperk\"ahler manifolds. We give a compenduum of
results from hyperk\"ahler geometry
which are due to F. \index{names}{Bogomolov, F. A.} Bogomolov (\cite{_Bogomolov_}) and
A. \index{names}{Beauville, A.} Beauville (\cite{_Beauville_}), and give an introduction
to the results of \cite{_Verbitsky:Hyperholo_bundles_},
\cite{_Verbitsky:hypercomple_}, \cite{_Verbitsky:Symplectic_II_}.

\item Section \ref{_hyperho_shea_Section_} contains a definition
and basic properties of 
\index{terms}{coherent sheaves!hyperholomorphic!basic properties of} 
hyperholomorphic sheaves. We prove that
a stable hyperholomorphic sheaf admits a 
\index{terms}{connections!hyperholomorphic!on coherent sheaves} 
hyperholomorphic connection,
and conversely, a 
\index{terms}{coherent sheaves!reflexive}
reflexive sheaf admitting a
hyperholomorphic connection is stable hyperholomorphic
(\ref{_hyperho_conne_exi_Theorem_}). 
This equivalence is constructed using
Bando-Siu theory of 
\index{terms}{connections!Yang--Mills!in coherent sheaves} 
Yang--Mills connections
on coherent sheaves.

We prove an analogue of Wirtinger's inequality for stable sheaves
(\ref{_stable_shea_degree_Corollary_}), 
which states that for any
induced complex structure $J\neq \pm I$, and
any stable 
\index{terms}{coherent sheaves!Wirtinger's inequality for}
reflexive sheaf $F$ on
$(M,I)$, we have
\[ 
   \deg_I\left(2c_2(F) - \frac{r-1}{r} c_1(F)^2\right) \geq 
   \left|\deg_J\left(2c_2(F) - \frac{r-1}{r} c_1(F)^2\right)\right|,
\]
and the equality holds if and only if 
$F$ is hyperholomorphic.

\item Section \ref{_cohomo_hype_Section_} 
contains the preliminary material used for the study
of \index{terms}{$C$-restricted complex structures} $C$-restricted complex structures in Section \ref{_C_restri_Section_}.
We give an exposition of various algebraic structures on 
the cohomology of a hyperk\"ahler manifold, which
were discovered in \cite{_so(5)_} 
and \cite{_Verbitsky:cohomo_}. 
In the last Subsection, we
apply the Wirtinger's inequality to prove that
the fundamental classes of complex subvarieties
and $c_2$ of stable 
\index{terms}{coherent sheaves!reflexive}
reflexive sheaves satisfy
a certain set of axioms. Cohomology classes satisfying these
axioms are called \index{terms}{CA-type (cohomology classes of)} 
classes of CA-type. 
This definition simplifies
the work on $C$-restricted complex structures in 
Section \ref{_C_restri_Section_}.

\item In Section \ref{_C_restri_Section_} we define \index{terms}{$C$-restricted complex structures} $C$-restricted
complex structures and prove the following. Consider a compact
hyperk\"ahler manifold and an $SU(2)$-invariant class $a\in H^4(M)$.
Then for all $C$-restricted complex
structures $I$, with $C> \deg a$, and all 
\index{terms}{coherent sheaves!semistable}  
semistable sheaves
$I$ on $(M, I)$ with $c_2(F) =a$, the sheaf $F$
is hyperholomorphic
(\ref{_sheaf_on_C_restr_hyperho_Theorem_}). 
This is used to show that a 
deformation of a hyperholomorphic sheaf is again
hyperholomorphic, over $(M,I)$ with $I$ 
a $C$-restricted complex structure, 
$c> \deg c_2(F)$.

\index{terms}{moduli space!of hyperk\"ahler structures} 
We define the moduli space of hyperk\"ahler structures,
and show that for a dense set $\c C$ of hyperk\"ahler
\index{terms}{hyperk\"ahler structures!admitting $C$-restricted complex structures} 
structures, all $\c H \in \c C$ admit a dense set of
$C$-induced complex structures, for all $C\in \R$
(\ref{_C_restri_dense_Theorem_}).

\item In Section \ref{_desingu_she_Section_} we give a proof
of Desingularization Theorem for stable 
\index{terms}{coherent sheaves!reflexive}
reflexive 
\index{terms}{coherent sheaves!hyperholomorphic!with isolated singularities} 
hyperholomorphic sheaves with isolated singularities
(\ref{_desingu_hyperho_Theorem_}).
We study the natural \index{terms}{$\C^*$-action on hyperk\"ahler manifolds} $\C^*$-action on a local ring of 
a hyperk\"ahler manifold defined in \cite{_Verbitsky:DesinguII_}.
We show that a sheaf $F$ admitting a hyperholomorphic
connection is equivariant with respect to this
$\C^*$-action. Then $F$ can be desingularized by a blow-up,
because any 
\index{terms}{$\C^*$-equivariant sheaves!desingularization of} 
$\C^*$-equivariant sheaf with an isolated
singularity can be desingularized by a blow-up
(\ref{_desingu_C^*_equi_Proposition_}).

\item Section \ref{_twisto_tra_Section_} is a primer
on twistor transform and 
\index{terms}{quaternionic-K\"ahler manifolds}
\index{terms}{twistor transform}
 quaternionic-K\"ahler geometry. 
We give an exposition of the works of A. \index{names}{Swann, A.} Swann
(\cite{_Swann_}), T. \index{names}{Nitta, T.} Nitta 
(\cite{_Nitta:bundles_}, \cite{_Nitta:Y-M_}) on quaternionic-K\"ahler
manifolds and explain the direct and inverse twistor transform
over hyperk\"ahler and qua\-ter\-ni\-onic-\--K\"ah\-ler manifolds. 

\item Section \ref{_C_equiv_twi_spa_Section_}
gives a correspondence between $B_2$-bundles
\index{terms}{$B_2$-bundles!and $\C^*$-equivariant holomorphic bundles} 
on a qua\-ter\-ni\-onic-\--K\"ah\-ler manifold, and $\C^*$-equivariant
holomorphic bundles on the twistor space of the corresponding
\index{terms}{twistor space!of a hyperk\"ahler manifold} 
hyperk\"ahler manifold constructed by A. \index{names}{Swann, A.} Swann. This is
called ``\index{terms}{Swann's formalism!for vector bundles}Swann's formalism for vector bundles''.
We use this correspondence to prove that an associate
graded sheaf of a 
\index{terms}{coherent sheaves!hyperholomorphic!connections on} 
hyperholomorphic sheaf is equipped with
a natural connection which is compatible with quaternions.

\item In Section \ref{_modu_hyperho_Section_},
we use the equivariant geometry of the bundle $\pi^* F\restrict C$
to show that a hyperholomorphic deformation of
a \index{terms}{vector bundles!hyperholomorphic!deformations of} 
hyperholomorphic bundle is again a bundle.
Together with results on \index{terms}{$C$-restricted complex structures} $C$-restricted complex structures
and \index{names}{Maruyama, M.} Maruyama's compactification (\cite{_Maruyama:Si_}),
this implies that the moduli of semistable bundles
\index{terms}{vector bundles!semistable!moduli of}
\index{terms}{vector bundles!stable!moduli of} 
are compact, under conditions of $C$-restrictness and
non-existence of \index{terms}{trianalytic subvarieties} trianalytic subvarieties 
(\ref{_space_semista_bu_compa_Theorem_}). 

\item In Section \ref{_new_exa_Section_},
we apply these results to the hyperk\"ahler geometry.
Using the \index{terms}{hyperk\"ahler desingularization} desingularization theorem for singular hyperk\"ahler
manifolds (\ref{_desingu_Theorem_}), we prove that the moduli of
stable deformations of a 
\index{terms}{vector bundles!hyperholomorphic!moduli of} 
hyperholomorphic bundle has a compact
hyperk\"ahler \index{terms}{hyperk\"ahler desingularization} desingularization
(\ref{_space_sta_bu_compa_hyperka_Theorem_}).
 We give an exposition of the theory of twistor
paths, which allows one to identify the categories of 
\index{terms}{vector bundles!stable!category of} 
stable bundles for different K\"ahler structures on the same 
hyperk\"ahler manifold
(\ref{_admi_twi_impli_Theorem_}). 
These results allow one to weaken the conditions
necessary for compactness of the moduli spaces of
vector bundles. Finally, we give a conjectural exposition
of how these results can be used to obtain new examples
of compact hyperk\"ahler manifolds.

\end{itemize}


\section{Hyperk\"ahler manifolds}
\label{_basics_hyperka_Section_}


\subsection{Hyperk\"ahler manifolds}
\label{_hyperka_Subsection_}

This subsection contains a compression of 
the basic and best known results 
and definitions from hyperk\"ahler geometry, found, for instance, in
\cite{_Besse:Einst_Manifo_} or in \cite{_Beauville_}.

\hfill

\definition \label{_hyperkahler_manifold_Definition_} 
(\cite{_Besse:Einst_Manifo_}) A {\bf hyperk\"ahler manifold} is a
Riemannian manifold $M$ endowed with three complex structures $I$, $J$
and $K$, such that the following holds.
 
\begin{description}
\item[(i)]  the metric on $M$ is K\"ahler with respect to these complex 
\index{terms}{hyperk\"ahler manifold!definition of} structures and
 
\item[(ii)] $I$, $J$ and $K$, considered as  endomorphisms
of a real tangent bundle, satisfy the relation 
$I\circ J=-J\circ I = K$.
\end{description}

\hfill 

The notion of a hyperk\"ahler manifold was 
introduced by E. \index{names}{Calabi, E.} Calabi (\cite{_Calabi_}).

\hfill

Clearly, a hyperk\"ahler manifold has a natural action of
the quaternion algebra ${\Bbb H}$ in its real tangent bundle $TM$. 
Therefore its complex dimension is even.
For each quaternion $L\in \Bbb H$, $L^2=-1$,
the corresponding automorphism of $TM$ is an almost complex
structure. It is easy to check that this almost 
complex structure is integrable (\cite{_Besse:Einst_Manifo_}).

\hfill

\definition \label{_indu_comple_str_Definition_} 
Let $M$ be a hyperk\"ahler manifold, and $L$ a quaternion satisfying
$L^2=-1$. The corresponding complex structure 
\index{terms}{induced complex structures!definition of}
on $M$ is called
{\bf an induced complex structure}. The $M$, considered as a K\"ahler
manifold, is denoted by $(M, L)$. In this case,
the hyperk\"ahler 
\index{terms}{hyperk\"ahler structures!compatible with a complex structure} 
structure is called {\bf compatible
with the complex structure $L$}.

Let $M$ be a compact complex manifold. We say
that $M$ is {\bf of hyperk\"ahler type}
\index{terms}{complex manifold of hyperk\"ahler type}
if $M$ admits a hyperk\"ahler structure
compatible with the complex structure.

\hfill

\hfill

\definition \label{_holomorphi_symple_Definition_} 
Let $M$ be a complex manifold and $\Theta$ a closed 
holomorphic 2-form over $M$ such that 
$\Theta^n=\Theta\wedge\Theta\wedge...$, is
a nowhere degenerate section of a canonical class of $M$
($2n=dim_\C(M)$).
Then $M$ is called {\bf holomorphically 
\index{terms}{holomorphic symplectic structure} 
symplectic}.

\hfill

Let $M$ be a hyperk\"ahler manifold; denote the
Riemannian form on $M$ by $<\cdot,\cdot>$.
Let the form $\omega_I := <I(\cdot),\cdot>$ be the usual K\"ahler
form  which is closed and parallel
(with respect to the Levi-Civita connection). Analogously defined 
forms $\omega_J$ and $\omega_K$ are
also closed and parallel. 
 
A simple linear algebraic
consideration (\cite{_Besse:Einst_Manifo_}) shows that the form
$\Theta:=\omega_J+\sqrt{-1}\omega_K$ is of
type $(2,0)$ and, being closed, this form is also holomorphic.
Also, the form $\Theta$ is nowhere degenerate, as another linear 
algebraic argument shows.
It is called {\bf the canonical holomorphic \index{terms}{holomorphic symplectic structure} symplectic form
of a manifold M}. Thus, for each hyperk\"ahler manifold $M$,
and an induced complex structure $L$, the underlying complex manifold
$(M,L)$ is holomorphically symplectic. The converse assertion
is also true:

\hfill

\theorem \label{_symplectic_=>_hyperkahler_Proposition_}
(\cite{_Beauville_}, \cite{_Besse:Einst_Manifo_})
Let $M$ be a compact holomorphically
\index{terms}{holomorphic symplectic structure} symplectic K\"ahler manifold with the holomorphic symplectic form
$\Theta$, a K\"ahler class 
$[\omega]\in H^{1,1}(M)$ and a complex structure $I$. 
Let $n=\dim_\C M$. Assume that
$\int_M \omega^n = \int_M (Re \Theta)^n$.
Then there is a unique hyperk\"ahler 
\index{terms}{hyperk\"ahler structures} 
structure $(I,J,K,(\cdot,\cdot))$
over $M$ such that the cohomology class of the symplectic form
$\omega_I=(\cdot,I\cdot)$ is equal to $[\omega]$ and the
canonical symplectic form $\omega_J+\1\omega_K$ is
equal to $\Theta$.

\hfill

\ref{_symplectic_=>_hyperkahler_Proposition_} 
follows from the conjecture of \index{names}{Calabi, E.} \index{terms}{Calabi--Yau Theorem} Calabi, pro\-ven by
\index{names}{Yau, S.-T.} Yau (\cite{_Yau:Calabi-Yau_}). 
\endproof

\hfill

Let $M$ be a hyperk\"ahler manifold. We identify the group $SU(2)$
with the group of unitary quaternions. This gives a canonical 
action of $SU(2)$ on the tangent bundle, and all its tensor
powers. In particular, we obtain a natural action of $SU(2)$
on the bundle of differential forms. 

\hfill

\lemma \label{_SU(2)_commu_Laplace_Lemma_}
The action of $SU(2)$ on differential forms commutes
with the Laplacian.
 
{\bf Proof:} This is Proposition 1.1
of \cite{_Verbitsky:Symplectic_II_}. \endproof
 
Thus, for compact $M$, we may speak of the natural action of
$SU(2)$ in cohomology.

\hfill

Further in this article, we use the following statement.

\hfill

\lemma \label{_SU(2)_inva_type_p,p_Lemma_} 
Let $\omega$ be a differential form over
a hyperk\"ahler manifold $M$. The form $\omega$ is $SU(2)$-invariant
if and only if it is of Hodge type $(p,p)$ with respect to all 
induced complex structures on $M$.

{\bf Proof:} This is \cite{_Verbitsky:Hyperholo_bundles_}, 
Proposition 1.2. \endproof

\subsection{Simple hyperk\"ahler manifolds}

\definition \label{_simple_hyperkahler_mfolds_Definition_} 
(\cite{_Beauville_}) A connected simply connected 
compact hy\-per\-k\"ah\-ler manifold 
$M$ is called {\bf simple} if $M$ cannot be represented 
as a product of two hyperk\"ahler manifolds:
\[ 
   M\neq M_1\times M_2,\ \text{where} \ dim\; M_1>0 \ \ \text{and} 
   \ dim\; M_2>0
\] 
\index{names}{Bogomolov, F. A.} 
\index{terms}{Bogomolov's theorem!on decomposition}
Bogomolov proved that every compact hyperk\"ahler manifold has a finite
covering which is a product of a compact torus
and several simple hyperk\"ahler manifolds. 
Bogomolov's theorem implies the following result
(\cite{_Beauville_}):

\hfill

\theorem\label{_simple_mani_crite_Theorem_}
Let $M$ be a compact hyperk\"ahler manifold.
Then the following conditions are equivalent.
\begin{description}
\item[(i)] $M$ is simple
\item[(ii)] $M$ satisfies $H^1(M, \R) =0$, $H^{2,0}(M) =\C$,
where $H^{2,0}(M)$ is the space of $(2,0)$-classes taken with
respect to any of induced complex structures.
\end{description}

\endproof

\subsection{Trianalytic subvarieties in hyperk\"ahler
manifolds.}
 
In this subsection, we give a definition and basic properties
of \index{terms}{trianalytic subvarieties} trianalytic subvarieties of hyperk\"ahler manifolds. 
We follow \cite{_Verbitsky:Symplectic_II_}, 
\cite{_Verbitsky:DesinguII_}.
 
\hfill
 
Let $M$ be a compact hyperk\"ahler manifold, $\dim_\R M =2m$.
 
\hfill
 
\definition\label{_trianalytic_Definition_} 
Let $N\subset M$ be a closed subset of $M$. Then $N$ is
called {\bf \index{terms}{trianalytic subvarieties!definition of} 
trianalytic} if $N$ is a complex analytic subset 
of $(M,L)$ for any induced complex structure $L$.
 
\hfill

\hfill
 
Let $I$ be an induced complex structure on $M$,
and $N\subset(M,I)$ be
a closed analytic subvariety of $(M,I)$, $dim_\C N= n$.
Consider the homology class 
represented by $N$. Let $[N]\in H^{2m-2n}(M)$ denote 
the Poincare dual cohomology class, so called
{\bf fundamental class} of $N$. Recall that
the hyperk\"ahler structure induces the action of 
the group $SU(2)$ on the space $H^{2m-2n}(M)$.
 
\hfill
 
\theorem\label{_G_M_invariant_implies_trianalytic_Theorem_} 
Assume that $[N]\in  H^{2m-2n}(M)$ is invariant with respect
to the action of $SU(2)$ on $H^{2m-2n}(M)$. Then $N$ is 
\index{terms}{trianalytic subvarieties} trianalytic.
 
{\bf Proof:} This is Theorem 4.1 of 
\cite{_Verbitsky:Symplectic_II_}.
\endproof

\hfill

The following assertion is the key to the proof of
\ref{_G_M_invariant_implies_trianalytic_Theorem_}
(see \cite{_Verbitsky:Symplectic_II_} for details).

\hfill

\proposition \label{_Wirti_hyperka_Proposition_}
(Wirtinger's inequality) 
Let $M$ be a compact hyperk\"ahler manifold,
$I$ an induced complex structure and $X\subset (M, I)$
a closed complex subvariety for complex dimension $k$.
Let $J$ be an induced complex structure, $J \neq \pm I$,
and $\omega_I$, $\omega_J$ the associated K\"ahler forms.
Consider the numbers
\[ 
   \deg_I X:= \int_X \omega_I^k, \ \  \deg_J X:= \int_X \omega_J^k
\]
(these numbers are called {\bf degree of a subvariety $X$
with respect to $I$, $J$}.
\index{terms}{degree!of a subvariety}
Then $\deg_I X\geq |\deg_J X|$, and the inequality is strict 
unless $X$ is \index{terms}{trianalytic subvarieties} trianalytic.

\endproof

\hfill
 
\remark \label{_triana_dim_div_4_Remark_}
Trianalytic subvarieties have an action of quaternion algebra in
the tangent bundle. In particular,
the real dimension of such subvarieties is divisible by 4.

\hfill

\definition \label{_generic_manifolds_Definition_} 
Let $M$ be a complex manifold admitting a hyperk\"ahler
\index{terms}{hyperk\"ahler structures}
\index{terms}{induced complex structure of general type!definition of} 
structure $\c H$. We say that $M$ is {\bf of general type}
or {\bf generic} with respect to $\c H$ if all elements of the group
\[ \bigoplus\limits_p H^{p,p}(M)\cap H^{2p}(M,\Z)\subset H^*(M)\] 
are $SU(2)$-invariant. 
We say that $M$ is {\bf Mumford--Tate generic} 
if for all $n\in \Z^{>0}$, all the cohomology classes
\[ \alpha \in
   \bigoplus\limits_p H^{p,p}(M^n)\cap H^{2p}(M^n,\Z)\subset H^*(M^n)
\] 
are $SU(2)$-invariant. In other words,
$M$ is Mumford--Tate generic if for all
$n\in {\Bbb Z}^{>0}$, the $n$-th power $M^n$ is
generic. Clearly, Mumford--Tate generic
implies generic.

\hfill

\proposition \label{_generic_are_dense_Proposition_} 
Let $M$ be a compact manifold, $\c H$ a hyperk\"ahler
\index{terms}{hyperk\"ahler structures} structure on $M$ and $S$
be the set of induced complex structures over $M$. Denote by 
$S_0\subset S$ the set of $L\in S$ such that 
$(M,L)$ is Mumford-Tate generic with respect to $\c H$. 
Then $S_0$ is dense in $S$. Moreover, the complement
$S\backslash S_0$ is countable.

{\bf Proof:} This is Proposition 2.2 from
\cite{_Verbitsky:Symplectic_II_}
\endproof

\hfill

\ref{_G_M_invariant_implies_trianalytic_Theorem_} has the following
immediate corollary:

\corollary \label{_hyperkae_embeddings_Corollary_} 
Let $M$ be a compact holomorphically \index{terms}{holomorphic symplectic structure} symplectic 
manifold. Assume that $M$ is of general type with respect to
a hyperk\"ahler \index{terms}{hyperk\"ahler structures} structure $\c H$.
Let $S\subset M$ be closed complex analytic
subvariety. Then $S$ is \index{terms}{trianalytic subvarieties} trianalytic 
with respect to $\c H$.

\endproof

\hfill

In \cite{_Verbitsky:hypercomple_},
\cite{_Verbitsky:Desingu_}, \cite{_Verbitsky:DesinguII_},
we gave a number of equivalent definitions of a singular hyperk\"ahler
and hypercomplex variety. We refer the reader to 
\cite{_Verbitsky:DesinguII_} for the precise definition;
for our present purposes it suffices to say that all
\index{terms}{trianalytic subvarieties!desingularization of} 
trianalytic subvarieties are hyperk\"ahler varieties.
The following Desingularization Theorem is very 
useful in the study of trianalytic subvarieties.

\hfill

\theorem \label{_desingu_Theorem_}
(\cite{_Verbitsky:DesinguII_})
Let $M$ be a hyperk\"ahler variety, and
$I$ an induced complex structure.
Consider the normalization \[ \widetilde{(M, I)}\stackrel n\arrow (M,I)\] 
of $(M,I)$. Then $\widetilde{(M, I)}$ is smooth and
has a natural hyperk\"ahler
\index{terms}{hyperk\"ahler desingularization} 
structure $\c H$, such that the associated
map $n:\; \widetilde{(M, I)} \arrow (M,I)$ agrees with $\c H$.
Moreover, the hyperk\"ahler 
manifold $\tilde M:= \widetilde{(M, I)}$
is independent from the choice of induced complex structure $I$.

\endproof

\hfill

Let $M$ be a K3 surface, and $M^{[n]}$ be a Hilbert scheme of points on 
$M$. Then $M^{[n]}$ admits a hyperk\"ahler 
\index{terms}{hyperk\"ahler structures!on a Hilbert scheme of points} 
structure
(\cite{_Beauville_}).
In \cite{_Verbitsky:Hilbert_}, we proved the following
theorem.

\hfill

\theorem\label{_no_triana_subva_of_Hilb_Theorem_}
Let $M$ be a complex K3 surface without automorphisms. Assume that
$M$ is Mumford-Tate generic with respect to some hyperka\"hler structure. 
Consider the Hilbert scheme $M^{[n]}$ of points on $M$.
Pick a hyperk\"ahler structure on $M^{[n]}$ which is compatible with
the complex structure. Then $M^{[n]}$ has no proper
\index{terms}{trianalytic subvarieties} trianalytic subvarieties.

\endproof

\hfill


\subsection{Hyperholomorphic bundles}
\label{_hyperholo_Subsection_}


This subsection contains several versions of a
definition of hyperholomorphic connection in a complex
vector bundle over a hyperk\"ahler manifold.
We follow \cite{_Verbitsky:Hyperholo_bundles_}.

\hfill

 Let $B$ be a holomorphic vector bundle over a complex
manifold $M$, $\nabla$ a  connection 
in $B$ and $\Theta\in\Lambda^2\otimes End(B)$ be its curvature. 
This connection
\index{terms}{connections!compatible with a holomorphic structure}
\index{terms}{connections!integrable}
is called {\bf compatible with a holomorphic structure} if
$\nabla_X(\zeta)=0$ for any holomorphic section $\zeta$ and
any antiholomorphic tangent vector field $X\in T^{0,1}(M)$. 
If there exists a holomorphic structure compatible with the given
Hermitian connection then this connection is called
{\bf integrable}.
 
\hfill
 
\index{terms}{Hodge decomposition}
One can define a {\bf Hodge decomposition} in the space of differential
forms with coefficients in any complex bundle, in particular,
$End(B)$.

\hfill

\theorem \label{_Newle_Nie_for_bu_Theorem_}
Let $\nabla$ be a Hermitian connection in a complex vector
bundle $B$ over a complex manifold. Then $\nabla$ is integrable
if and only if $\Theta\in\Lambda^{1,1}(M, \End(B))$, where
$\Lambda^{1,1}(M, \End(B))$ denotes the forms of Hodge
type (1,1). Also,
the holomorphic structure compatible with $\nabla$ is unique.

{\bf Proof:} This is Proposition 4.17 of \cite{_Kobayashi_}, 
Chapter I.
$\blacksquare$

\hfill

This result has the following more general version:

\hfill

\proposition \label{_Newle_Nie_for_NH_bu_Proposition_}
Let $\nabla$ be an arbitrary (not necessarily Hermitian)
connection in a complex vector bundle $B$. Then
$\nabla$ is integrable
if and only its $(0,1)$-part has square zero.

$\blacksquare$

\hfill

This proposition is a version of Newlander-Nirenberg theorem.
For vector bundles, it was proven by \index{names}{Atiyah, M.} Atiyah and \index{names}{Bott, R.} Bott.

\hfill

\definition \label{_hyperho_conne_Definition_}
Let $B$ be a Hermitian vector bundle with
a connection $\nabla$ over a hyperk\"ahler manifold
$M$. Then $\nabla$ is called {\bf hyperholomorphic} if 
$\nabla$ is
integrable with respect to each of the complex structures induced
by the hyperk\"ahler \index{terms}{hyperk\"ahler structures} structure. 
 
As follows from 
\ref{_Newle_Nie_for_bu_Theorem_}, $\nabla$ is hyperholomorphic
if and only if its curvature $\Theta$ is of Hodge type (1,1) with
respect to any of complex structures induced by a hyperk\"ahler 
structure.

As follows from \ref{_SU(2)_inva_type_p,p_Lemma_}, 
$\nabla$ is hyperholomorphic
if and only if $\Theta$ is a $SU(2)$-invariant differential form.

\hfill

\example \label{_tangent_hyperholo_Example_} 
(Examples of \index{terms}{vector bundles!hyperholomorphic!examples of} 
hyperholomorphic bundles)

\begin{description}

\item[(i)]
Let $M$ be a hyperk\"ahler manifold, and $TM$ be its tangent bundle
equi\-p\-ped with the Levi--Civita connection $\nabla$. Consider a complex
structure on $TM$ induced from the quaternion action. Then $\nabla$
is a Hermitian connection
which is integrable with respect to each induced complex structure,
and hence, is hyperholomorphic.

\item[(ii)] For $B$ a hyperholomorphic bundle, all its tensor powers
are also 
\index{terms}{vector bundles!hyperholomorphic!tensor powers of} 
hyperholomorphic.

\item[(iii)] Thus, the bundles of differential forms on a hyperk\"ahler
manifold are also hyperholomorphic.

\end{description}


\subsection{Stable bundles and Yang--Mills connections.}
\label{_sta_bu_and_YM_Subsection_}


This subsection is a compendium of the most
basic results and definitions from the \index{terms}{connections!Yang--Mills} Yang--Mills theory
over K\"ahler manifolds, concluding in the fundamental
theorem of \index{terms}{Uhlenbeck--Yau Theorem} Uhlenbeck--\index{names}{Yau, S.-T.}Yau \cite{_Uhle_Yau_}.

\hfill

\definition\label{_degree,slope_destabilising_Definition_} 
Let $F$ be a coherent sheaf over
an $n$-dimensional compact K\"ahler manifold $M$. We define
{\bf the degree} $\deg(F)$ (sometimes the degree
is also denoted by $\deg c_1(F)$) as
\index{terms}{degree!of a sheaf}
\[ 
   \deg(F)=\int_M\frac{ c_1(F)\wedge\omega^{n-1}}{vol(M)}
\] 
\index{terms}{slope!of a coherent sheaf} 
and $\text{slope}(F)$ as
\[ 
   \text{slope}(F)=\frac{1}{\text{rank}(F)}\cdot \deg(F). 
\]
The number $\text{slope}(F)$ depends only on a
cohomology class of $c_1(F)$. 

Let $F$ be a coherent sheaf on $M$
and $F'\subset F$ its proper subsheaf. Then $F'$ is 
called {\bf destabilizing subsheaf} 
if $\text{slope}(F') \geq \text{slope}(F)$

A coherent sheaf $F$ is called {\bf 
\index{terms}{coherent sheaves!stable!definition of} stable}
\footnote{In the sense of Mumford-Takemoto}
if it has no destabilizing subsheaves. 
A coherent sheaf $F$ is called {\bf 
\index{terms}{coherent sheaves!semistable!definition of} 
semistable}
if for all destabilizing subsheaves $F'\subset F$,
we have $\text{slope}(F') = \text{slope}(F)$.
 
\hfill

Later on, we usually consider the bundles $B$ with $deg(B)=0$.

\hfill

Let $M$ be a K\"ahler manifold with a K\"ahler form $\omega$.
For differential forms with coefficients in any vector bundle
there is a Hodge operator $L: \eta\arrow\omega\wedge\eta$.
There is also a fiberwise-adjoint Hodge operator $\Lambda$
(see \cite{_Griffi_Harri_}).
 
\hfill

\definition \label{Yang-Mills_Definition_}
Let $B$ be a holomorphic bundle over a K\"ahler manifold $M$
with a holomorphic Hermitian connection $\nabla$ and a 
curvature $\Theta\in\Lambda^{1,1}\otimes End(B)$.
The Hermitian metric on $B$ and the connection $\nabla$
defined by this metric are called {\bf 
\index{terms}{vector bundles!Yang--Mills!definition} 
\index{terms}{connections!Yang--Mills!definition} 
\index{terms}{metric!Yang--Mills!definition}
Yang-Mills} if 

\[
   \Lambda(\Theta)=constant\cdot \Id\restrict{B},
\]
where $\Lambda$ is a Hodge operator and $\Id\restrict{B}$ is 
the identity endomorphism which is a section of $End(B)$.

Further on, we consider only these 
\index{terms}{connections!Yang--Mills} 
Yang--Mills connections for which this constant is zero.

\hfill

A holomorphic bundle is called  {\bf indecomposable} 
if it cannot be decomposed onto a direct sum
of two or more holomorphic bundles.

\hfill

The following fundamental 
theorem provides examples of 
\index{terms}{vector bundles!Yang--Mills!existence of} 
Yang-\--Mills \linebreak bundles.

\theorem \label{_UY_Theorem_} 
(\index{terms}{Uhlenbeck--Yau Theorem} Uhlenbeck-\index{names}{Yau, S.-T.}Yau)
Let B be an indecomposable
holomorphic bundle over a compact K\"ahler manifold. Then $B$ admits
a Hermitian \index{terms}{connections!Yang--Mills!and stability} 
Yang-Mills connection if and only if it is 
\index{terms}{vector bundles!stable!Yang--Mills connections on} 
stable, and this connection is unique.
 
{\bf Proof:} \cite{_Uhle_Yau_}. \endproof

\hfill

\proposition \label{_hyperholo_Yang--Mills_Proposition_}
Let $M$ be a hyperk\"ahler manifold, $L$
an induced complex structure and $B$ be a complex vector
bundle over $(M,L)$. 
Then every \index{terms}{connections!hyperholomorphic!are Yang-Mills} 
hyperholomorphic connection $\nabla$ in $B$
is \index{terms}{vector bundles!Yang--Mills} 
Yang-Mills and satisfies $\Lambda(\Theta)=0$
where $\Theta$ is a curvature of $\nabla$.
 
\hfill

{\bf Proof:} We use the definition of a hyperholomorphic 
connection as one with $SU(2)$-invariant curvature. 
Then \ref{_hyperholo_Yang--Mills_Proposition_}
follows from the

\hfill

\lemma \label{_Lambda_of_inva_forms_zero_Lemma_}
Let $\Theta\in \Lambda^2(M)$ be a $SU(2)$-invariant 
differential 2-form on $M$. Then
$\Lambda_L(\Theta)=0$ for each induced complex structure
$L$.\footnote{By $\Lambda_L$ we understand the Hodge operator 
$\Lambda$ associated with the K\"ahler complex structure $L$.}

{\bf Proof:} This is Lemma 2.1 of \cite{_Verbitsky:Hyperholo_bundles_}.
\endproof
 
\hfill

Let $M$ be a compact hyperk\"ahler manifold, $I$ an induced 
complex structure. 
For any \index{terms}{vector bundles!stable!Yang--Mills connections on} 
stable holomorphic bundle on $(M, I)$ there exists a unique
Hermitian \index{terms}{connections!Yang--Mills} Yang-Mills connection 
which, for some bundles, turns out to be hyperholomorphic. 
It is possible to tell when
this happens.

\hfill

\theorem \label{_inva_then_hyperho_Theorem_}
Let $B$ be a 
\index{terms}{vector bundles!stable!hyperholomorphic connections on} 
stable holomorphic bundle over
$(M,I)$, where $M$ is a hyperk\"ahler manifold and $I$
is an induced complex structure over $M$. Then 
$B$ admits a compatible 
\index{terms}{connections!hyperholomorphic!existence of} 
hyperholomorphic connection if and only
if the first two Chern classes $c_1(B)$ and $c_2(B)$ are 
$SU(2)$-invariant.\footnote{We use \ref{_SU(2)_commu_Laplace_Lemma_}
to speak of action of $SU(2)$ in cohomology of $M$.}

{\bf Proof:} This is Theorem 2.5 of
 \cite{_Verbitsky:Hyperholo_bundles_}. \endproof

\subsection{Twistor spaces}

Let $M$ be a hyperk\"ahler manifold. Consider the product manifold $X = M
\times S^2$. Embed the sphere $S^2 \subset {\Bbb H}$ 
into the quaternion algebra
${\Bbb H}$ 
as the subset of all quaternions $J$ with $J^2 = -1$. For every point
$x = m \times J \in X = M \times S^2$ the tangent space $T_xX$ is
canonically decomposed $T_xX = T_mM \oplus T_JS^2$. Identify $S^2 = \C P^1$
and let $I_J:T_JS^2 \to T_JS^2$ be the complex structure operator. Let
$I_m:T_mM \to T_mM$ be the complex structure on $M$ induced by $J \in S^2
\subset {\Bbb H}$.

The operator $I_x = I_m \oplus I_J:T_xX \to T_xX$ satisfies $I_x \circ I_x =
-1$. It depends smoothly on the point $x$, hence defines an almost complex
structure on $X$. This almost complex structure is known to be integrable
(see \cite{_Salamon_}). 

\hfill

\definition\label{_twistor_Definition_}
The complex manifold $\langle X, I_x \rangle$ is called {\it the twistor
space} for the hyperk\"ahler manifold $M$, denoted by $\Tw(M)$.
\index{terms}{$\Tw(M)$}
This manifold is equipped with a real analytic projection
$\sigma:\; \Tw(M)\arrow M$ and a complex analytic
projection $\pi:\; \Tw(M) \arrow \C P^1$.
\index{terms}{twistor space!of a hyperk\"ahler manifold!definition of} 

\hfill

The twistor space $\Tw(M)$ is not, generally speaking,
a K\"ahler manifold. For $M$ compact,
it is easy to show that $\Tw(M)$ does not admit
a K\"ahler metric. We consider $\Tw(M)$
as a Hermitian manifold with the product metric.

\hfill

\definition \label{_Li_Yau_condi_Definition_}
Let $X$ be an $n$-dimensional Hermitian manifold and let
$\sqrt{-1}\omega$ be the imaginary part of the metric on $X$. Thus $\omega$
is a real $(1,1)$-form.
Assume that the form $\omega$ satisfies the following condition
of \index{names}{Li, J.} Li and \index{names}{Yau, S.-T.}Yau (\cite{_Li_Yau_}). 
\begin{equation}\label{_Li_Yau_condi_Equation_}
\omega^{n-2} \wedge d\omega = 0.
\end{equation}
Such Hermitian metrics are called {\bf metrics
satisfying the condition of 
\index{terms}{Li--Yau condition}
\index{names}{Yau, S.-T.}
Li--Yau}.

For a closed real $2$-form $\eta$ let 
$$
\deg\eta = \int_X \omega^{n-1} \wedge \eta.
$$
\index{terms}{degree!of a form}
The condition \eqref{_Li_Yau_condi_Equation_} 
ensures that $\deg\eta$ depends only on the
cohomology class of $\eta$.  Thus it defines a degree functional
$\deg:H^2(X,\R) \to \R$. This functional allows one to repeat verbatim the
Mumford-Takemoto definitions of 
\index{terms}{vector bundles!stable!on Li--Yau manifolds} 
stable and semistable bundles in this
more general situation. Moreover, the Hermitian 
\index{terms}{connections!Yang--Mills!on Li--Yau manifolds} 
Yang-Mills equations also
carry over word-by-word. \index{names}{Li, J.} Li and \index{names}{Yau, S.-T.}Yau proved a version of 
\index{terms}{Uhlenbeck--Yau Theorem} Uhlenbeck--Yau theorem in this situation 
(\cite{_Li_Yau_}; see also \ref{_UY_for_shea_Theorem_}).

\hfill

\proposition
Let $M$ be a hyperk\"ahler manifold and
$\Tw(M)$ its twistor space, considered 
as a Hermitian manifold. Then $\Tw(M)$
satisfies the conditions of \index{names}{Li, J.} \index{terms}{Li--Yau condition}
Li--\index{names}{Yau, S.-T.}Yau.

{\bf Proof:} \cite{_NHYM_}, Proposition 4.5.
\endproof


\section{Hyperholomorphic sheaves}
\label{_hyperho_shea_Section_}


\subsection{Stable sheaves and Yang-Mills connections}

In \cite{_Bando_Siu_}, S. \index{names}{Bando, S.} Bando and Y.-T. \index{names}{Siu, Y.-T.} Siu developed the
machinery allowing one to apply the methods of Yang-Mills
theory to torsion-free coherent sheaves. In the course of 
this paper, we apply their work to generalise the results of
\cite{_Verbitsky:Hyperholo_bundles_}. In this Subsection,
we give a short exposition of their results.

\hfill

\definition\label{_refle_Definition_}
Let $X$ be a complex manifold, and $F$ a coherent sheaf on $X$.
Consider the sheaf $F^*:= \c Hom_{\calo_X}(F, \calo_X)$.
There is a natural functorial map 
$\rho_F:\; F \arrow F^{**}$. The sheaf $F^{**}$
is called {\bf a reflexive hull}, or {\bf \index{terms}{reflexization}
reflexization}
of $F$. The sheaf $F$ is called {\bf reflexive} if the map
$\rho_F:\; F \arrow F^{**}$ is an isomorphism. 

\hfill

\remark
For all coherent sheaves $F$, the map
$\rho_{F^*}:\; F^* \arrow F^{***}$ is an isomorphism
(\cite{_OSS_}, Ch. II, the proof of Lemma 1.1.12).
Therefore, a 
\index{terms}{coherent sheaves!reflexive}
reflexive hull of a sheaf is always 
reflexive.

\hfill

\claim 
Let $X$ be a K\"ahler manifold, and $F$ a torsion-free coherent sheaf over 
$X$. Then $F$ (semi)stable if and only if $F^{**}$ 
is (semi)stable.

{\bf Proof:} This is
\cite{_OSS_}, Ch. II, Lemma 1.2.4.
\endproof

\hfill

\definition 
Let $X$ be a K\"ahler manifold, and $F$ a coherent sheaf over 
$X$. The sheaf $F$ is called
{\bf polystable} if $F$ is a direct sum of 
\index{terms}{coherent sheaves!polystable!definition of} 
stable sheaves.

\hfill

The admissible 
Hermitian metrics, introduced by \index{names}{Bando, S.} Bando and \index{names}{Siu, Y.-T.} Siu
in \cite{_Bando_Siu_}, play the role of the
ordinary Hermitian metrics for vector bundles.

\hfill

Let $X$ be a K\"ahler manifold.
In Hodge theory, one considers the operator 
$\Lambda:\; \Lambda^{p, q}(X) \arrow\Lambda^{p-1, q-1}(X)$
acting on differential forms on $X$, which is adjoint to the
multiplication by the K\"ahler form. This operator is defined
on differential forms with coefficient in every bundle.
Considering a curvature $\Theta$ of a bundle $B$
as a 2-form with coefficients in $\End(B)$, we define
the expression $\Lambda\Theta$ which is a section of
$\End(B)$.

\hfill

\definition \label{_admi_metri_Definition_}
Let $X$ be a K\"ahler manifold, and $F$ a 
\index{terms}{coherent sheaves!reflexive}
reflexive 
coherent sheaf over $X$. Let $U\subset X$ be the set of all
points at which $F$ is locally trivial. By definition,
the restriction $F\restrict U$ of $F$ to $U$ is a bundle.
An {\bf admissible metric} on $F$ is a Hermitian metric $h$
on the bundle $F\restrict U$ which satisfies the following
assumptions
\begin{description}
\item[(i)] the curvature $\Theta$ of $(F, h)$ is square integrable, and 
\item[(ii)] the corresponding section $\Lambda \Theta\in \End(F\restrict U)$ 
is uniformly bounded. 
\end{description}

\hfill

\definition \label{_Yang-Mills_sheaves_Definition_}
Let $X$ be a K\"ahler manifold, $F$ a 
\index{terms}{coherent sheaves!reflexive}
reflexive 
coherent sheaf over $X$, and $h$ an admissible metric on $F$.
Consider the corresponding Hermitian connection
$\nabla$ on $F\restrict U$. The metric $h$ and
the connection $\nabla$ are called {\bf 
\index{terms}{connections!Yang--Mills!admissible} Yang-Mills}
if its curvature satisfies
\[ \Lambda \Theta\in \End(F\restrict U) = c\cdot \id 
\]
where $c$ is a constant and $\id$ the unit section 
$\id \in \End(F\restrict U)$. 

\hfill

Further in this
paper, we shall only consider 
\index{terms}{connections!Yang--Mills!with $\Lambda\nabla^2=0$} 
Yang-Mills connections with $\Lambda \Theta=0$. 

\hfill

\remark
By Gauss-Bonnet formule, the 
constant $c$ is equal to $\deg(F)$, where $\deg(F)$
is the degree of $F$ (\ref{_degree,slope_destabilising_Definition_}).
\index{terms}{degree!of a sheaf}

\hfill

One of the main results of \cite{_Bando_Siu_}
is the following analogue of \index{terms}{Uhlenbeck--Yau Theorem} Uhlenbeck--\index{names}{Yau, S.-T.}Yau theorem
(\ref{_UY_Theorem_}).

\hfill

\theorem\label{_UY_for_shea_Theorem_}
Let $M$ be a compact K\"ahler manifold, or a compact
Hermitian manifold satisfying conditions of \index{names}{Li, J.}
 \index{terms}{Li--Yau condition}
Li-\index{names}{Yau, S.-T.}Yau 
(\ref{_Li_Yau_condi_Definition_}), and $F$ a coherent
sheaf without torsion. Then $F$ admits an admissible 
\index{terms}{metric!Yang--Mills!admissible} Yang--Mills 
metric is and only if $F$ is polystable. Moreover, if $F$
is \index{terms}{coherent sheaves!stable!admissible Yang--Mills metric on} 
stable, then this metric is unique, up to a constant
multiplier. 

{\bf Proof:} In \cite{_Bando_Siu_}, \ref{_UY_for_shea_Theorem_} 
is proved for K\"ahler $M$ (\cite{_Bando_Siu_}, Theorem 3). 
It is easy to adapt this proof for 
Hermitian manifolds satisfying conditions of Li--\index{names}{Yau, S.-T.}Yau.
\endproof

\hfill

\remark
Clearly, the connection, corresponding to a metric on $F$,
does not change when the metric is multiplied by a scalar.
The \index{terms}{metric!Yang--Mills!on a polystable sheaf} 
Yang--Mills metric on a polystable sheaf is unique up to
a componentwise multiplication by scalar multipliers.
Thus, the Yang--Mills connection of \ref{_UY_for_shea_Theorem_}
is unique.

\hfill

Another important theorem of \cite{_Bando_Siu_} is the following.

\hfill

\theorem\label{_YM_can_be_exte_Theorem_}
Let $(F, h)$ be a holomorphic vector bundle with a
Hermitian metric $h$ defined on a K\"ahler manifold
$X$ (not necessary compact nor complete) outside a
closed subset $S$ with locally finite Hausdorff measure of
real co-dimension $4$.    
Assume that the curvature tensor
of $F$ is locally square integrable on $X$. Then
$F$ extends to the
whole space $X$ as a 
\index{terms}{coherent sheaves!reflexive}
reflexive sheaf $\c F$. 
Moreover, if the metric $h$ is
\index{terms}{metric!Yang--Mills} 
Yang-Mills, then $h$ can be smoothly extended 
as a Yang-Mills metric over the place where ${\cal F}$ is
locally free. 

{\bf Proof:} This is \cite{_Bando_Siu_}, Theorem 2.
\endproof

\subsection{Stable and semistable sheaves over
hyperk\"ahler manifolds}

Let $M$ be a compact hyperk\"ahler manifold, $I$ an induced
complex structure, $F$ a torsion-free coherent sheaf
over $(M,I)$ and $F^{**}$ its \index{terms}{reflexization}
reflexization. Recall that the cohomology of
$M$ are equipped with a natural $SU(2)$-action
(\ref{_SU(2)_commu_Laplace_Lemma_}). The motivation for
the following definition is \ref{_inva_then_hyperho_Theorem_} 
and \ref{_UY_for_shea_Theorem_}.

\hfill

\definition \label{_hyperho_shea_Definition_}
Assume that the first two Chern classes
of the sheaves $F$, $F^{**}$ are $SU(2)$-invariant.
Then $F$ is called {\bf stable 
\index{terms}{coherent sheaves!hyperholomorphic!definition of} 
hyperholomorphic} if
$F$ is stable, and {\bf semistable hyperholomorphic}
if $F$ can be obtained as a successive 
extension of stable hyperholomorphic
sheaves.

\hfill

\remark \label{_slope_hyperho_Remark_}
\index{terms}{slope!of a hyperholomorphic sheaf} 
The slope of a hyperholomorphic sheaf is zero, because
a degree of an $SU(2)$-invariant 2-form
is zero (\ref{_Lambda_of_inva_forms_zero_Lemma_}).

\hfill

\claim \label{_hyperho_suppo_of_F^**/F_Claim_}
Let $F$ be a \index{terms}{coherent sheaves!semistable} semistable coherent sheaf over $(M,I)$. Then the following 
conditions are equivalent.
\begin{description}
\item[(i)] $F$ is stable 
\index{terms}{coherent sheaves!hyperholomorphic!singularities of} 
hyperholomorphic
\item[(ii)] Consider the support $S$ of the sheaf $F^{**}/F$
as a complex subvariety of $(M,I)$. Let $X_1$, ... , $X_n$ be 
the set of irreducible components of $S$ of codimension 2. 
Then $X_i$ is \index{terms}{trianalytic subvarieties} trianalytic for all $i$, 
and the sheaf $F^{**}$ is stable hyperholomorphic.
\end{description}

{\bf Proof:} Consider an exact sequence
\[ 0 \arrow F \arrow F^{**} \arrow F^{**}/ F\arrow 0. \]
Let $[F / F^{**}]\in H^4(M)$ be the fundamental class
of the union of all co\-di\-men\-sion-2 components of support of
the sheaf $F / F^{**}$, taken with appropriate multiplicities.
Then, $c_2(F^{**}/ F) =- [F / F^{**}]$.
From the product formula for Chern classes,
it follows that 
\begin{equation} \label{_c_2(F)_and_F^**_Equation_}
   c_2(F)= c_2(F^{**}_i) + c_2(F^{**}/ F) 
   = c_2(F^{**}_i) - [F / F^{**}].
\end{equation}
Clearly, if all $X_i$ are \index{terms}{trianalytic subvarieties} trianalytic then
the class $[F / F^{**}]$ is $SU(2)$-invariant.
Thus, if the sheaf $F^{**}$ is hyperholomorphic
and all $X_i$ are trianalytic, then the second
Chern class of $F$ is $SU(2)$-invariant, and $F$
is hyperholomorphic. Conversely, assume that
$F$ is hyperholomorphic. We need to show that
all $X_i$ are trianalytic.
By definition,
\[ [F / F^{**}] = \sum_i \lambda_i [X_i]
\]
where $[X_i]$ denotes the fundamental class of $X_i$,
and $\lambda_i$ is the multiplicity of $F / F^{**}$
at $X_i$. By \eqref{_c_2(F)_and_F^**_Equation_},
($F$ hyperholomorphic) implies that the class
$[F / F^{**}]$ is $SU(2)$-invariant. Since
$[F / F^{**}]$ is $SU(2)$-invariant, we have
\[ \sum_i \lambda_i\deg_J(X_i) = \sum_i \lambda_i\deg_I(X_i).
\]
By Wirtinger's inequality
(\ref{_Wirti_hyperka_Proposition_}),
\[ \deg_J(X_i) \leq\deg_I(X_i),
\]
and the equality is reached only if $X_i$ is \index{terms}{trianalytic subvarieties} trianalytic.
By definition, all the numbers $\lambda_i$ are positive.
Therefore,
\[ \sum_i \lambda_i\deg_J(X_i) \leq \sum_i \lambda_i\deg_I(X_i).
\]
and the equality is reached only if all
the subvarieties $X_i$ are trianalytic.
This finishes the proof of 
\ref{_hyperho_suppo_of_F^**/F_Claim_}.
\endproof

\hfill

\claim
Let $M$ be a compact hyperk\"ahler
manifold, and $I$ an induced complex
structure of general type. Then 
all torsion-free coherent sheaves 
over $(M, I)$ are semistable hyperholomorphic.

{\bf Proof:} Let $F$ be a torsion-free coherent sheaf
over $(M, I)$. Clearly from the definition
of induced complex
structure of general type, the sheaves
$F$ and $F^{**}$ have $SU(2)$-invariant 
Chern classes. Now, all $SU(2)$-invariant
2-forms have degree zero (\ref{_Lambda_of_inva_forms_zero_Lemma_}), 
\index{terms}{degree!of a form}
and thus, $F$ is semistable.
\endproof

\subsection{Hyperholomorphic connection in 
torsion-free sheaves}

Let $M$ be a hyperk\"ahler manifold, $I$ an induced complex
structure, and $F$ a torsion-free sheaf over $(M,I)$.
Consider the natural $SU(2)$-action in the bundle
$\Lambda^i (M,B)$ of the differential $i$-forms
with coefficients in a vector bundle $B$. Let
$\Lambda^i_{inv}(M, B)\subset \Lambda^i (M, B)$
be the bundle of $SU(2)$-invariant $i$-forms.

\hfill

\definition \label{_hyperholo_co_Definition_}
Let $X\subset (M, I)$ be a complex subvariety of
codimension at least 2, such that $F\restrict{M\backslash X}$
is a bundle, $h$ be an admissible metric on 
$F\restrict{M\backslash X}$ and $\nabla$ the associated
connection. Then $\nabla$ is called {\bf hyperholomorphic} if its
curvature 
\[ \Theta_\nabla = 
   \nabla^2 \in 
   \Lambda^2\left(M, \End\left(F\restrict{M\backslash X}\right)\right)
\]
is $SU(2)$-invariant, i. e. belongs to 
$\Lambda^2_{inv}\left(M, 
\End\left(F\restrict{M\backslash X}\right)\right)$.

\hfill

\claim\label{_singu_triana_Claim_}
The singularities of a 
\index{terms}{connections!hyperholomorphic!singular set of} 
hyperholomorphic connection form a 
\index{terms}{trianalytic subvarieties} trianalytic
subvariety in $M$.

{\bf Proof:} Let $J$ be an induced complex structure on $M$,
and $U$ the set of all points of $(M,I)$ where $F$ is non-singular.
Clearly, $(F, \nabla)$ is a bundle with admissible connection on
$(U,J)$. Therefore, the holomorphic structure on $F\restrict{(U,J)}$ 
can be extended to $(M,J)$. Thus, the singular set of $F$ is 
holomorphic with respect to $J$. This proves 
\ref{_singu_triana_Claim_}.
\endproof

\hfill

\proposition \label{_conne_=>_hyperho_Proposition_}
Let $M$ be a compact hyperk\"ahler manifold, $I$ an
induced complex structure and  $F$ 
a 
\index{terms}{coherent sheaves!reflexive}
reflexive sheaf admitting a 
\index{terms}{coherent sheaves!hyperholomorphic!connections on} 
hyperholomorphic
connection. Then $F$ is a polystable hyperholomorphic sheaf.

\hfill

{\bf Proof:} By \ref{_hyperho_co_YM_Remark_} and 
\ref{_UY_for_shea_Theorem_}, $F$ is polystable. We 
need only to show that the Chern classes $c_1(F)$ and $c_2(F)$
are $SU(2)$-invariant. Let $U\subset M$ be the maximal
open subset of $M$ such that $F\restrict U$ is locally
trivial. By \ref{_YM_can_be_exte_Theorem_}, the metric $h$
and the connection $\nabla$ can be extended to $U$. 
Let $\Tw U\subset \Tw M$ be the corresponding twistor space,
and $\sigma:\; \Tw U \arrow U$ the standard map.
Consider the bundle $\sigma^* F\restrict U$, equipped with
a connection $\sigma^* \nabla$. It is well known \footnote{See 
for instance the section 
``Direct and inverse twistor \index{terms}{twistor transform} transform'' in
\cite{_NHYM_}.}
that $\sigma^* F\restrict U$ is a bundle
with an admissible Yang-Mills metric (we use Yang-Mills in the
sense of 
\index{terms}{metric!Yang--Mills!in the sense of Li--Yau} 
\index{terms}{Li--Yau!Yang--Mills metric, in the sense of}
\index{names}{Yau, S.-T.}
Li-Yau, see \ref{_Li_Yau_condi_Definition_}). By 
\ref{_YM_can_be_exte_Theorem_}, $\sigma^* F\restrict U$
can be extended to a 
\index{terms}{coherent sheaves!reflexive}
reflexive sheaf $\c F$ on $\Tw M$.
Clearly, this extension coincides with the push-forward of 
$\sigma^* F\restrict U$.
The singular set $\tilde S$ of $\c F$ is a pull-back of the singular set
$S$ of $F$. Thus, $S$ is \index{terms}{trianalytic subvarieties} trianalytic. By \index{terms}{hyperk\"ahler desingularization} desingularization theorem
(\ref{_desingu_Theorem_}), 
$S$ can be desingularized to a hyperk\"ahler manifold
in such a way that its twistors form a desingularization of $\c S$.
{}From the exact description of the singularities of $\c S$, provided by the
desingularization theorem, we obtain that the standard projection
$\pi:\; \c S \arrow \C P^1$ is flat.
By the following lemma,  
the restriction of $\c F$ to the fiber $(M,I) = \pi^{-1}(\{I\})$
of $\pi$ coincides with $F$.

\hfill

\lemma \label{_exte_flat_Lemma_}
Let $\pi:\; X \arrow Y$ be a map of complex varieties, and
$S\hookrightarrow X$ a subvariety of $X$ of codimension
at least 2, which is flat over $Y$.
Denote by $U\stackrel j \hookrightarrow X$ 
the complement $U = (X\backslash S)$. Let $F$ be a vector bundle over
$U$, and $j_* F$ its push-forward. Then 
the restriction of
$j_* F$ to the fibers of $\pi$ is 
\index{terms}{coherent sheaves!reflexive}
reflexive.

\hfill

{\bf Proof:} Let $Z= \pi^{-1}(\{y\})$ be a fiber of $\pi$. 
Since $S$ is flat over $Y$ and of codimension at least 2, we have
$j_*(\calo_{Z \cap U}) = \calo_{Z}$.
Clearly, for an open embedding $\gamma:\; T_1 \arrow T_2$ and coherent sheaves
$A, B$ on $T_1$, we have $\gamma_*(A\otimes B) = \gamma_* A \otimes \gamma_* B$.
Thus, for all coherent sheaves $A$ on $U$, we have
\begin{equation} \label{_j_*_commu_tenso_Equation_} 
  j_* A \otimes \calo_{Z} = j_*(A \otimes \calo_{Z\cap U}). 
\end{equation}
This implies that $j_*(F \restrict Z) = j_* F\restrict Z$.
It is well known (\cite{_OSS_}, Ch. II, 1.1.12; see also
\ref{_normal_refle_Lemma_}) that a push-forward of a
\index{terms}{coherent sheaves!reflexive!singularities of}
reflexive sheaf under an open embedding $\gamma$ is reflexive, provided that
the complement of the image of $\gamma$ has codimension at least 2.
Therefore, $j_* F\restrict Z$ is a reflexive sheaf over $Z$.
This proves \ref{_exte_flat_Lemma_}. \endproof 

\hfill

Return to the proof of \ref{_conne_=>_hyperho_Proposition_}.
Consider
the sheaf $\c F$ on the twistor space constructed above.
Since $\c F$ is reflexive, its singularities have 
codimension at least 3
(\cite{_OSS_}, Ch. II, 1.1.10). 
Therefore, $\c F$ is flat in codimension 2, and
the first two Chern classes of 
$F= \c F\restrict{\pi^{-1}(I)}$ can be obtained by restricting
the first two Chern classes of $\c F$ to the subvariety
$(M,I) = \pi^{-1}(I) \subset \Tw(M)$. It remains to show
that such restriction is $SU(2)$-invariant. 
Clearly, $H^2((M, I)) = H^2((M, I)\backslash S)$,
and $H^4((M, I)) = H^4((M, I)\backslash S)$. 
Therefore, 
\[ c_1\left(\c F \restrict {(M,I)}\right) 
     = c_1\left(\c F \restrict {(M,I)\backslash S}\right) 
\]
and
\[ c_2\left(\c F \restrict {(M,I)}\right) = 
c_2\left(\c F \restrict {(M,I)\backslash S}\right).
\]
On the other hand, the restriction
$\c F \restrict {\Tw(M)\backslash S}$ is a bundle. 
Therefore, the classes
\[
 c_1\left(\c F \restrict {(M,I)\backslash S}\right), \ \ 
 c_2\left(\c F \restrict {(M,I)\backslash S}\right)
\]
are independent from $I\in \C P^1$. 
On the other hand, these classes are of type $(p,p)$
with respect to all induced complex structures
$I\in \C P^1$. By \ref{_SU(2)_inva_type_p,p_Lemma_},
this implies that the classes $c_1(\c F \restrict {(M,I)})$,
$c_1(\c F \restrict {(M,I)})$ are $SU(2)$-invariant.
As we have shown above, these two classes 
are equal to the first Chern classes of $F$.
\ref{_conne_=>_hyperho_Proposition_}
is proven.
\endproof

\hfill

\subsection{Existence of hyperholomorphic connections}
\index{terms}{connections!hyperholomorphic!existence of} 

The following theorem is the main result of this section.

\hfill

\theorem\label{_hyperho_conne_exi_Theorem_}
Let $M$ be a compact hyperk\"ahler manifold, $I$
an induced complex structure and $F$ a 
\index{terms}{coherent sheaves!reflexive!hyperholomorphic connections on}
\index{terms}{coherent sheaves!hyperholomorphic!connections on} 
reflexive sheaf on
$(M,I)$. Then $F$ admits a hyperholomorphic connection if
and only if $F$ is polystable hyperholomorphic
\index{terms}{coherent sheaves!hyperholomorphic!polystable} 
in the sense of \ref{_hyperho_shea_Definition_}.

\hfill

\remark \label{_hyperho_co_YM_Remark_}
{}From \ref{_Lambda_of_inva_forms_zero_Lemma_}, it is clear that
a hyperholomorphic connection is always 
\index{terms}{connections!Yang--Mills!hyperholomorphic} 
Yang-Mills. Therefore,
such a connection is unique (\ref{_UY_for_shea_Theorem_}). 

\hfill

The ``only if'' part of \ref{_hyperho_conne_exi_Theorem_} 
 is \ref{_conne_=>_hyperho_Proposition_}.
The proof of ``if'' part of 
\ref{_hyperho_conne_exi_Theorem_} takes the rest of
this subsection. 

\hfill

Let $I$ be an induced complex structure.
We denote the  corresponding Hodge decomposition on differential forms
by $\Lambda^*(M)= \oplus \Lambda^{p,q}_I(M)$, and
the standard Hodge operator by 
$\Lambda_I:\;  \Lambda^{p,q}_I(M) \arrow  \Lambda^{p-1,q-1}_I(M)$.
All these structures are defined on the differential forms
with coefficients in a bundle.
\index{terms}{degree!of a form}
Let $\deg_I \eta:= \int_M Tr (\Lambda_I)^k(\eta)$, for 
$\eta\in \Lambda^k(M, \End B)$.
The following claim is implied by an elementary linear-algebraic 
computation.

\hfill

\claim\label{_degree_2-forms-in-End(B)_line-alge_Claim_}
Let $M$ be a hyperk\"ahler manifold, $B$ a Hermitian 
vector bundle over $M$, and $\Theta$ a 2-form on $M$
with coefficients in $\goth{su}(B)$. Assume that
\[ \Lambda_I \Theta =0, \ \ \  \Theta\in \Lambda^{1,1}_I(M, \End B)
\]
for some induced complex structure $I$. 
Assume, moreover, that $\Theta$ is square-integrable.
Let $J$ be another induced complex structure, $J\neq \pm I$.
Then 
\[ \deg_I \Theta^2 \geq |\deg_J \Theta^2|, \] 
and the equality is reached only if $\Theta$ is 
$SU(2)$-invariant.

{\bf Proof:} 
The following general argument is used.

\hfill

\sublemma\label{_deg_coeff_End(B)_Sublemma_}
Let $M$ be a K\"ahler manifold, $B$ a Hermitian 
vector bundle over $M$, and $\Xi$ a square-integrable 
2-form on $M$ with coefficients in $\goth{su}(B)$. 
Then: 
\begin{description}
\item[(i)] For 
\[ \Lambda_I \Xi =0, \ \ \ \Xi\in \Lambda^{1,1}_I(M, \End B)
\]
we have
\[ \deg_I \Xi^2 = C \int_M|\Xi|^2 \Vol M , \]
where $C= (4\pi^2 n (n-1))^{-1} $  $M$.
\item[(ii)] For 
\[  \Xi\in \Lambda^{2,0}_I(M, \End B)\oplus \Lambda^{0,2}_I(M, \End B)
\]
we have
\[ \deg_I \Xi^2 = -C \int_M |\Xi|^2\Vol M, \]
where $C$ is the same constant as appeared in (i).
\end{description}

{\bf Proof:} The proof is based on a linear-algebraic computation
(so-called L\"ubcke-type argument). The same computation is used to
prove Hodge-Riemann bilinear relations.
\endproof

\hfill

Return to the proof of \ref{_degree_2-forms-in-End(B)_line-alge_Claim_}.
Let $\Theta= \Theta^{1,1}_J + \Theta^{2,0}_J + \Theta^{0,2}_J$
be the Hodge decomposition associated with $J$. 
The following Claim shows that
$\Theta^{1,1}_J$ satisfies conditions of
\ref{_deg_coeff_End(B)_Sublemma_} (i).

\hfill

\claim
Let $M$ be a hyperk\"ahler manifold,
$I$, $L$ induced complex structures and
$\Theta$ a 2-form on $M$ satisfying 
\[ \Lambda_I \Theta =0, \ \ \ \Theta\in \Lambda^{1,1}_I(M).
\]
Let $\Theta^{1,1}_L$ be the $(1,1)$-component of $\Theta$
taken with respect to $L$. Then $\Lambda_L \Theta^{1,1}_L =0$.

\hfill

{\bf Proof:} Clearly, $\Lambda_L \Theta^{1,1}_L= \Lambda_L\Theta$.
Consider the natural Hermitian structure
on the space of 2-forms. Since $\Theta$ is of type $(1,1)$ with respect
to $I$, $\Theta$ is fiberwise orthogonal to 
the holomorphic \index{terms}{holomorphic symplectic structure} symplectic form 
$\Omega= \omega_J +\1 \Omega_K\in \Lambda^{2,0}_I(M)$.
By the same reason, $\Theta$ is orthogonal to $\bar\Omega$.
Therefore, $\Theta$ is orthogonal to 
$\omega_J$ and $\omega_K$. Since $\Lambda_I \Theta =0$,
$\Theta$ is also orthogonal to $\omega_I$. The map
$\Lambda_L$ is a projection to the form $\omega_L$
which is a linear combination of $\omega_I$, $\omega_J$ 
and $\omega_K$. Since $\Theta$ is fiberwise orthogonal 
to $\omega_L$, we have $\Lambda_L \Theta =0$. \endproof

\hfill

By
\ref{_deg_coeff_End(B)_Sublemma_},
we have
\[ \deg_J \left(\Theta^{1,1}_J\right)^2 = C \int_M|\Theta^{1,1}_J|^2
\]
and 
\[ \deg_J \left(\Theta^{2,0}_J + 
   \Theta^{0,2}_J\right)^2 = - C \int_M|\Theta^{2,0}_J + 
   \Theta^{0,2}_J|^2.
\]
Thus,
\[ \deg_J \Theta^2  = C \int_M\left|\Theta^{1,1}_J\right|~2- C \int_M\left|\Theta^{2,0}_J + 
   \Theta^{0,2}_J\right|~2.
\]
On the other hand,
\[ \deg_I \Theta^2 = C \int_M\left|\Theta\right|~2 =
   C \int_M\left|\Theta^{1,1}_J\right|~2 + C \int_M\left|\Theta^{2,0}_J + 
   \Theta^{0,2}_J\right|~2.
\]
Thus, $\deg_I \Theta^2> |\deg_J \Theta^2|$ unless
$\Theta^{2,0}_J + 
   \Theta^{0,2}_J=0$. On the other hand,
$\Theta^{2,0}_J + 
   \Theta^{0,2}_J=0$ means that $\Theta$ is of type $(1,1)$
with respect to $J$. Consider the standard $U(1)$-action on 
differential forms associated with the complex structures 
$I$ and $J$. These two $U(1)$-actions generate
the whole Lie group $SU(2)$ acting on $\Lambda^2(M)$
(here we use that $I\neq \pm J$). Since $\Theta$ is of type $(1,1)$
with respect to $I$ and $J$, this form is 
$SU(2)$-invariant. This proves 
\ref{_degree_2-forms-in-End(B)_line-alge_Claim_}.\endproof

\hfill

Return to the proof of \ref{_hyperho_conne_exi_Theorem_}.
Let $\nabla$ be the admissible 
\index{terms}{connections!Yang--Mills!admissible} 
Yang-Mills connection in $F$,
and $\Theta$ its curvature.
Recall that the form $Tr \Theta^2$ represents
the cohomology class 
$2c_2(F) - \frac{r-1}{r} c_1(F)^2$, where
$c_i$ are Chern classes of $F$.
Since the form $Tr \Theta^2$ is square-integrable, the integral
\[ \deg_J \Theta^2= \int_M Tr \Theta^2\omega_J^{n-2} \]
makes sense. In \cite{_Bando_Siu_}, it was shown
how to approximate the connection $\nabla$ by smooth
connections, via the heat equation.
This argument, in particular, was used to show that
the value of integrals like 
$\int_M Tr \Theta^2\omega_J^{n-2}$
can be computed through cohomology classes
and the Gauss--Bonnet formula  
\[ Tr \Theta^2 = 2c_2(F) - \frac{r-1}{r} c_1(F)^2.\]
Since the classes $c_2(F)$, $c_1(F)$ are 
$SU(2)$-invariant, we have
\[ \deg_I\Theta^2= \deg_J \Theta^2 \]
for all induced complex structures $I$, $J$.
By \ref{_degree_2-forms-in-End(B)_line-alge_Claim_},
this implies that $\Theta$ is $SU(2)$-invariant.
\ref{_hyperho_conne_exi_Theorem_} is proven.
\endproof

\hfill

The same argument implies the following corollary.

\hfill

\corollary\label{_stable_shea_degree_Corollary_}
Let $M$ be a compact hyperk\"ahler manifold, $I$
an induced complex structure, $F$ a stable 
\index{terms}{coherent sheaves!stable} 
\index{terms}{coherent sheaves!reflexive}
reflexive sheaf on
$(M,I)$, and $J$ be an induced complex structure,
$J\neq \pm I$. Then
\[ 
   \deg_I\left(2c_2(F) - \frac{r-1}{r} c_1(F)^2\right) \geq 
   \left|\deg_J\left(2c_2(F) - \frac{r-1}{r} c_1(F)^2\right)\right|,
\]
\index{terms}{degree!of a sheaf}
and the equality holds if and only if 
$F$ is hyperholomorphic.

\index{terms}{coherent sheaves!Wirtinger's inequality for}

\endproof

\subsection{Tensor category of hyperholomorphic sheaves}

This subsection is extraneous. Further on, we do not use the
tensor structure on the category of hyperholomorphic sheaves. 
However, we need the canonical identification of the categories
of hyperholomorphic sheaves associated with different 
induced complex structures.

{}From \index{names}{Bando, S.} Bando-Siu (\ref{_UY_for_shea_Theorem_}) 
it follows that on a compact K\"ahler manifold
a tensor product of stable 
\index{terms}{coherent sheaves!stable!tensor product of}
reflexive sheaves is polystable. Similarly, 
\ref{_hyperho_conne_exi_Theorem_} implies that a tensor product of polystable
hyperholomorphic sheaves is polystable hyperholomorphic. 
\index{terms}{coherent sheaves!hyperholomorphic!tensor category of}
We define the following
category. 

\hfill

\definition
Let $M$ be a compace hyperk\"ahler manifold and $I$ an induced complex
structure. Let $\c F_{st}(M,I)$ be a category with objects
\index{terms}{coherent sheaves!reflexive!category of}
reflexive polystable
hyperholomorphic sheaves and morphisms as in category
of coherent sheaves. This category is obviously additive. The tensor product
on $\c F_{st}(M,I)$ is induced from the tensor product of coherent sheaves.

\hfill

\claim
The category $\c F_{st}(M,I)$ is abelian. Moreover, it is a Tannakian tensor
category.
\index{terms}{Tannakian categories}

{\bf Proof:} Let $\phi:\; F_1 \arrow F_2$
be a morphism of hyperholomorphi sheaves. In 
\ref{_degree,slope_destabilising_Definition_} , we introduced
{\bf a slope} of a coherent sheaf. Clearly, 
\index{terms}{slope!of a hyperholomorphic sheaf} 
$sl(F_1)\leq sl(\im \phi)\leq sl(F_2)$.
All \index{terms}{coherent sheaves!hyperholomorphic!slope of} 
hyperholomorphic sheaves have slope 0 by 
\ref{_slope_hyperho_Remark_}. Thus,
$sl(\im \phi)=0$ and the subsheaf $\im\phi\subset F_2$ is 
destabilizing. Since $F_2$ is polystable, this sheaf is decomposed:
\[ F_2 = \im \phi \oplus \coker \phi.\]
A similar argument proves that
$F_1 = \ker\phi \oplus \coim \phi$, with all summands hyperholomorphic. 
This proves that $\c F_{st}(M,I)$ is abelian. The Tannakian properties 
\index{terms}{Tannakian categories}
are clear. \endproof

\hfill

The category $\c F_{st}(M,I)$ does not depend from the choice of 
induced complex structure $I$:

\hfill

\theorem\label{_equi_cate_Theorem_}
Let $M$ be a compact hyperk\"ahler 
manifold, $I_1$, $I_2$ induced complex structures and
$\c F_{st}(M,I_1)$, $\c F_{st}(M,I_2)$ the associated categories
of polystable \index{terms}{coherent sheaves!reflexive}
reflexive hyperholomorphic sheaves. 
Then, there exists a natural equivalence of tensor categories
\[ 
   \Phi_{I_1, I_2}:\; \c F_{st}(M,I_1)\arrow \c F_{st}(M,I_2).
\]

{\bf Proof:} Let $F\in \c F_{st}(M,I_1)$ be a reflexive
polystable hyperholomorphic sheaf
and $\nabla$ the canonical admissible 
\index{terms}{connections!Yang--Mills!admissible} 
Yang-Mills connection.
Consider the sheaf $\c F$ on the \index{terms}{twistor space} twistor space $\Tw(M)$ constructed
as in the proof of \ref{_conne_=>_hyperho_Proposition_}. 
Restricting $\c F$ to $\pi^{-1}(I_2)\subset \Tw(M)$,
we obtain a coherent sheaf $F'$ on $(M,I_2)$. As we have shown in
the proof of \ref{_conne_=>_hyperho_Proposition_}, 
the sheaf $(F')^{**}$ is polystable hyperholomorphic. 
Let $\Phi_{I_1, I_2}(F):= (F')^{**}$.
It is easy to check that thus constructed map
of objects gives a functor 
\[ 
   \Phi_{I_1, I_2}:\; \c F_{st}(M,I_1)\arrow \c F_{st}(M,I_2),
\]
and moreover, $\Phi_{I_1, I_2}\circ \Phi_{I_2, I_1} = Id$.
This shows that $\Phi_{I_1, I_2}$ is an equivalence.
\ref{_equi_cate_Theorem_} is proven. \endproof

\hfill

\definition \label{_hyperho_shea_on_M_Definition_}
By \ref{_equi_cate_Theorem_}, the category
$\c F_{st}(M,I_1)$ is independent from the choice of
induced complex structure. We call this category
{\bf the category of polystable 
\index{terms}{coherent sheaves!hyperholomorphic!category of} 
hyperholomorphic reflexive 
sheaves on $M$} and denote it by $\c F(M)$.
The objects of $\c F(M)$ are called {\bf hyperholomorphic sheaves on
$M$}. For a hyperholomorphic sheaf on $M$, 
we denote by $F_I$ the corresponding
sheaf from $\c F_{st}(M,I_1)$.

\hfill

\remark
Using the same argument as proves \ref{_admi_twi_impli_Theorem_} (ii),
it is easy to check that
the category $\c F(M)$ is a deformational invariant of $M$. That is,
for two hyperk\"ahler manifolds $M_1$, $M_2$ which are
deformationally equivalent, the categories $\c F(M_i)$ are
also equivalent, assuming that $Pic(M_1)= Pic(M_2)=0$. 
The proof of this result is essentially contained in
\cite{_coho_announce_}.

\hfill

\remark
As \index{names}{Deligne, P.} Deligne proved (\cite{_Deli:Tanna_}),
\index{terms}{Tannakian categories}
for a each Tannakian category $\c C$ equipped with a fiber functor,
there exists a natural pro-algebraic group $G$ such that
$\c C$ is a group of representations of $G$. For $\c F(M)$, there
are several natural fiber functors. The simplest one is defined for
each induced complex structure $I$ such that $(M,I)$ is algebraic
(such complex structures always exist, as proven in
\cite{_Fujiki_}; see also \cite{_Verb:alge_} and 
Subsection \ref{_alge_indu_Subsection_}). 
Let $\c K(M,I)$ is the space of rational functions on $(M,I)$.
For $F\in \c F_{st}(M, I)$, consider the functor
$F\arrow \eta_I(F)$, where $\eta_I(F)$ is
the space of global sections of $F\otimes \c K(M,I)$.
This is clearly a fiber functor, which associates
to $\c F(M)$ the group $G_I$. The corresponding pro-algebraic
group $G_I$ is a deformational, that is, topological,
invariant of the hyperk\"ahler
manifold.


\section{Cohomology of hyperk\"ahler manifolds}
\label{_cohomo_hype_Section_}


This section contains a serie of preliminary results which
are used further on to define and study the \index{terms}{$C$-restricted complex structures} $C$-restricted
complex structures.

\subsection{Algebraic induced complex structures}
\label{_alge_indu_Subsection_}

This subsection contains a recapitulation of results
of \cite{_Verb:alge_}.

\hfill

A more general version of the 
following theorem was proven by
A. \index{names}{Fujiki, A.} Fujiki (\cite{_Fujiki_}, Theorem 4.8 (2)).

\hfill

\theorem \label{_alge_dense_Theorem_}
Let $M$ be a compact simple hyperk\"ahler manifold and
$\c R$ be the set of induced complex structures 
$\c R \cong \C P^1$. Let $\c R_{alg}\subset \c R$
be the set of all algebraic induced complex structures. 
Then $\c R_{alg}$ is countable and dense in $\c R$. 

{\bf Proof:} This is \cite{_Verb:alge_}, Theorem 2.2. \endproof

\hfill

In the proof of \ref{_alge_dense_Theorem_}, the following
important lemma was used.

\hfill

\lemma \label{_K_proje_on_R_Lemma_}
\begin{description}
\item [(i)] Let $\calo\subset H^2(M, \R)$ be the set of all
cohomology classes 
which are K\"ahler with respect to some induced comples structure.
Then $\calo$ is open in $H^2(M, \R)$. Moreover, for all
$\omega\in \calo$, the class $\omega$ is {\it not} $SU(2)$-invariant.
\item[(ii)]
Let $\eta\in H^2(M, \R)$ 
be a cohomology class which is not $SU(2)$-invariant. Then 
   there exists a unique up to a sign induced complex structure
$I\in \c R/\{\pm 1\}$ such that $\eta$ belongs to $H^{1,1}_I(M)$.
\end{description}

{\bf Proof:} This statement is a form of
\cite{_Verb:alge_}, Lemma 2.3. \endproof

\subsection{The action of $\goth{so}(5)$ on the 
cohomology of a hyperk\"ahler manifold}
\label{_so(5)_Subsection_}

This subsection is a recollection of data from
\cite{_so(5)_} and \cite{_Verbitsky:Symplectic_II_}.

\hfill

Let $M$ be a hyperk\"ahler manifold. For an induced complex structure
$R$ over $M$, consider the K\"ahler form 
$\omega_R=(\cdot,R\cdot)$, where $(\cdot,\cdot)$ is the Riemannian form.
 As usually, $L_R$ denotes the operator of exterior
multiplication by $\omega_R$, which is acting on the differential
forms $A^*(M,\C)$ over $M$. Consider the adjoint operator to $L_R$, 
denoted by $\Lambda_R$.

One may ask oneself, what algebra is generated by $L_R$ and 
$\Lambda_R$ for all induced complex structures $R$? The answer was given
in \cite{_so(5)_}, where the following theorem was proven.

\hfill

\theorem
 (\cite{_so(5)_}) Let $M, \c H$ be a hyperk\"ahler manifold,
and $\goth{a}_{\c H}$ be a Lie algebra generated by $L_R$ and $\Lambda_R$
for all induced complex structures $R$ over $M$. Then the
Lie algebra  $\goth{a}_{\c H}$ is isomorphic to $\goth{so}(4,1)$.
\endproof

\hfill

The following facts about a structure of $\goth{a}_{\c H}$
were also proven in \cite{_so(5)_}.

Let $I$, $J$ and $K$ be three induced complex structures on
$M$, such that $I\circ J=-J\circ I=K$. For an induced complex structure
$R$, consider an operator $ad R$ on cohomology, acting on $(p,q)$-forms
as a multiplication by $(p-q)\1$. The operators $ad R$ generate a 3-dimensional
Lie algebra $\g_{\c H}$, which is isomorphic to $\goth{su}(2)$.
This algebra coincides with the Lie algebra associated to the
standard $SU(2)$-action on $H^*(M)$.
The algebra $\goth{a}_{\c H}$ contains
$\g_{\c H}$ as a subalgebra, as follows: 
\begin{equation}\label{_ad_I_as_commu_Equation_} 
   [\Lambda_J,L_K]=[L_J,\Lambda_K]= ad\: I\;\; \text{(etc)}.
\end{equation}
The algebra  $\goth{a}_{\c H}$ is 10-dimensional. 
It has the following basis: $L_R,\Lambda_R$, $ad\: R$
$(R=I,J,K)$ and the element $H=[L_R,\Lambda_R]$.
The operator $H$ is a standard Hodge operator; it acts on $r$-forms
over $M$ as multiplication by a scalar $n-r$,
where $n=dim_\C M$.

\hfill

\definition \label{_isoty_Definition_}
Let $\g$ be a semisimple Lie algebra, $V$ its representation
and $V= \oplus V_{\alpha}$ a $\g$-invariant decompostion of $V$,
such that for all $\alpha$, $V_\alpha$ is a direct sum of
isomorphic finite-dimensional representations $W_\alpha$ of $V$, and
all $W_\alpha$ are distinct. Then the decomposition $V= \oplus V_{\alpha}$
is called {\bf the isotypic decomposition of $V$}.

\hfill

It is clear that for all
finite-dimensional representations, isotypic decomposition always
exists and is unique.

\hfill

Let $M$ be a compact hyperk\"ahler manifold.
Consider the cohomology space $H^*(M)$ equipped with the natural
action of $\goth{a}_{\c H}=
\goth{so}(5)$. Let $H_o^*\subset H^*(M)$ be the isotypic
component containing $H^0(M)\subset H^*(M)$. 
Using the root system written explicitly for 
$\goth{a}_{\c H}$ in \cite{_so(5)_}, \cite{_Verbitsky:cohomo_},
it is easy to check that $H^*_o(M)$ is an irreducible
representation of $\goth{so}(5)$. Let 
$p:\; H^*(M) \arrow H^*_o(M)$ be the unique
$\goth{so}(5)$-invariant projection, and
$i:\; H^*_o(M) \hookrightarrow H^*(M)$ the natural embedding.

\hfill

Let $M$ be a compact hyperk\"ahler manifold, $I$ an induced complex
structure, and $\omega_I$ the corresponding K\"ahler form.
Consider the {\bf degree 
\index{terms}{degree!of a form}
map} $\deg_I:\; H^{2p}(M) \arrow \C$,
$\eta\arrow \int_M \eta \wedge \omega_I^{n-p}$,
where $n= \dim_C M$.

\hfill

\proposition\label{_degree_isotypic_Proposition_}
The space \[ H^*_o(M)\subset H^*(M)\] is a subalgebra of $H^*(M)$,
which is invariant under the $SU(2)$-action. Moreover,
for all induced complex structures $I$, the degree map
\[\deg_I:\; H^*(M) \arrow \C\] satisfies
\[ \deg_I(\eta) =\deg_I(i(p(\eta)),
\]
where $i:\; H^*_o(M) \hookrightarrow H^*(M)$, $p:\; H^*(M) \arrow H^*_o(M)$
are the $\goth{so}(5)$-invariant maps defined above.
And finally, the projection $p:\; H^*(M) \arrow H^*_o(M)$
is $SU(2)$-invariant.

\hfill

{\bf Proof:} The space $H^*_o(M)$ is generated from
$\pmb 1\in H^0(M)$ by operators $L_R$, $\Lambda_R$. 
To prove that $H^*_o(M)$ is closed under multiplication,
we have to show that $H^*_o(M)$ is generated (as a linear space)
by expressions of type $L_{r_1} \circ L_{R_2} \circ ... \circ \pmb 1$.
By \eqref{_ad_I_as_commu_Equation_},
the commutators of $L_R$, $\Lambda_R$ map such expressions
to linear combinations of such expressions. On the other hand,
the operators $\Lambda_R$ map $\pmb 1$ to zero. Thus,
the operators $\Lambda_R$ map expressions of type 
$L_{r_1} \circ L_{R_2} \circ ... \circ \pmb 1$ to linear combinations
of such expressions. This proves that $H^*_o(M)$ is closed under multiplication.
The second statement of \ref{_degree_isotypic_Proposition_}
is clear (see, e. g. \cite{_Verbitsky:Symplectic_II_}, 
proof of Proposition 4.5). 
It remains to show that $H^*_o(M)\subset H^*(M)$ is an $SU(2)$-invariant 
subspace and that  $p:\; H^*(M) \arrow H^*_o(M)$
is compatible with the $SU(2)$-action. From 
\eqref{_ad_I_as_commu_Equation_}, we obtain that
the Lie group $G_A$ associated with $\goth a_{\c H}\cong \goth{so}(1,4)$ 
contains $SU(2)$ acting in a standard way
on $H^*(M)$. Since the map $p:\; H^*(M) \arrow H^*_o(M)$
commutes with $G_A$-action, $p$ also commutes with $SU(2)$-action.
We proved \ref{_degree_isotypic_Proposition_}.
\endproof

\subsection{Structure of the cohomology ring}
\label{_cohomo_stru_Subsection_}

In \cite{_Verbitsky:cohomo_} (see also \cite{_coho_announce_}),
we have computed explicitly the subalgebra of cohomology of $M$ generated
by $H^2(M)$. This computation can be summed up as follows.

\hfill

\theorem \label{_S^*H^2_is_H^*M_intro-Theorem_} 
(\cite{_Verbitsky:cohomo_}, Theorem 15.2)
Let $M$ be a compact hyperk\"ahler manifold, $H^1(M)=0$,
$\dim_\C M=2n$, and $H^*_r(M)$ the subalgebra of cohomology of $M$ generated
by $H^2(M)$.
Then 
\[\bigg\{\begin{array}{lr}
 H^{2i}_r(M)\cong S^i H^2(M)&
 \mbox{\ \ for $i\leq n$, and}\\
 H^{2i}_r(M)\cong S^{2n-i} H^2(M) & 
 \mbox{\ \ for $i\geq n$ }
 \end{array}
\]

\endproof

\hfill

\theorem\label{_gene_all_SU(2)_Theorem_}
Let $M$ be a simple hyperk\"ahler manifold.
Consider the group $G$ generated by a union of all $SU(2)$ for
all hyperk\"ahler 
\index{terms}{hyperk\"ahler structures} 
structures on $M$. Then the Lie algebra of 
$G$ is isomorphic to $\goth{so}(H^2(M))$,
for a certain natural integer bilinear symmetric form on $H^2(M)$,
called \index{terms}{Bogomolov--Beauville pairing}
Bogomolov-Beauville form. 

{\bf Proof:} \cite{_Verbitsky:cohomo_} (see also \cite{_coho_announce_}).
\endproof

\hfill

The key element in the proof of
\ref{_S^*H^2_is_H^*M_intro-Theorem_} and \ref{_gene_all_SU(2)_Theorem_}
is the following algebraic computation. 

\hfill

\theorem \label{_r_proper_Theorem_}
Let $M$ be a simple hyperk\"ahler manifold,
and $\c H$ a hyperk\"ahler 
\index{terms}{hyperk\"ahler structures!$\protect\goth{so}(1,4)$-action associated with} 
structure on $M$.
Consider
the Lie subalgebra 
\[ {\goth a}_{\c H}\subset \End(H^*(M)), 
   \ \ {\goth a}_{\c H} \cong \goth{so}(1,4),
\]
associated with the hyperk\"ahler structure 
(Sub\-sec\-tion \ref{_so(5)_Subsection_}). Let 
 \[ \g\subset \End(H^*(M)) \]
be the Lie algebra generated by subalgebras 
${\goth a}_{\c H}\subset \End(H^*(M))$,
for all hyperk\"aher structures $\c H$ on $M$.
Then 
\begin{description} 
\item[(i)] The algebra $\g$ is naturally isomorphic to 
the Lie algebra $\goth{so}(V\oplus \goth H)$,
where $V$ is the linear space $H^2(M, \R)$ 
equipped with the \index{terms}{Bogomolov--Beauville pairing}
Bogomolov--Beauville pairing,
and $\goth H$ is a 2-dimensional vector space 
with a quad\-ra\-tic form of signature $(1, -1)$.
\item[(ii)] 
The space $H_r^*(M)$ is invariant under the action of $\g$, 
Moreover, \[ H_r^*(M)\subset H^*(M)\] is an isotypic
\footnote{See \ref{_isoty_Definition_} for the definition of
isotypic decomposition.} component of
the space $H^*(M)$ considered as a representation of $\g$.
\end{description}

{\bf Proof:} \cite{_Verbitsky:cohomo_} (see also \cite{_coho_announce_}).
\endproof

\hfill

As one of the consequences of 
\ref{_S^*H^2_is_H^*M_intro-Theorem_}, we obtain the
following lemma, which will be used further on in this paper.

\hfill

\lemma\label{_p_multi_Lemma_}
Let $M$ be a simple hyperk\"ahler manifold,
$\dim_{\Bbb H} M =n$, and $p:\; H^*(M) \arrow H^*_o(M)$ the
map defined in Subsection \ref{_so(5)_Subsection_}.
Then, for all $x, y \in H^*_r(M)$, we have
\[ p(x) p(y) = p(xy), \text{\ \  whenever\ \ } 
   xy \in \bigoplus\limits_{i\leq 2n} H^i(M).
\]
{\bf Proof:}  Let $\omega_I$, $\omega_J$, 
$\omega_K$,
$x_1$, ..., $x_n$ be an orthonormal basis in $H^2(M)$.
Clearly, the vectors $x_1$, ..., $x_k$ are $SU(2)$-invariant.
Therefore, these vectors are highest vectors of the corresponding 
${\goth a}_{\c H}$-representations, with respect
to the root system and Cartan subalgebra for ${\goth a}_{\c H}$
which is written in \cite{_so(5)_} or \cite{_Verbitsky:cohomo_}.
We obtain that the monomials 
\[ P_{k_1, k_2, k_3, \{n_i\} }
   = \omega_I^{k_1} \omega_J^{k_2} \omega_K^{k_3} \prod x_i^{n_i},
   \sum n_i = N, \ \ P_{k_1, k_2, k_3, \{n_i\} }
   \in \bigoplus\limits_{i\leq 2n} H^i(M)
\] 
belong to the different isotypical components for different $N$'s.
By \ref{_S^*H^2_is_H^*M_intro-Theorem_}, a product of 
two such monomials $P_{k_1, k_2, k_3, \{n_i\} }$ and
$P_{k_1', k_2', k_3', \{n_i'\} }$ is equal to 
$P_{k_1+k_1', k_2+k_2', k_3+k_3', \{n_i+n_i'\} }$,
assuming that 
\[ P_{k_1, k_2, k_3, \{n_i\} }
   P_{k_1', k_2', k_3', \{n_i'\} }\in \bigoplus\limits_{i\leq 2n} H^i(M).
\]
Thus,
the isotypical decomposition associated with the 
$\goth{a}_{\c H}$-action is compatible
with multiplicative structure on $H^*(M)$, 
for low-dimensional cycles.
This implies \ref{_p_multi_Lemma_}.
\endproof

\hfill

We shall use the following corollary of 
\ref{_p_multi_Lemma_}.

\hfill

\corollary \label{_p_prods_H^2_Corollary_}
Let $M$ be a simple hyperk\"ahler manifold,
$\dim_{\Bbb H} M>1$, and $\omega_1, \omega_2 \in H^2(M)$
cohomology classes which are $SU(2)$-invariant. 
Then, for all induced complex structures $I$, we have
$\deg_I(\omega_1 \omega_2) =0$.
\index{terms}{degree!of a form}

{\bf Proof:} By definition, the classes $\omega_1, \omega_2$
satisfy $\omega_i\in \ker p$. 
By \ref{_p_multi_Lemma_}, we have
$\omega_1 \omega_2\in \ker p$. By
\ref{_degree_isotypic_Proposition_}, 
$\deg_I\omega_1 \omega_2 =0$.
\endproof

\hfill

Let $\omega$ be a rational K\"ahler form. The corresponding
$\goth{sl}(2)$-action on $H^*(M)$ is clearly compatible with the
rational structure on $H^*(M)$. It is easy to see (using, for instance,
\ref{_K_proje_on_R_Lemma_}) that $\g$ is generated by 
$\goth{sl}(2)$-triples associated with rational K\"ahler forms
$\omega$. Therefore,
the action of $\g$ on $H^*(M)$ is compatible with the rational structure
on $H^*(M)$. Using the isotypic decomposition, we define
a natural $\g$-invariant map $r:\; H^*(M) \arrow H^*_r(M)$.
Further on, we shall use the following properties of this map.

\hfill

\claim \label{_r_H^*_to_H_r_Claim_}
\begin{description}
\item[(i)] The map $r:\; H^*(M) \arrow H^*_r(M)$
is compatible with the rational structure on $H^*(M)$.
\item[(ii)] For every $x\in \ker r$, and every hyperk\"ahler
\index{terms}{hyperk\"ahler structures} structure $\c H$, the corresponding map
$p:\; H^*(M) \arrow H^*_o(M)$ satisfies $p(x)=0$.
\item[(iii)] For every $x\in \ker r$, every hyperk\"ahler
structure $\c H$, and every induced complex
structure $I$ on $M$,  we have $\deg_I x =0$.
\end{description}

{\bf Proof:}
\ref{_r_H^*_to_H_r_Claim_} (i) is clear, because 
the action of $\g$ on $H^*(M)$ is compatible with the rational structure
on $H^*(M)$. 
To prove \ref{_r_H^*_to_H_r_Claim_} (ii), we notice that 
the space $H^*_r(M)$ is generated from $H^0(M)$ by
the action of $\g$, and $H^*_o(M)$ is generated from $H^0(M)$ by
the action of ${\goth a}_{\c H}$. Since ${\goth a}_{\c H}$ is by definition a 
subalgebra in $\g$, we have $H^*_o(M) \subset H^*_r(M)$.
The isotypic projection
$r:\; H^*(M) \arrow H^*_r(M)$ is by definition compatible with 
the $\g$-action. Since ${\goth a}_{\c H}\subset \g$, the map $r$
is also compatible with the ${\goth a}_{\c H}$-action.
Therefore, $\ker r\subset \ker p$.
\ref{_r_H^*_to_H_r_Claim_} (iii) is implied
by \ref{_r_H^*_to_H_r_Claim_} (ii) and
\ref{_degree_isotypic_Proposition_}.
\endproof

\hfill

Let $x_i$ be an basis in $H^2(M, {\Bbb Q})$ which is 
rational and orthonormal 
with respect to \index{terms}{Bogomolov--Beauville pairing}
Bogomolov-Beauville pairing,
$(x_i, x_i)_{\c B} = \epsilon_i = \pm 1$. 
Consider the cohomology class 
$\theta':= \epsilon_i x_i^2\in H^4(M, {\Bbb Q})$.
Let $\theta \in H^4(M, \Z)$ be a non-zero integer cohomology
class which is proportional to $\theta'$.
{}From results of \cite{_Verbitsky:cohomo_} (see also \cite{_coho_announce_}),
the following proposition can be easily deduced.

\hfill

\proposition\label{_theta_SU(2)_inva_Proposition_}
The cohomology class $\theta\in H^4(M, \Z)$ 
is $SU(2)$-invariant for all hyperk\"ahler 
\index{terms}{hyperk\"ahler structures!generic} structures on $M$.
Moreover, for a generic hyperk\"ahler structure, 
the group of $SU(2)$-invariant integer classes $\alpha\in H^4_r(M)$
has rank one, where $H^*_r(M)$ is the subalgebra of 
cohomology generated by $H^2(M)$.

\hfill

{\bf Proof:} Clearly, if an integer class $\alpha$
is $SU(2)$-invariant for a generic hyperk\"ahler structure, then
$\alpha$ is $G$-invariant, where $G$ is the group defined in
\ref{_gene_all_SU(2)_Theorem_}. On the other hand, $H^4_r(M)\cong S^2(H^2(M))$,
as follows from \ref{_S^*H^2_is_H^*M_intro-Theorem_}. Clearly,
the vector 
$\theta\in H^4_r(M)\cong S^2(H^2(M))$ is $\goth{so}(H^2(M))$-invariant.
Moreover, the space of $\goth{so}(H^2(M))$-invariant vectors in
$S^2(H^2(M))$ is one-dimensional. Finally, from an explicit computation
of $G$ it follows that $G$ acts on $H^4(M)$ as $SO(H^2(M))$,
and thus, the Lie algebra invariants coincide with invariants
of $G$. We found that the space of $G$-invariants in
$H^4_r(M)$ is one-dimensional and generated by $\theta$.
This proves \ref{_theta_SU(2)_inva_Proposition_}.
\endproof

\hfill

\remark
It is clear how to generalize \ref{_theta_SU(2)_inva_Proposition_}
from dimension 4 to all dimensions. 
The space $H^{2d}_r(M)^G$ of $G$-invariants 
in $H^{2d}_r(M)$ is 1-dimensional 
for $d$ even and zero-dimensional for $d$ odd.

\subsection{Cohomology classes of CA-type}

Let $M$ be a compact hyperk\"ahler manifold, and $I$ an induced
complex structure.
All cohomology classes which appear as fundamental 
classes of complex subvarieties of $(M, I)$ 
satisfy certain properties. Classes satisfying
these properties are called classes of CA-type, from
Complex Analytic. Here is the definition of CA-type.

\hfill

\definition
\index{terms}{CA-type (cohomology classes of)!definition} 
 Let $\eta\in H^{2,2}_I(M) \cap H^4(M, \Z)$ be an
integer (2,2)-class. Assume that 
for all induced complex structures $J$, 
satisfying $I\circ J = - J\circ I$,
we have $\deg_I(\eta)\geq \deg_J(\eta)$,
and the equality is reached only if
$\eta$ is $SU(2)$-invariant.
Assume, moreover, that
$\deg_I(\eta)\geq |\deg_J(\eta)|$. 
Then $\eta$ is called {\bf a class
of CA-type}.

\hfill

\theorem \label{_funda_and_Chern_CA_Theorem_}
Let $M$ be a simple hyperk\"ahler manifold, 
of dimension $\dim_{\Bbb H} M>1$, $I$ an induced
complex structure, and $\eta \in H^{2,2}_I(M) \cap H^4(M, \Z)$
an integer (2,2)-class. Assume that one of the following
conditions holds.
\begin{description}
\item[(i)] There exists a complex subvariety 
$X\subset (M, I)$ such that $\eta$ is the fundamental class of $X$
\item[(ii)] There exists a 
\index{terms}{coherent sheaves!stable!Chern classes of} 
stable coherent torsion-free sheaf $F$
over $(M,I)$, such that the first Chern class of $F$
is zero, and $\eta=c_2(F)$.
\end{description}
Then $\eta$ is of CA-type.
\index{terms}{CA-type (cohomology classes of)} 

\hfill

{\bf Proof:} \ref{_funda_and_Chern_CA_Theorem_} (i)
is a direct consequence of Wirtinger's inequality
(\ref{_Wirti_hyperka_Proposition_}).
It remains to prove \ref{_funda_and_Chern_CA_Theorem_} (ii).

We assume, temporarily, that $F$ is reflexive.
By \ref{_stable_shea_degree_Corollary_},
we have
\begin{equation} \label{_c_2,1_ineq_Equation_}
   \deg_I(2c_2(F) - \frac{r-1}{r} c_1(F)^2) \geq 
   \left|\deg_J(2c_2(F) - \frac{r-1}{r} c_1(F)^2)\right|,
\end{equation}
and the equality happens only if $F$ is hyperholomorphic.
Since $c_1(F)$ is $SU(2)$-invariant, we have
$\deg_I(c_1(F)^2) = \deg_J (c_1(F)^2) =0$
(\ref{_p_prods_H^2_Corollary_}).
Thus, \eqref{_c_2,1_ineq_Equation_}
implies that 
\[ \deg_I 2c_2(F) \geq |\deg_J 2c_2(F)| \]
and the inequality is strict unless $F$ is hyperholomorphic,
in which case, the class $c_2(F)$ is $SU(2)$-invariant by definition.
We have proven \ref{_funda_and_Chern_CA_Theorem_} (ii)
for the case of \index{terms}{coherent sheaves!reflexive}
reflexive $F$.

\hfill

For $F$ not necessary reflexive sheaf,
we have shown in the proof of \ref{_hyperho_suppo_of_F^**/F_Claim_}
that
\[ c_2(F) = c_2(F^{**}) + \sum n_i [X_i], \]
where $n_i$ are positive integers, and $[X_i]$ are the 
fundamental classes of irreducible components of support
of the sheaf $F^{**}/F$. Therefore, the class $c_2(F)$
is a sum of classes of CA-type. Clearly, a sum of
\index{terms}{CA-type (cohomology classes of)} 
cohomology classes of CA-type is again a class of
CA-type. This proves
\ref{_funda_and_Chern_CA_Theorem_}. \endproof


\section{$C$-restricted complex structures on hy\-per\-k\"ah\-ler manifolds}
\label{_C_restri_Section_}


\subsection{Existence of $C$-restricted complex structures}

We assume from now till the end of this section
that the hyperk\"ahler manifold $M$ is simple
(\ref{_simple_hyperkahler_mfolds_Definition_}). 
This assumption can be avoided, but
it simplifies notation. 

We assume from now till the end of this section
that the hyperk\"ahler manifold $M$ is compact 
of real dimension $\dim_\R M \geq 8$,
i. e. $\dim_{\Bbb H}M \geq 2$. This assumption is absolutely necessary.
The case of hyperk\"ahler surfaces with $\dim_{\Bbb H}M =1$ (torus and
K3 surface) is trivial and for our purposes not interesting.
It is not difficult to extend our definitions and
results to the case of a compact hyperk\"ahler manifold
which is a product of simple hyperk\"ahler manifolds
with $\dim_{\Bbb H}M \geq 2$.

\hfill

\index{terms}{$C$-restricted complex structures!definition of} 
\index{terms}{induced complex structures!$C$-restricted}
\definition\label{_C_restri_Definition_}
Let $M$ be a compact hyperk\"ahler manifold,
and $I$ an induced complex structure. As usually, we denote by
$\deg_I:\; H^{2p}(M) \arrow \C$ the associated degree 
map, and by $H^*(M)= \oplus H^{p,q}_I(M)$ 
the Hodge decomposition. Assume that $I$ is algebraic. Let $C$ be a positive
real number. We say that the induced complex structure
$I$ is {\bf $C$-restricted} if the following conditions hold.
\begin{description}
\item[(i)] For all non-$SU(2)$-invariant cohomology classes 
classes $\eta\in H^{1,1}_I(M) \cap H^2(M, \Z)$,
we have $|\deg_I(\eta)| > C$. 
\item[(ii)] Let $\eta \in H^{2,2}_I(M)$ be a cohomology class
of CA-type which is not $SU(2)$-invariant. 
\index{terms}{CA-type (cohomology classes of)} 
Then $|\deg_I(\eta)| > C$.
\end{description}

\hfill

The heuristic (completely informal) meaning of 
this definition is the following.
The degree map plays the role of the metric on the cohomology.
Cohomology classes with small degrees are ``small'', the rest is ``big''.
Under reasonably strong assumptions,
there are only finitely many ``small'' integer classes, and the
rest is ``big''. For each non-$SU(2)$-invariant cohomology class $\eta$
there exists at most two induced complex structures for which
$\eta$ is of type $(p,p)$. Thus, for most induced complex structures,
all non-$SU(2)$-invariant integer $(p,p)$ classes are ``big''.
Intuitively, the $C$-restriction means that
all non-$SU(2)$-invariant integer (1,1) and (2,2)-cohomology 
classes are ``big''. This definition is needed for the
study of first and second Chern classes of sheaves.
The following property of $C$-restricted complex structures
is used (see \ref{_funda_and_Chern_CA_Theorem_}):
for every subvariety $X\subset (M,I)$ 
of complex codimension 2, either $X$ 
is \index{terms}{trianalytic subvarieties} trianalytic or $\deg_I(X)>C$.

\hfill

\index{terms}{hyperk\"ahler structures!admitting $C$-restricted complex structures!definition of} 
\index{terms}{$C$-restricted complex structures} 

\definition \label{_admitti_C_restri_Definition_} 
Let $M$ be a compact manifold,
and $\c H$ a hyperk\"ahler structure on $M$. We say that
$\c H$ {\bf admits $C$-restricted complex structures}
if for all $C>0$, the set of all $C$-restricted algebraic 
complex structures is dense in the set $R_{\c H}= \C P^1$
of all induced complex structures.

\hfill

\proposition \label{_restri_for_H^11_1-dim-Proposition_}
Let $M$ be a compact simple hyperk\"ahler manifold, $\dim_{\Bbb H}(M) >1$,
and $r:\; H^4(M) \arrow H_r^4(M)$ be the map defined in 
\ref{_r_H^*_to_H_r_Claim_}. Assume that 
for all algebraic induced complex structures $I$, 
the group $H^{1,1}_I(M)\cap H^2(M, \Z)$ has rank one, and
the group \[ H^{2,2}_I(M)\cap H^4(M, \Z)/(\ker r)\] has rank 2.
Then $M$ admits 
\index{terms}{$C$-restricted complex structures!existence of} 
$C$-restricted complex structures.

\hfill

{\bf Proof:}
The proof of \ref{_restri_for_H^11_1-dim-Proposition_} takes the
rest of this section.

Denote by $\c R$ the set $\c R\cong \C P^1$ of all induced complex
structures on $M$. Consider 
the set $\c R/\{\pm 1\}$ of induced complex structures up to a sign
(\ref{_K_proje_on_R_Lemma_}). Let $\alpha \in H^2(M)$ be a cohomology class
which is not $SU(2)$-invariant. According to
\ref{_K_proje_on_R_Lemma_}, there exists a unique
element $c(\alpha)\in \c R/\{\pm 1\}$ such that
$\alpha \in H^{1,1}_{c(\alpha)}(M)$. 
This defines a map
\[ 
   c:\; \left( H^2(M, \R) \backslash H^2_{inv}(M)\right) \arrow \c R/\{\pm 1\}, 
\]
where $H^2_{inv}(M)\subset H^2(M)$ is the set of all
$SU(2)$-invariant cohomology classes.
For induced complex structures $I$ and $-I$, and $\eta\in H^{2p}(M)$,
the degree maps satisfy 
\index{terms}{degree!for different induced complex structures}
\begin{equation} \label{_deg_-I_Equation_}
\deg_I(\eta) = (-1)^{p}\deg_{-I}(\eta).
\end{equation}
Thus, the number $|\deg_I(\eta)|$ is independent from the sign
of $I$. 

Let $\eta\in H^*(M, \Z)$ be a cohomology class. The {\bf largest
divisor} of $\eta$ is the biggest positive integer number $k$
such that the cohomology class $\frak \eta k$ is also integer.

Let $\alpha\in H^2(M, \Z)$ be an integer cohomology class,
which is not $SU(2)$-invariant, 
$k$ its largest divisor and $\tilde \alpha:= \frak \alpha k$
the corresponding integer class. Denote by $\widetilde\deg(\alpha)$
the number
\[ 
   \widetilde\deg(\alpha) := \left| \deg_{c(\alpha)}(\tilde \alpha)\right|.
\]
The induced complex structure $c(\alpha)$ is defined up to a sign,
but from \eqref{_deg_-I_Equation_} it is clear that $\widetilde\deg(\alpha)$
is independent from the choice of a sign.

\hfill

\lemma \label{_C_restri_from_A,d_Lemma_}
Let $M$ be a compact hyperk\"ahler manifold, and
$I$ be an algebraic induced complex structure, such that
the group $H^{1,1}_I(M)\cap H^2(M, \Z)$ has rank one,
and the group $H^{2,2}_I(M)\cap H^4(M, \Z)/(\ker r)$ has rank 2.
Denote by $\alpha$ the generator of $H^{1,1}_I(M)\cap H^2(M, \Z)$.
Since the class $\alpha$ is proportional to a K\"ahler form, $\alpha$
is not $SU(2)$-invariant (\ref{_K_proje_on_R_Lemma_}, (i)).
Let $d:= \widetilde \deg\alpha$. 
Then, there exists a positive real constant $A$ depending
on volume of $M$, its topology
 and its dimension, such that $I$ is $d\cdot A$-restricted.

\hfill

{\bf Proof:} This lemma is a trivial calculation based on
results of \cite{_Verbitsky:cohomo_} (see also \cite{_coho_announce_} 
and Subsection \ref{_cohomo_stru_Subsection_}). 

Since $H^{1,1}_I(M)\cap H^2(M, \Z)$ has rank one,
for all $\eta \in H^{1,1}_I(M)\cap H^2(M, \Z)$, $\eta\neq 0$,
we have $| \deg_I \eta| \geq d$. This proves the first condition of
\ref{_C_restri_Definition_}.

Let $\theta$ be the $SU(2)$-invariant integer cycle $\theta\in H^4(M)$
defined in \ref{_theta_SU(2)_inva_Proposition_}. 
By \ref{_SU(2)_inva_type_p,p_Lemma_}, 
$\theta\in H^{2,2}_I(M)$. Consider
$\alpha^2\in H^{2,2}_I(M)$, where $\alpha$ is the 
generator of $H^{1,1}_I(M)\cap H^2(M, \Z)$. 

\hfill

\sublemma\label{_degrees_theta_alpha^2_Sublemma_}
Let $J$ be an induced complex structure,
$J\circ I=-J\circ I$, and $\deg_I$, $\deg_J$ the degree maps
associated with $I$, $J$. 
Then
\[ \deg_I \alpha^2 >0, \deg_J\alpha^2 =0, \deg_I\theta=\deg_J\theta >0. 
\]
\index{terms}{degree!for different induced complex structures}

{\bf Proof:}
Since $\alpha$ is a K\"ahler class with respect to $I$, we have
$\deg_I \alpha^2 >0$. Since the cohomology class $\theta$ is $SU(2)$-invariant,
and $SU(2)$ acts transitively on the set of induced complex structures, 
we have $\deg_I\theta=\deg_J\theta$. It remains to show that
$\deg_J\alpha^2 =0$ and $\deg_J\theta >0$. The manifold $M$ is
by our assumptions simple; thus, $\dim H^{2,0}(M) =1$
(\cite{_Besse:Einst_Manifo_}). Therefore, in the natural $SU(2)$-invariant
decomposition 
\begin{equation}\label{_H^2_isoty_Equation_} 
    H^2(M) = H_{inv}^2(M) \oplus H^2_+(M),
\end{equation}
we have $\dim H^2_+(M) = 3$. In particular, the intersection 
$H^2_+(M)\cap H^{1,1}_I(M)$ is 1-dimensional. Consider the
decomposition of $\alpha$, associated
with \eqref{_H^2_isoty_Equation_}: $\alpha = \alpha_+ + \alpha_{inv}$.
Since $\alpha$ is of type $(1,1)$ with respect to $I$,
the class $\alpha_+$ is proportional to the K\"ahler class
$\omega_I$, with positive coefficient. A similar argument leads to
the following decomposition for $\theta$:
\[
  \theta = \omega_I^2 + \omega_J^2 +\omega_K^2 + \sum x_i^2,
\]
where $K= I\circ J$ is an induced complex structure, and 
the classes $x_i$ belong to $H_{inv}^2(M)$.
{}From \ref{_p_prods_H^2_Corollary_},
we obtain that the classes $x_i^2$ satisfy
$\deg_I(x_i^2)=0$
(here we use $\dim_{\Bbb H}(M)>1$). 
Thus, 
\[ \deg_I(\theta) = \deg_I(\omega_I^2 + \omega_J^2 +\omega_K^2) =
   \deg_I(\omega_I^2) >0.
\]
Similarly one checks that
\[ \deg_J(\alpha^2) = \deg_J((\alpha_+ + \alpha_{inv})^2) 
   = \deg_J(\alpha_+^2) = \deg_J (c^2 \omega_I) =0.
\]
This proves \ref{_degrees_theta_alpha^2_Sublemma_}.
\endproof

\hfill

Return to the proof of \ref{_C_restri_from_A,d_Lemma_}.
Since $\deg_I\alpha^2 \neq \deg_J \alpha^2$, the class 
$\alpha^2$ is {\it not}
$SU(2)$-invariant. Since $\theta$ {\it is} $SU(2)$-invariant,
$\theta$ is not collinear with $\alpha^2$. 
The degrees $\deg_I$ of $\theta$ and $\alpha^2$ are non-zero;
we have $\deg_I(\theta)=\deg_J(\theta)$,
$\deg_I(\alpha^2) \neq \deg_J(\alpha^2)$ for
$J$ an induced complex structure, $J\neq \pm I$.
By \ref{_degree_isotypic_Proposition_}, no non-trivial 
linear combination of $\theta$, $\alpha^2$ belongs to
$\ker p$. By \ref{_r_H^*_to_H_r_Claim_} (ii), the classes 
$\theta$, $\alpha^2$ generate a 2-dimensional subspace in
$H^4(M, {\Bbb Q})/\ker r$. 

By assumptions of \ref{_C_restri_from_A,d_Lemma_}, 
the group $H^{2,2}_I(M)\cap H^4(M, \Z)/(\ker r)$ has rank 2. Therefore
$\omega$ and $\alpha^2$ generate the space
\[  H^{2,2}_I(M)\cap H^4(M, {\Bbb Q})/(\ker r).
\]
To prove \ref{_C_restri_from_A,d_Lemma_}
it suffices to show that for all integer classes
\[ \beta = a\alpha^2 + b \theta,  \ \  a \in {\Bbb Q}\backslash 0, \ \  
    \deg_I \beta \geq \deg_J \beta,
\]
we have $|\deg_I \beta| >A \cdot d$,
for a constant $A$ depending only on volume, topology and dimension of $M$.
Since $\deg_I \beta \geq |\deg_J \beta|$,
and $\deg_J \alpha^2 =0$ (\ref{_degrees_theta_alpha^2_Sublemma_}),
we have
\[ 
   \deg_I (a \alpha^2 + b \theta) \geq |\deg_J b \theta|.
\]
Therefore, either $a$ and $b$ have the same sign,
 or $\deg_I (a \alpha^2) > 2\deg_I(b\theta)$.
In both cases, 
\begin{equation}\label{_deg_beta_geq_deg_alpha_Equation_}
  |\deg_I \beta|\geq \frac{1}{2}\deg_I (a \alpha^2).
\end{equation}
Let $x\in {\Bbb Q}^{>0}$ be the smallest positive rational value of
$a$ for which there exists an integer class
$\beta = a\alpha^2 + b \theta$. We have an integer
lattice $L_1$ in $H^4_r(M)$ provided by the products of
integer classes; the integer lattice $L_2\supset L_1$ provided by
integer cycles might be different from that one.
Clearly, $x$ is greater than determinant 
$\det (L_1 /L_2)$ of $L_1$ over $L_2$,
and this determinant is determined by the topology of $M$.

Form the definition of $x$ and \eqref{_deg_beta_geq_deg_alpha_Equation_},
we have $|\deg_I \beta| > x^2  \deg_I (\alpha^2)$.
On the other hand, $\deg_I (\alpha^2)> C \deg_I(\alpha)$,
where $C$ is a constant depending on volume and
dimension of $M$. Setting $A:= x^2 \cdot C$, we obtain
$|\deg_I \beta| > x^2 \cdot C \cdot d$. This proves
\ref{_C_restri_from_A,d_Lemma_}. \endproof

\hfill

Consider the maps
\[ \widetilde\deg:\; H^2(M, \Z)\backslash H^2_{inv}(M)
   \arrow \R,
\]
\[
    c:\;  H^2(M)\backslash H^2_{inv}(M) \arrow \c R/\{\pm 1\}
\]
introduced in the beginning of the proof of
\ref{_restri_for_H^11_1-dim-Proposition_}.

\hfill

\lemma \label{_dense_big_tilde_deg_Lemma_}
In assumptions of \ref{_restri_for_H^11_1-dim-Proposition_},
let \[ \calo\subset H^2(M, \R)\backslash H^2_{inv}(M) \]
be an open subset of $H^2(M, \R)$,
such that for all $x\in \calo$, $k\in \R^{>0}$, we have
$k\cdot x \in \calo$. Assume that
$\calo$ contains the K\"ahler class $\omega_I$ for all
induced complex structures $I\in \c R$.
For a positive number $C\in \R^{>0}$,
consider the set $X_C\subset \calo$ 
\[ X_C := \left\{ \vphantom\prod\alpha \in \calo \cap H^2(M, \Z)\;\; | \;\; 
          \widetilde \deg(\alpha) \geq C \right\}.
\]
Then 
$c(X_C)$ is dense in $\c R/\{\pm 1\}$ for all 
$C\in \R^{>0}$.

\hfill

{\bf Proof:}  
The map $\widetilde \deg$ can be expressed in the
following wey. We call an integer cohomology class 
$\alpha\in H^2(M, \Z)$ {\bf indivisible}
if its largest divisor is 1, that is, there are no
integer classes $\alpha'$, and numbers $k\in \Z$, $k>1$, such that
$\alpha = k \alpha'$.

\hfill

\sublemma
 Let $\alpha \in H^2(M)$ be an 
 non-$SU(2)$-invariant cohomology class and
$\alpha= \alpha_{inv} + \alpha_+$ be a decomposition
associated with \eqref{_H^2_isoty_Equation_}.
Assume that $\alpha$ is indivisible.
Then 
\begin{equation}\label{_tilde_deg_Equation_} 
   \widetilde \deg (\alpha) = C\sqrt{((\alpha_+,\alpha_+)_{\c B})}, 
\end{equation}
where $(\cdot,\cdot)_{\c B}$ is the \index{terms}{Bogomolov--Beauville pairing}
Bogomolov-Beauville pairing
on $H^2(M)$ (\cite{_coho_announce_}; see also
\ref{_gene_all_SU(2)_Theorem_}),
and $C$ a constant depending on $\dim M$, $\Vol M$.

{\bf Proof:}
By \ref{_degree_isotypic_Proposition_},
\[ \deg_I(\alpha) = \deg_I(\alpha_+)
\]
(clearly, $p(\alpha) = \alpha_+$). By definition of
$(\cdot,\cdot)_{\c B}$,  we have 
\[ \deg_I(\alpha_+)=(\alpha_+,\omega_{c(\alpha)})_{\c B} \]
On the other hand, $\alpha_+$ is collinear with $\omega_{c(\alpha)}$
by definition of the map $c$. Now \eqref{_tilde_deg_Equation_} 
follows trivially from routine properties of bilinear forms.
\endproof

\hfill

Let $I$ be an induced complex structure
such that the cohomology class $\omega_I$ is irrational:
$\omega_I \notin H^2(M, {\Bbb Q})$. 
\begin{equation}\label{_sequence_x_i_Equation_}
\begin{minipage}[m]{0.8\linewidth}
To prove \ref{_dense_big_tilde_deg_Lemma_},
we have to produce a sequence $x_i \in \calo \cap H^2(M, \Z)$
such that 
\begin{description}
\item[(i)] $c(x_i)$ 
converges to $I$, 
\item[(ii)] and $\lim\tilde \deg (x_i) = \infty$.
\end{description}
\end{minipage}
\end{equation}
We introduce a metric $(\cdot,\cdot)_{\c H}$ on $H^2(M, \R)$,
\[ (\alpha, \beta)_{\c H}:= (\alpha_+,\beta_+)_{\c B} -
   (\alpha_{inv},\beta_{inv})_{\c B}.
\]
It is easy to check that $(\cdot,\cdot)_{\c H}$
is positive definite (\cite{_Verbitsky:cohomo_}).
For every $\epsilon$, there exists
a rational class $\omega_\epsilon \in H^2(M, {\Bbb Q})$
which approximates $\omega_I$ with precision
\[ 
   (\omega_\epsilon - \omega_I, \omega_\epsilon - \omega_I)_{\c H} <\epsilon.
\]
Since $\calo$ is open and contains $\omega_I$, we may assume that
$\omega_\epsilon$ belongs to $\calo$.
Take a sequence $\epsilon_i$ converging to $0$, and let
$\tilde x_i:= \omega_{\epsilon_i}$ be the corresponding 
sequence of rational cohomology cycles. Let 
$x_i:= \lambda_i\tilde x_i$ be the minimal
positive integer such that $x_i \in H^2(M, \Z)$.
We are going to show that the sequence $x_i$
satisfies the conditions of \eqref{_sequence_x_i_Equation_}.
First of all, $\tilde x_i$ converges to $\omega_I$,
and the map
\[
    c:\;  H^2(M)\backslash H^2_{inv}(M) \arrow \c R/\{\pm 1\}
\] 
is continuous. Therefore, $\lim c(\tilde x_i) = c(\omega_I) = I$.
By construction of $c$, $c$ satisfies $c(x) = c(\lambda x)$, and thus,
$c(x_i) = c(\tilde x_i)$. This proves the condition (i) of
\eqref{_sequence_x_i_Equation_}. On the other hand,
since $\omega_I$ is irrational, the sequence $\lambda_i$
goes to infinity. Therefore,
\[ \lim (x_i, x_i)_{\c H} = \infty. \]
It remains to compare $(x_i, x_i)_{\c H}$ with
$\tilde \deg x_i$. By \eqref{_tilde_deg_Equation_},
\[ \tilde \deg x_i = \sqrt{((x_i)_+, (x_i)_+)_{\c B}}. \]
On the other hand, since $(x_i)_+\in H^2_+(M)$, we have
\[ ((x_i)_+, (x_i)_+)_{\c B} = ((x_i)_+, (x_i)_+)_{\c H}.\]
To prove \eqref{_sequence_x_i_Equation_} (ii),
it remains to show that
\[ \lim ((x_i)_+, (x_i)_+)_{\c H} = \lim (x_i, x_i)_{\c H}. \]
Since the cohomology class $\tilde x_i\in H^2(M, {\Bbb Q})$ 
$\epsilon$-approximates $\omega_I$, and $\omega_I$ belongs to
$H^2_+(M)$, we have
\[ 
   (\tilde x_i - (\tilde x_i)_+, \tilde x_i - (\tilde x_i)_+)_{\c H} <\epsilon_i.
\]
Therefore,
\begin{equation}\label{_x_i_close_x_i_+_Equation_} 
  ( x_i - ( x_i)_+,  x_i - ( x_i)_+)_{\c H} <\lambda_i \epsilon_i.
\end{equation}
On the other hand, for $i$ sufficiently big, the cohomology 
class $\tilde x_i$ approaches $\omega_I$, and
\begin{equation}\label{_x_i_bigger_Equation_}
   (x_i, x_i)_{\c H} > 
   \frac 1 2 \lambda_i (\omega_I, \omega_I)_{\c H}
\end{equation}
Comparing \eqref{_x_i_close_x_i_+_Equation_} and
\eqref{_x_i_bigger_Equation_} and using the distance property
for the distance given by $\sqrt{(\cdot,\cdot)_{\c H}}$, we find that 
\begin{equation}\label{_x_i_+_bigger_Equation_}
   \sqrt{(x_i)_+, (x_i)_+} > 
   \sqrt{\frac 1 2 \lambda_i (\omega_I, \omega_I)_{\c H}} 
   - \sqrt{\lambda_i \epsilon_i} = 
     \sqrt{\lambda_i}\cdot
        \left(\sqrt{\frac 1 2 (\omega_I, \omega_I)_{\c H}} -
        \sqrt{\epsilon_i} \right).
\end{equation}
Since $\epsilon_i$ converges to 0 and $\lambda_i$ converges
to infinity, the right hand side of 
\eqref{_x_i_+_bigger_Equation_} converges to infinity.
On the other hand, by 
\eqref{_tilde_deg_Equation_} 
the left hand side of \eqref{_x_i_+_bigger_Equation_}
  is equal constant times $\tilde \deg x_i$, so
$\lim \tilde \deg x_i = \infty$. This proves the second condition of
\eqref{_sequence_x_i_Equation_}.
\ref{_dense_big_tilde_deg_Lemma_} is proven.
\endproof

\hfill

We use \ref{_C_restri_from_A,d_Lemma_} and 
\ref{_dense_big_tilde_deg_Lemma_} in order to finish the proof of
\ref{_restri_for_H^11_1-dim-Proposition_}.

Let $M$ be a compact hyperk\"ahler manifold,
and $\calo\subset H^2(M, \R)$ be the set of all K\"ahler classes
for the K\"ahler metrics compatible with one of induced complex
structures. By \ref{_K_proje_on_R_Lemma_}, $\calo$ is open in
$H^2(M, \R)$. Applying \ref{_dense_big_tilde_deg_Lemma_} to
$\calo$, we obtain the following. In assumptions of 
\ref{_restri_for_H^11_1-dim-Proposition_},
let $Y_C\subset \c R$ be the set of all 
algebraic induced complex structures
$I$ with $\tilde \deg \alpha >C$, where
$\alpha$ is a rational K\"ahler class, $\alpha \in H^{1,1}(M) \cap H^2(M, \Z)$.
Then $Y_C$ is dense in $\c R$. 
Now, \ref{_C_restri_from_A,d_Lemma_},
implies that for all $I \in Y_C$, the
induced complex structure 
$I$ is $A \cdot C$-restricted, where $A$ is the
universal constant of \ref{_C_restri_from_A,d_Lemma_}.
Thus, for all $C$ the set of $C$-restricted induced complex 
structures is dense in $\c R$. This proves that
$M$ admits $C$-restricted complex structures.
We finished the proof of 
\ref{_restri_for_H^11_1-dim-Proposition_}.
\endproof

\subsection{Hyperk\"ahler structures admitting
$C$-restricted complex structures}
\label{_modu_and_C-restri_Subsection_}

Let $M$ be a compact complex manifold admitting a hyperk\"ahler
\index{terms}{hyperk\"ahler structures} structure $\c H$.  Assume that $(M, \c H)$ is a
simple hyperk\"ahler manifold of dimension $\dim_{\Bbb H} M >1$.
The following definition of (coarse, marked)  moduli space for
complex and hyperk\"ahler 
\index{terms}{moduli space!of complex structures!definition of}
\index{terms}{moduli space!of hyperk\"ahler structures!definition of}
\index{terms}{hyperk\"ahler structures!moduli of}
 structures on $M$ is standard.

\hfill

\definition\label{_moduli_hyperka_Definition_}
Let $M_{C^\infty}$ be the $M$ considered as a 
differential manifold,  $\widetilde{Comp}$ be the set
of all integrable complex structures, 
and $\widetilde{\mbox{\it Hyp}}$ be the set of all hyperk\"ahler
structures on $M_{C^\infty}$. The set $\widetilde{\mbox{\it Hyp}}$
is equipped with a natural topology. Let $\widetilde{\mbox{\it Hyp}}^0$
be a connected component of $\widetilde{\mbox{\it Hyp}}$ containing $\c H$
and $\widetilde{Comp}^0$ be a set of all complex structures
$I\in \widetilde{Comp}$ which are compatible with some hyperk\"ahler
structure $\c H_1\in \widetilde{\mbox{\it Hyp}}^0$.
Let $\mbox{\it Diff}$ be the group of diffeomorphisms
of $M$ which act trivially on the cohomology. The
coarse, marked moduli $\mbox{\it Hyp}$ of hyperk\"ahler structures on
$M$ is the quotient 
$\mbox{\it Hyp}:= \widetilde{\mbox{\it Hyp}}^0/\mbox{\it Diff}$ 
equipped with a natural topology. The coarse, marked moduli $Comp$ of
complex structures on $M$ is defined as 
$Comp:= \widetilde{Comp}^0/\mbox{\it Diff}$.
For a detailed discussion of various aspects of this definition,
see \cite{_Verbitsky:cohomo_}.

\hfill

Consider the variety 
\[ X \subset {\Bbb P} H^2(M, \C),\] consisting of all lines 
$l\in {\Bbb P} H^2(M, \C)$ which are isotropic with respect to the
\index{terms}{Bogomolov--Beauville pairing}
Bogomolov-Beauville's pairing:
\[ 
   X:= \{ l \in H^2(M, \C) \; \; | \; (l, l)_{\c B} =0 \}. 
\]
Since $M$ is simple, $\dim H^{2,0}(M, I) =1$ for all induced complex structures.
Let $P_c:\; Comp \arrow {\Bbb P} H^2(M, \C)$ map $I$ to the line
$H^{2,0}_I(M) \subset  H^2(M, \C)$. The map $P_c$ is
called {\bf the period map}.  It is well known that
$Comp$ is equipped with a natural complex structure. From general
properties of the period map it follows that
$P_c$ is compatible with this complex structure.
Clearly from the definition
of \index{terms}{Bogomolov--Beauville pairing}
Bogomolov-Beauville's form, $P_c(I)\in X$ for all induced complex
structures $I\in Comp$ (see \cite{_Beauville_} for details).

\hfill

\theorem \label{_Bogomo_etale_Theorem_}
\cite{_Besse:Einst_Manifo_} (Bogomolov) 
\index{names}{Bogomolov, F. A.}
\index{terms}{Bogomolov's theorem!on period mapping}
The complex analytic
map \[ P_c:\; Comp \arrow X\] is locally an etale covering.
\footnote{
The space $Comp$ is smooth, as follows from
\ref{_Bogomo_etale_Theorem_}. This space
is, however, in most cases not separable 
(\cite{_Huybrechts_}). The space $\mbox{\it Hyp}$ has no
natural complex structures, and can be odd-dimensional.
}

\endproof

\hfill

It is possible to formulate a similar statement about hyperk\"ahler 
\index{terms}{hyperk\"ahler structures!period map for} 
structures.
For a hyperk\"ahler structure $\c H$, consider the set $\c R_{\c H}\subset Comp$
of all induced complex structures associated with this hyperk\"ahler
structure. The subset $\c R_{\c H}\subset Comp$ is a complex analytic
subvariety, which is isomorphic to $\C P^1$. Let $S:= P_c(\c R_{\c H})$
be the corresponding projective line in $X$, and $\bar L(X)$ be the 
space of smooth deformations of $S$ in $X$. The points of $L(X)$ correspond
to smooth rational curves of degree 2 in ${\Bbb P} H^2(M, \C)$. 
For every such curve $s$, there exists a unique 3-dimensional plane 
$L(s) \subset H^2(M, \C)$, such that $s$ is contained in ${\Bbb P} L$.
Let $Gr$ be the Grassmanian manifold of all 3-dimensional
planes in $H^2(M, \C)$ and $Gr_0\subset Gr$ the set of all planes 
$L\in Gr$ such that the restriction of the \index{terms}{Bogomolov--Beauville pairing}
Bogomolov-Beauville
form to $L$ is non-degenerate.
Let $L(X)\subset \bar L(X)$ be the space of all rational curves
$s\in \bar L(X)$ such that the restriction of the \index{terms}{Bogomolov--Beauville pairing}
Bogomolov-Beauville
form to $L(s)$ is non-degenerate: $L(s)\in Gr_0$.
The correspondence $s\arrow L(s)$ gives a map 
$\kappa:\; L(X) \arrow Gr_0$. 

\hfill

\lemma
The map $\kappa:\; L(X) \arrow Gr_0$ is an isomorphism of complex
varieties.

{\bf Proof:} For every plane $L\in Gr_0$, consider the set $s(L)$
of all isotropic lines $l\in L$, that is, lines 
satisfying $(l, l)_{\c B}=0$.
Since $(\cdot,\cdot)_{\c B}\restrict L$ is non-degenerate,
the set $s(L)$ is a rational curve in ${\Bbb P}L$.
Clearly, this curve has degree 2. Therefore,
$s(L)$ belongs to $X(L)$. The map $L\arrow s(L)$
is inverse to $\kappa$. \endproof

\hfill

Consider the standard anticomplex involution 
\[ \iota:\; H^2(M, \C) \arrow H^2(M, \C), \ \ \ \eta\arrow \bar \eta. \]
Clearly, $\iota$
is compatible with the \index{terms}{Bogomolov--Beauville pairing}
Bogomolov-Beauville form.
Therefore, $\iota$ acts on $L(X)$ as an anticomplex involution.
Let $L(X)_\iota\subset L(X)$ be the set of all $S\in L(X)$ fixed
by $\iota$. 

\hfill

Every hyperk\"ahler 
\index{terms}{hyperk\"ahler structures!twistor curves associated with} 
\index{terms}{hyperk\"ahler structures!Bogomolov's theorem for} 
structure \[ \c H\in \mbox{\it Hyp} \]
gives a rational curve $\c R_{\c H}\subset Comp$ with points
corresponding all induced complex structures. Let 
$P_h(\c H)\subset X$ be the line $P_c(\c R_{\c H})$.
Clearly from the definition, 
$P_h(\c H)$ belongs to $L(X)_\iota$. We have constructed
a map $P_h:\; \mbox{\it Hyp} \arrow L(X)_\iota$.
Let $L(Comp)$ be the space of deformations if
$\c R_{\c H}$ in $Comp$. 
Denote by \[ \gamma:\; \mbox{\it Hyp} \arrow L(Comp) \]
the map $\c H \arrow \c R_{\c H}$.
The following result gives a hyperk\"ahler analogue
of \index{names}{Bogomolov, F. A.} \index{terms}{Bogomolov's theorem!on period mapping}
Bogomolov's theorem (\ref{_Bogomo_etale_Theorem_}).

\hfill

\theorem \label{_hyperka_etale_Theorem_}
The map $\gamma:\; \mbox{\it Hyp} \arrow L(Comp)$
is an embedding. The map $P_h:\; \mbox{\it Hyp} \arrow L(X)_\iota$
is locally a covering. 

{\bf Proof:} The first claim is an immediate
consequence of \index{terms}{Calabi--Yau Theorem} Calabi-\index{names}{Yau, S.-T.}Yau Theorem 
(\ref{_symplectic_=>_hyperkahler_Proposition_}). 
Now, \ref{_hyperka_etale_Theorem_}
follows from the \index{terms}{Bogomolov's theorem!on period mapping}
Bogomolov's theorem 
(\ref{_Bogomo_etale_Theorem_}) and dimension count.
\endproof

\hfill

Let $I\in Comp$ be a complex structure on $M$.
Consider the groups
\[ 
   H_h^2(M, I):= H^{1,1}(M, I) \cap H^2(M, \Z) 
\]
and 
\[ 
   H_h^2(M, I):= H^{2,2}_r(M, I) \cap H^4(M, \Z). 
\]
For a general $I$, 
$H^2_h(M, I)=0$ and $H^4_h(M, I)=\Z$ as follows from 
\ref{_theta_SU(2)_inva_Proposition_}. Therefore,
the set of all $I$ with $\rk H^2_h(M, I) =1$, 
$\rk H^4_h(M, I)=2$ is a union of countably 
many subvarieties of codimension 1 in $Comp$.
Similarly, the set $V\subset Comp$ 
of all $I$ with $\rk H^2_h(M, I) >1$, 
$\rk H^4_h(M, I)>2$ is a union of countably 
many subvarieties of codimension more than 1.
Together with \ref{_hyperka_etale_Theorem_},
this implies the following.

\hfill

\claim\label{_hype_1_2_dense_Claim_}
Let $U\subset \mbox{\it Hyp}$ be the set of 
all $\c H\in Hyp$ such that $\c R_{\c H}$ does not intersect $V$.
Then $U$ is dense in $\mbox{\it Hyp}$.

\hfill

{\bf Proof:} Consider a natural involution $i$ of $Comp$
which is compatible with the involution $\iota:\; X \arrow X$ inder 
the period map $P_c:\; Comp \arrow X$. This involution maps
the complex structure $I$ to $-I$. 
\begin{equation}\label{_Hyp_identi_Equation_}
\begin{minipage}[m]{0.8\linewidth}
By \ref{_hyperka_etale_Theorem_},
$\mbox{\it Hyp}$ is identified with an open subset  in 
the set $L(X)_\iota$ of real points of $L(Comp)$. 
\end{minipage}
\end{equation}
Let $L_U\subset L(Comp)$ be the set of all lines
which do not intersect $V$. Since $V$ is a union of
subvarieties of codimension at least 2, a general rational
line $l \in L(Comp)$ does not intersect $V$. Therefore,
$L_U$ is dense in $L(Comp)$. Thus, the set of real points
of $L_U$ is dense $L(X)_\iota$. Using the identification
\eqref{_Hyp_identi_Equation_}, 
we obtain the statement of \ref{_hype_1_2_dense_Claim_}.
\endproof

\hfill

\ref{_hype_1_2_dense_Claim_}
together with \ref{_restri_for_H^11_1-dim-Proposition_}
imply the following theorem.

\hfill

\theorem\label{_C_restri_dense_Theorem_}
Let $M$ be a compact simple hyperk\"ahler manifold,
$\dim_{\Bbb H}M >1$, and $\mbox{\it Hyp}$ its coarse
marked moduli of hyperk\"ahler 
\index{terms}{moduli space!of hyperk\"ahler structures} 
\index{terms}{$C$-restricted complex structures!are dense} 
\index{terms}{hyperk\"ahler structures!admitting $C$-restricted complex structures} 
structures.
Let $U\subset \mbox{\it Hyp}$ be the set of
all hyperk\"ahler structures which admit $C$-restricted complex
structures (\ref{_admitti_C_restri_Definition_}). 
Then $U$ is dense in $\mbox{\it Hyp}$.

\endproof

\subsection{Deformations of coherent sheaves over manifolds with
$C$-res\-t\-ric\-ted complex structures}

The following theorem shows that a 
\index{terms}{coherent sheaves!semistable} 
semistable
deformation of a 
\index{terms}{coherent sheaves!hyperholomorphic!deformations of} 
hyperholomorphic sheaf on $(M,I)$ 
is again hyperholomorphic, provided that $I$ is a
\index{terms}{$C$-restricted complex structures!and coherent sheaves} 
$C$-restricted complex structure and $C$ is sufficiently big.

\hfill

\theorem \label{_sheaf_on_C_restr_hyperho_Theorem_}
Let $M$ be a compact hyperk\"ahler manifold, and
$\c F \in \c F(M)$ a polystable hyperholomorphic sheaf
on $M$ (\ref{_hyperho_shea_on_M_Definition_}).
Let $I$ be a $C$-restricted induced complex structure,
for $C= \deg_I c_2(\c F)$,\footnote{Clearly, since
$\c F$ is hyperholomorphic, the class 
$c_2(\c F)$ is $SU(2)$-invariant, and the number
$\deg_I c_2(\c F)$ independent from $I$.}
and $F'$ be a semistable torsion-free coherent sheaf
on $(M,I)$ with the same rank and Chern classes as 
$\c F$. Then the sheaf $F'$ is hyperholomorphic.

\hfill

{\bf Proof:} Let $F_1$, ..., $F_n$ be the Jordan-H\"older
series for the sheaf $F'$. Since $\c F$ is hyperholomorphic,
we have $\text{slope}(\c F)=0$ (\ref{_slope_hyperho_Remark_}). 
\index{terms}{slope!of a hyperholomorphic sheaf} 
Therefore, $\slope(F_i)=0$, and $\deg_I(c_1(F_i))=0$. 
By \ref{_C_restri_Definition_} (i), then, the class
$c_1(F_i)$ is $SU(2)$ invariant for all $i$.
 To prove that $F'$ is hyperholomorphic
it remains to show that the classes $c_2(F_i)$, $c_2(F_i^{**})$
are $SU(2)$-invariant for all $i$.

\hfill

Consider an exact sequence 
\[ 
   0 \arrow F_i \arrow F^{**}_i \arrow F_i / F_i^{**} \arrow 0. 
\]
Let $[F_i / F_i^{**}]\in H^4(M)$ be the fundamental class
of the union of all components of $Sup(F_i / F_i^{**})$
of complex codimension 2, taken with appropriate multiplicities.
Clearly,
$c_2(F_i)= c_2(F^{**}_i) + [F_i / F_i^{**}]$.
Since $[F_i / F_i^{**}]$ is an effective cycle,
$\deg_I([F_i / F_i^{**}])\geq0$. By the \index{names}{Bogomolov, F. A.} Bogomolov-Miyaoka-\index{names}{Yau, S.-T.}Yau 
inequality (see \ref{_stable_shea_degree_Corollary_}),
we have $\deg_I (c_2(F^{**}_i)\geq 0$. Therefore, 
\begin{equation}\label{_c_2_greater_Equation_}
\deg_I c_2(F_i)\geq \deg_I c_2(F^{**}_i)\geq 0.
\end{equation}

\hfill

Using the product formula for Chern classes, we obtain
\begin{equation}\label{_pro_Chern_2_Equation_}
c_2(F) = \sum_i  c_2(F_i) + \sum_{i, j} c_2(F_i)\wedge c_2(F_j).
\end{equation}
By \ref{_p_prods_H^2_Corollary_},
$\deg_I(\sum_{i, j} c_2(F_i)\wedge c_2(F_j)) =0$.
Since the numbers $\deg_I c_2(F_i)$ are non-negative,
we have $\deg_I c_2(F_i)\leq \deg_I c_2(F) =C$.
By \ref{_funda_and_Chern_CA_Theorem_}, the classes
$c_2(F_i)$, $c_2(F_i^{**})$
are of CA-type. By \ref{_C_restri_Definition_} (ii),
\index{terms}{CA-type (cohomology classes of)} 
then, the inequality $\deg_I c_2(F_i)\leq C$ implies that
the class $c_2(F_i)$ is $SU(2)$-invariant. By 
\eqref{_c_2_greater_Equation_}, 
$\deg_I c_2(F^{**}_i)\leq\deg_I c_2(F_i)$,
so the class $c_2(F_i^{**})$ is also $SU(2)$-invariant.
\ref{_sheaf_on_C_restr_hyperho_Theorem_}
is proven.
\endproof


\section{Desingularization of  hyperholomorphic sheaves}
\label{_desingu_she_Section_}
\index{terms}{coherent sheaves!hyperholomorphic!desingularization of}

The aim of this section is the following theorem.

\hfill

\theorem\label{_desingu_hyperho_Theorem_}
Let $M$ be a hyperk\"ahler manifold, not necessarily
compact, $I$ an induced complex structure,
and $F$ a \index{terms}{coherent sheaves!reflexive}
reflexive coherent sheaf over $(M, I)$
equipped with a \index{terms}{connections!hyperholomorphic} 
hyperholomorphic connection 
(\ref{_hyperholo_co_Definition_}). Assume that
 $F$ has isolated singularities. Let $\tilde M\stackrel \sigma\arrow M$ 
be a blow-up of $(M,I)$ in the singular set of $F$, and
$\sigma^* F$ the pullback of $F$. Then
$\sigma^* F$ is a locally trivial sheaf,
that is, a holomorphic vector bundle.

\hfill

We prove \ref{_desingu_hyperho_Theorem_}
in Subsection \ref{_desingu_she_Subsection_}.

\hfill

The idea of the proof is the following. We apply to $F$ the 
methods used in the proof of Desingularization Theorem
(\ref{_desingu_Theorem_}). The main ingredient in the proof of
Desingularization Theorem is the existence of a natural \index{terms}{$\C^*$-action on hyperk\"ahler manifolds} $\C^*$-action
on the completion $\hat \calo_x(M, I)$ of the local ring $\calo_x(M, I)$, for
all $x\in M$. This $\C^*$-action identifies $\hat \calo_x(M, I)$
with a completion of a graded ring. Here we show that a sheaf
$F$ is \index{terms}{$\C^*$-equivariant sheaves} 
$\C^*$-equivariant. Therefore, a germ of $F$ at $x$ has a
grading, which is compatible with the natural $\C^*$-action on
 $\hat \calo_x(M, I)$. Singularities of such 
\index{terms}{coherent sheaves!reflexive!singularities of}
reflexive sheaves
can be resolved by a single blow-up.

\subsection{Twistor lines and  complexification}
\label{_twi_lines_C^*_Subsection_}

Further on, we need the following definition.

\hfill

\definition
\index{terms}{complexification!definition of}
Let $X$ be a real analytic variety,
which is embedded to a complex variety $X_\C$.
Assume that the sheaf of complex-valued real analytic 
functions on $X$ coincides with the restriction
of $\calo_{X_\C}$ to $X\subset X_\C$. Then
$X_\C$ is called {\bf a complexification of $X$}.

\hfill

For more details on complexification, the reader is referred to
\cite{_GMT_}. There are the most important properties.

\hfill

\claim \label{_complexi_Claim_}
In a neighbourhood of $X$, the manifold $X_\C$ has an anti-complex
involution. The variety $X$ is identified with the set of 
fixed points of this involution, considered as a
real analytic variety.

Let $Y$ be a complex variety, and $X$ the underlying
real analytic variety. Then the product of $Y$ and its complex
conjugate is a \index{terms}{complexification!properties of} 
complexification of $X$, with embedding
$X\hookrightarrow Y \times \bar Y$ given by the diagonal.

The complexification is unique in the following weak sense.
For $X_\C$, $X'_\C$ complexifications of $C$, the complex manifolds
$X_\C$, $X'_\C$ are naturally identified in a neighbourhood
of $X$.

\endproof

\hfill

Let $M$ be a hyperk\"ahler manifold, $\Tw(M)$ its \index{terms}{twistor space} twistor space,
and $\pi:\; \Tw(M) \arrow \C P^1$ the \index{terms}{twistor space} twistor projection.
Let $l\subset \Tw(M)$ be a rational curve, such that
the restriction of $\pi$ to $l$ is an identity. 
Such a curve gives a section of $\pi$, and vice versa, every
section of $\pi$ corresponds to such a curve.
The set of sections of the projection $\pi$ is called
{\bf the space of twistor lines}, denoted by $\Lin$, or $\Lin(M)$.
\index{terms}{twistor lines!space of}
\index{terms}{twistor lines!definition of}
This space is equipped with complex structure,
by \index{names}{Douady, A.} Douady (\cite{_Douady_}). 

Let $m\in M$ be a point. Consider a twistor line 
${s_m}:\; I \arrow (I \times m)\in \C P^1 \times M = \Tw$. 
Then $s_m$ is called {\bf a horizontal twistor
\index{terms}{twistor lines!horizontal!definition of}
line}. The space of horizontal twistor lines is a 
real analytic subvariety in $\Lin$, denoted by $\Hor$, or $\Hor(M)$.
Clearly, the set $\Hor$ is naturally identified with $M$.

\hfill

\proposition\label{_Lin_is_MxM_Proposition_}
(\index{names}{Hitchin, N. J.} Hitchin, \index{names}{Karlhede, A.} Karlhede, Lindstr\"{o}m, Ro\v{c}ek)
Let $M$ be a hyperk\"ahler manifold,
$\Tw(M)$ its \index{terms}{twistor space} twistor space, $I$, $J\in \C P^1$  induced complex
structures, and $\Lin$ the space of \index{terms}{twistor lines} twistor lines.
The complex manifolds $(M,I)$ and $(M, J)$ are naturally
embedded to $\Tw(M)$: 
\[ (M,I) = \pi^{-1}(I),\ \  (M,J) = \pi^{-1}(J).
\] 
Consider a point $s\in \Lin$, $s:\; \C P^1\arrow \Tw(M)$.
Let \[ ev_{I,J}:\; \Lin(M) \arrow (M, I)\times (M,J) \]
be the map defined by $ev_{I,J}(s) = (s(I), s(J))$.
Assume that $I\neq J$. Then there exists a
neighbourhood $U$ of $\Hor\subset \Lin$, such that
the restriction of $ev_{I,J}$ to $U$ is an open embedding.

{\bf Proof:} \cite{_HKLR_}, \cite{_Verbitsky:hypercomple_}.
\endproof

\hfill

Consider the anticomplex involution $i$ of $\C P^1 \cong S^2$
which corresponds to the central symmetry of $S^2$.
Let $\iota:\; \Tw \arrow \Tw$ be the corresponding
involution of the 
\index{terms}{twistor space!involution of} 
twistor space $\Tw(M) = \C P^1\times M$, 
$(x, m) \arrow (i(x), m)$. It is clear that
$\iota$ maps holomorphic subvarieties of $\Tw(M)$ to 
holomorphic subvarieties. Therefore, $\iota$ acts
on $\Lin$ as an anticomplex involution.
For $J= -I$, we obtain a local identification of
$\Lin$ in a neighbourhood of $\Hor$ with 
$(M, I) \times (M, -I)$, that is, with $(M, I)$ times its complex
conjugate. Therefore, the space of 
\index{terms}{twistor lines!space of} 
twistor lines is a 
\index{terms}{complexification!the space of twistor lines as} 
complexification of $(M,I)$. The natural anticomplex involution
of \ref{_complexi_Claim_} coincides with $\iota$. This
gives an identification of $\Hor$ and the real analytic
manifold underlying $(M,I)$.

\hfill

We shall explain how to construct the natural 
\index{terms}{$\C^*$-action on hyperk\"ahler manifolds!construction of} 
$\C^*$-action on
a local ring of a hyperk\"ahler manifold, using the machinery
of \index{terms}{twistor lines} twistor lines.

\hfill

Fix a point $x_0\in M$ and induced complex structures
$I$, $J$, such that $I\neq \pm J$. Let $V_1$, $V_2$ be neighbourhoods of
$s_{x_0}\in \Lin$, and $U_1$, $U_2$ be neighbourhoods of
$(x_0,x_0)$ in $(M, I)\times (M,-I)$, $(M, J)\times (M,-J)$,
such that the evaluation maps $ev_{I, -I}$, $ev_{J, -J}$
induce isomorphisms
\[ ev_{I, -I}:\; V_1 \oldtilde\arrow  U_1, \ \ \ \ 
   ev_{J, -J}:\; V_2\oldtilde\arrow  U_2.
\]
Let $B$ be an open neighbourhood of $x_0\in M$, 
such that $(B, I)\times (B,-I)\subset U_1$
and $(B, I)\times (B,-I)\subset U_2$. Denote by
$V_I\subset V_1$ be the preimage of $(B, I)\times (B,-I)$
under $ev_{I, -I}$, and by 
$V_J\subset V_2$ be the preimage of $(B, J)\times (B,-J)$
under $ev_{J, -J}$. Let $p_I:\; V_I \arrow (B, I)$ be the evaluation,
$s\arrow s(I)$, and $e_I:\; (B, I)\arrow V_1$ the map associating
to $x\in B$ the unique twistor line passing through
$(x, x_0) \subset (B, I)\times (B,-I)$.
In the same fashion, we define $e_J$ and $p_J$.
We are interested in the composition
\[ \Psi_{I, J}:= e_I\circ p_J \circ e_J \circ p_I:\; (B_0, I) \arrow (B,I) \]
which is defined in a smaller neighbourhood 
$B_0\subset B$ of $x_0\in M$.

\hfill

The following proposition is the focal point of this Subsection:
we explain the map $\Psi_{I, J}$ of \cite{_Verbitsky:DesinguII_},
\cite{_Verbitsky:hypercomple_} is geometric terms (in
\cite{_Verbitsky:DesinguII_}, \cite{_Verbitsky:hypercomple_}
this map was defined algebraically).

\hfill

\proposition \label{_Psi_acts_on_TM_Proposition_}
Consider the map $\Psi_{I, J}:\; (B_0, I) \arrow (B,I)$
defined above. By definition, $\Psi_{I, J}$
preserves the point $x_0\in B_0\subset B$. 
Let $d\Psi_{I, J}$ be the differential of $\Psi_{I, J}$ acting 
on the tangent space $T_{x_0}B_0$. Assume that
$I\neq \pm J$. Then 
$d\Psi_{I, J}$ is a multiplication by a scalar
$\lambda\in \C$, $0<|\lambda| <1$.

\hfill

{\bf Proof:} 
The map $\Psi_{I, J}$ was defined in \cite{_Verbitsky:DesinguII_},
\cite{_Verbitsky:hypercomple_} using the identifications
between the real analytic varieties underlying $(M,I)$ and $(M, J)$. 
We proved that $\Psi_{I, J}$ defined this way acts on $T_{x_0}B_0$
as a multiplication by the scalar $\lambda\in \C$, $0<|\lambda| <1$.
It remains to show that the map $\Psi_{I, J}$ defined in  
\cite{_Verbitsky:DesinguII_}, \cite{_Verbitsky:hypercomple_} 
coincides with $\Psi_{I, J}$ defined above.

Consider the natural identification
\[ (B, I)\times (B,-I)\sim (B, J)\times (B,-J), \]
which is defined in a neighbourhood $B_\C$ of $(x_0, x_0)$.
There is a natural projection 
$a_I:\; B_\C \arrow (M,I)$. Consider the embedding
$b_I:\; (B, I) \arrow B_\C$, $x\arrow (x, x_0)$, 
defined in a neighbourhood of $x_0\in (B,I)$.
In a similar way we define $a_J$, $b_J$. 
In \cite{_Verbitsky:DesinguII_},
\cite{_Verbitsky:hypercomple_} we defined 
$\Psi_{I, J}$ as a composition
$b_I\circ a_J \circ b_J \circ a_I$. 
Earlier in this Subsection, we described
a local identification of $(B, I)\times (B,-I)$
and $\Lin(B)$. Clearly, under this identification,
the maps $a_I$, $b_I$ correspond to $p_I$, $e_I$. Therefore,
the definition of $\Psi_{I,J}$ given in this paper
is equivalent to the definition given in
\cite{_Verbitsky:DesinguII_},
\cite{_Verbitsky:hypercomple_}.
\endproof

\subsection[The automorphism $\Psi_{I,J}$ acting on hyperholomorphic
sheaves]{The automorphism $\Psi_{I,J}$ acting on hyperholomorphic \\
sheaves}
\label{_Psi_on_shea_Subsection_}

In this section, we prove that 
\index{terms}{coherent sheaves!hyperholomorphic!$\C^*$-equivariant} 
hyperholomorphic sheaves
are equivariant with respect to the map
$\Psi_{I, J}$, considered as an automorphism of the
local ring $\calo_{x_0}(M,I)$.

\hfill

\theorem \label{_Psi_equiv_hyperho_Theorem_}
Let $M$ be a hyperk\"ahler manifold, not necessarily compact,
$x_0\in M$ a point, 
$I$ an induced complex structure and $F$ a reflexive sheaf
over $(M,I)$ equipped with a \index{terms}{connections!hyperholomorphic} hyperholomorphic connection.
Let $J\neq \pm I$ be another induced complex structure,
and $B_0$, $B$ the neighbourhoods of $x_0\in M$ for which
the map $\Psi_{I,J}:\; B_0 \arrow B$ was defined in 
\ref{_Psi_acts_on_TM_Proposition_}. 
Assume that $\Psi_{I,J}:\; B_0 \arrow B$ is an isomorphism.
Then there exists
a canonical functorial isomorphism of coherent sheaves
\[ \Psi_{I,J}^F:\; F\restrict {B_0} \arrow \Psi_{I,J}^*(F\restrict B).
\]
{\bf Proof:}
Return to the notation introduced in Subsection
\ref{_twi_lines_C^*_Subsection_}. Let $W:= V_I \cap V_J$.
By definition of $V_I$, $V_J$, the evaluation maps
produce open embeddings 
\[ ev_{I, -I}:\; \Lin(W) \hookrightarrow (W,I)\times (W, -I),\]
and 
\[ ev_{J, -J}:\; \Lin(W) \hookrightarrow (W,J)\times (W, -J),\]
Let $S\subset W$ be the singular set of $F\restrict W$,
$\Tw(S) \subset \Tw(W)$ the corresponding embedding, 
and $L_0\subset \Lin(W)$ be the set of all
lines $l\in \Lin(W)$ which do not intersect
$\Tw(S)$. Consider the maps
\[ p_I:\; L_0 \hookrightarrow 
   (W,I)\backslash S 
\]
and 
\[  p_J:\; L_0 \hookrightarrow 
   (W,J)\backslash S 
\]
obtained by restricting the evaluation map
$p_I:\; \Lin(M) \arrow (M, I)$ to $L_0\subset \Lin(M)$.
Since $F$ is equipped with a \index{terms}{connections!hyperholomorphic} hyperholomorphic connection,
the vector bundle $F\restrict{(M, J)\backslash S}$ has a natural
holomorphic structure. Let 
$\underline F_1:= p_I^*\left(F\restrict{(M, I)\backslash S}\right)$
and $\underline F_2:= p_J^*\left(F\restrict{(M, J)\backslash S}\right)$
be the corresponding pullback sheaves over $L_0$,
and $F_1$, $F_2$ the sheaves on $\Lin(W)$ obtained as direct
images of $\underline F_1$, $\underline F_2$ under the open embedding
$L_0 \hookrightarrow \Lin(W)$.

\hfill

\lemma \label{_F_1_=F_2_Lemma_}
Under these assumptions, the sheaves 
$F_1$, $F_2$ are coherent \index{terms}{coherent sheaves!reflexive}
reflexive sheaves. Moreover,
there exists a natural isomorphism of coherent sheaves
$\Psi_{1,2}:\;F_1\arrow F_2$.

\hfill

{\bf Proof:}
The complex codimension of the singular set $S$ in $(M,I)$
is at least 3, because $F$ is reflexive
(\cite{_OSS_}, Ch. II, 1.1.10). Since $S$
is \index{terms}{trianalytic subvarieties} trianalytic (\ref{_singu_triana_Claim_}),
this codimension is even. Thus,
$\codim_\C (S, (M,I)) \geq 4$. Therefore,
\[ \codim_\C (\Tw(S), \Tw(M)) \geq 4. \]
Consider the set $L_S$ of all twistor lines 
$l\in \Lin(W)$ passing through $\Tw(S)$.
For generic points $x,y\in \Tw(W)$, there exists
a unique line $l\in \Lin(W)$ passing through $x, y$. Therefore,
\[ \codim_\C (L_S, \Lin(W)) = \codim_\C (\Tw(S), \Tw(M))-1\geq 3. \]
By definition, $L_0:= \Lin(W)\backslash L_S$. Since
$F_1$, $F_2$ are direct images of bundles 
$\underline F_1$, $\underline F_2$ over a subvariety
$L_S$ of codimension 3, these sheaves are coherent and 
\index{terms}{coherent sheaves!reflexive}
reflexive
(\cite{_OSS_}, Ch. II, 1.1.12; see also \ref{_normal_refle_Lemma_}).
To show that they are naturally isomorphic it remains
to construct an isomorphism between $\underline F_1$
and $\underline F_2$.

Let $\c F$ be a coherent sheaf on $\Tw(W)$ obtained 
from $F\restrict W$ as in the proof of
\ref{_conne_=>_hyperho_Proposition_}.
The singular set of $\c F$ is $\Tw(S)\subset \Tw(W)$.
Therefore, the restriction
$\c F\restrict{\Tw(W)\backslash \Tw(S)}$ is a holomorphic vector bundle.
For all horizontal twistor lines $l_x\subset \Tw(W)\backslash \Tw(S)$,
the restriction $\c F\restrict {l_x}$ is clearly a trivial vector
bundle over $l_x\cong \C P^1$. A small deformation of a trivial
vector bundle is again trivial. Shrinking
$W$ if necessary, we may assume that for all
lines $l\in L_0$, the restriction of $\c F$ to $l\cong \C P^1$
is a trivial vector bundle.

The isomorphism 
$\underline \Psi_{1,2}:\;\underline F_1\arrow\underline F_2$ 
is constructed as follows. Let $l\in L_0$ be a \index{terms}{twistor lines} twistor line. The restriction
$\c F\restrict {l}$ is trivial. Consider $l$ as a map
$l:\; \C P^1 \arrow \Tw(M)$. We identify $\C P^1$ with the set
of induced complex structures on $M$. By definition,
the fiber of $F_1$ in $l$
is naturally identified with the space $\c F\restrict {l(I)}$,
and the fiber of $F_2$ in $l$ is identifies with 
$\c F\restrict {l(J)}$. Since $\c F\restrict {l}$ is trivial,
the fibers of the bundle $\c F\restrict {l}$ are naturally
identified. This provides a vector bundle isomorphism
$\underline \Psi_{1,2}:\;\underline F_1\arrow\underline F_2$
mapping $\underline F_1 \restrict l = \c F\restrict {l(I)}$
to $\underline F_2 \restrict l = \c F\restrict {l(J)}$.
It remains to show that this isomorphism is compatible
with the holomorphic structure. 

Since the bundle $\c F\restrict l$ is trivial, we have
an identification 
\[ \c F\restrict {l(I)}\cong\c F\restrict{l(J)}= \Gamma(\c F\restrict l),
\]
where $\Gamma(\c F\restrict l)$ is the space of global sections
of $\c F\restrict l$.
Thus, $F_i \restrict l = \Gamma(\c F\restrict l)$, and this
identification is clearly holomorphic. This proves
\ref{_F_1_=F_2_Lemma_}. \endproof

\hfill

We return to the proof of \ref{_Psi_equiv_hyperho_Theorem_}.
Denote by $F_J$ the restriction of $\c F$ to
$(M, J)= \pi^{-1}(J) \subset \Tw(M)$.
The map $\Psi_{I, J}$ was defined as a composition
$e_I\circ p_J \circ e_J \circ p_I$.
The sheaf $p_I^* F$ is by definition isomorphic to
$F_1$, and $p_J^* F_J$ to $F_2$. On the other hand, clearly,
$e_J^* F_2= F_J$. Therefore,
$(p_J \circ e_J)^* F_2 \cong F_2$.
Using the isomorphism $F_1\cong F_2$, we obtain
$(p_J \circ e_J)^* F_1 \cong F_1$.
To sum it up, we have the following isomorphisms:
\begin{equation*}
\begin{split}
p_I^* F &\cong F_1,\\
(p_J \circ e_J)^* F_1 &\cong F_1,\\
e_I^* F_1 &\cong F.
\end{split}
\end{equation*}
A composition of these isomorphisms gives
an isomorphism
\[ \Psi_{I,J}^F:\; F\restrict {B_0} \arrow \Psi_{I,J}^*(F\restrict B).
\]
This proves \ref{_Psi_equiv_hyperho_Theorem_}.
\endproof

\subsection{A $\C^*$-action on a local ring of a hyperk\"ahler manifold}
\label{_C^*_action_on_loca_ring_Subsection_}

Let $M$ be a hyperk\"ahler manifold, non necessarily compact, 
$x\in M$ a point and
$I$, $J$ induced complex structures, $I\neq J$. 
Consider the complete local ring 
$\calo_{x, I}:= \hat \calo_x(M, I)$.
Throughout this section we consider the map $\Psi_{I,J}$ 
(\ref{_Psi_acts_on_TM_Proposition_})
as an automorphism of the ring $\calo_{x, I}$.
Let $\goth m$ be the maximal ideal of
$\calo_{x, I}$, and $\goth m/\goth m^2$ 
the Zariski cotangent space of $(M,I)$ in $x$.
\begin{equation} \label{_Psi_acts_on_cota_Equation_}
\begin{minipage}[m]{0.8\linewidth}
By \ref{_Psi_acts_on_TM_Proposition_},
$\Psi_{I,J}$ acts on $\goth m/\goth m^2$ 
as a multiplication by a number $\lambda\in\C$,
$0<|\lambda|<1$.
\end{minipage}
\end{equation}
Let $V_{\lambda^n}$ be the eigenspace corresponding to 
the eigenvalue $\lambda^n$,
\[ V_{\lambda^n}:= \{ v\in \calo_{x, I}\ \ 
   |\ \ \Psi_{I,J}(v) = \lambda^n v\}.
\]
Clearly, $\oplus V_{\lambda^i}$ is a graded subring in $\calo_{x, I}$.
In \cite{_Verbitsky:DesinguII_},
(see also \cite{_Verbitsky:hypercomple_}) we proved that
the ring $\oplus V_{\lambda^i}$ is dense in $\calo_{x, I}$
with respect to the adic topology. Therefore, the ring
$\calo_{x, I}$ is identified with the adic completion of
$\oplus V_{\lambda^i}$.

\hfill

Consider at action of $\C^*$ on $\oplus V_{\lambda^i}$, with
$z\in \C^*$ acting on $V_{\lambda^i}$ as a multiplication by
$z^i$. This \index{terms}{$\C^*$-action on hyperk\"ahler manifolds} 
$\C^*$-action is clearly continuous, with respect to
the adic topology.  Therefore, it can be extended to
\[\calo_{x, I} = \widehat {\oplus V_{\lambda^i}}.\]

\hfill

\definition \label{_Psi(z)_Definition_}
Let $M$, $I$, $J$, $x$, $\calo_{x, I}$ be as in the beginning
of this Subsection. Consider the $\C^*$-action 
\[ \Psi_{I,J}(z):\; \calo_{x, I} \arrow \calo_{x, I} \]
constructed as above. Then $\Psi_{I,J}(z)$ is called
{\bf the canonical $\C^*$-action associated with
$M$, $I$, $J$, $x$.} 

\hfill

In the above notation, consider a \index{terms}{coherent sheaves!reflexive}
reflexive sheaf $F$ on $(M,I)$
equipped with a \index{terms}{connections!hyperholomorphic} hyperholomorphic connection. Denote the germ of
$F$ at $x$ by $F_x$, $F_x:= F\otimes_{\calo_{(M,I)}} \calo_{x, I}$. 
From \ref{_Psi_equiv_hyperho_Theorem_}, we obtain an isomorphism
$F_x \cong \Psi^*_{I,J} F_x$. This isomorphism can be interpreted
as an automorphism
\[ \Psi^F_{I,J}:\; F_x\arrow F_x\]
satisfying
\begin{equation}\label{_Psi^F_and_multi_Equation_}
 \Psi^F_{I,J}(\alpha v) = \Psi_{I,J}(\alpha) v, 
\end{equation}
for all $\alpha \in \calo_{x, I}$, $v\in F_x$.

By \eqref{_Psi^F_and_multi_Equation_}, 
the automorphism $\Psi^F_{I,J}$ respects the filtration
\[ F_x \supset \goth m F_x \supset \goth m^2 F_x \supset ... \]
Thus, it makes sense to speak of $\Psi^F_{I,J}$-action
on $\goth m^i F_x /\goth m^{i+1} F_x$.

\hfill

\lemma\label{_Psi^F_on_m^iF/m^i+1F_Lemma_}
The automorphism $\Psi^F_{I,J}$ acts on 
$\goth m^i F_x /\goth m^{i+1} F_x$
as a multiplication by $\lambda^i$, where
$\lambda\in \C$ is the number considered in 
\eqref{_Psi_acts_on_cota_Equation_}.

\hfill

{\bf Proof:} By \eqref{_Psi^F_and_multi_Equation_},
it suffices to prove \ref{_Psi^F_on_m^iF/m^i+1F_Lemma_}
for $i=0$. In other words, we have to show that
$\Psi^F_{I,J}$ acts as identity on $F_x /\goth m F_x$.
We reduced \ref{_Psi^F_on_m^iF/m^i+1F_Lemma_} to the following
claim.

\hfill

\claim \label{_Psi^F_identi_on_F/mF_Claim_}
In the above assumptions, 
the automorphism $\Psi^F_{I,J}$ 
acts as identity on $F_x /\goth m F_x$.

\hfill

{\bf Proof:} 
In the course of defining the map
 $\Psi^F_{I,J}$, we identified the space $\Lin(M)$ with
a \index{terms}{complexification} complexification of $(M,I)$, and defined 
the maps 
\[ p_I:\; \Lin(M) \arrow (M, I),\ \ 
    p_J:\; \Lin(M) \arrow (M, I)
\] 
(these maps are smooth, in a neighbourhood of $\Hor\subset \Lin(M)$,
by \ref{_Lin_is_MxM_Proposition_}),
and 
\[ e_I:\; (B, I) \arrow \Lin(M),\ \ 
   e_J:\; (B, J) \arrow \Lin(M)
\] 
(these maps are locally 
closed embeddings). Consider $(M,J)$ as a subvariety
of $\Tw(M)$, $(M, J) = \pi^{-1}(J)$.
Let $\c F$ be the lift of $F$ to $\Tw(M)$ 
(see the proof of \ref{_conne_=>_hyperho_Proposition_} 
for details). Denote the completion of
$\calo_x(M, J)$ by $\calo_{x, J}$.
Let $F_J$ denote the $\calo_{x, J}$-module
$\left(\c F\restrict{(M,J)} \right)\otimes_{\calo_{(M,J)}} \hat \calo_{x, J}$.
Consider the horizontal \index{terms}{twistor lines} twistor line
$l_x\in \Lin(M)$. Let $\Lin_x(M)$ be the spectre of the
completion $\calo_{x,\Lin}$ of the local ring of
holomorphic functions on $\Lin(M)$ in $l_x$.
The maps $p_I$, $p_J$, $e_I$, $e_J$ can be considered
as maps of corresponding formal manifolds:
\begin{equation*}
\begin{split}
p_I:\; \Lin_x(M) & \arrow \Spec(\calo_{x, I}),\\
p_J:\; \Lin_x(M) & \arrow \Spec(\calo_{x, J}),\\
e_I:\; \Spec(\calo_{x, I})& \arrow\Lin_x(M),\\
e_J:\; \Spec(\calo_{x, J})& \arrow\Lin_x(M),
\end{split}
\end{equation*}
As in Subsection \ref{_Psi_on_shea_Subsection_},
we consider the $\calo_{x,\Lin}$-modules
$F_1:= p_I^* F_x$ and $F_2:= p_J^* F_J$. By \ref{_F_1_=F_2_Lemma_},
there exists a natural isomorphism $\Psi_{1,2}:\; F_1\arrow F_2$.

Let $\goth m_{l_x}$ be the maximal ideal of $\calo_{x,\Lin}$.
Since the morphism $p_I$ is smooth, the space 
$F_1/ \goth m_{l_x} F_1$ is naturally isomorphic to
$F_x/\goth m F_x$. Similarly, the space
$F_2/ \goth m_{l_x} F_2$ is isomorphic to
$F_J / \goth m_J F_J$, where $\goth m_J$ is the
maximal ideal of $\calo_{x, J}$. We have a chain of isomorphisms
\begin{equation} \label{_chain_iso_Psi_F/mF_Equation_}
\begin{split}
   F_x/\goth m F_x &\stackrel{p_I^*}\arrow F_1/ \goth m_{l_x} F_1
   \stackrel {\Psi_{1,2}} \arrow F_2/ \goth m_{l_x} F_2\\
   &\stackrel{e_J^*}\arrow F_J/ \goth m_{J} F_J
   \stackrel{p_J^*}\arrow  F_2/ \goth m_{l_x} F_2\\
   &\stackrel {\Psi_{1,2}^{-1}} \arrow F_1/ \goth m_{l_x} F_1
   \stackrel{e_I^*}\arrow F_x/ \goth m F_x.
\end{split}
\end{equation}
By definition, for any $f\in F_x/\goth m F_x$,
the value of $\Psi^F_{I,J}(f)$ is given by
the composide map of \eqref{_chain_iso_Psi_F/mF_Equation_} 
applied to $f$. The composition 
\begin{equation} \label{_restri_F/mF_F_2_and_back_Equation_}
    F_2/ \goth m_{l_x} F_2
   \stackrel{e_J^*}\arrow F_J/ \goth m_{J} F_J
   \stackrel{p_J^*}\arrow  F_2/ \goth m_{l_x} F_2
\end{equation}
is identity, because the spaces
$F_2/ \goth m_{l_x} F_2$ and 
$F_J / \goth m_J F_J$ are canonically identified, and this
identification can be performed via 
$e_J^*$ or $p_J^*$. Thus,
the map \eqref{_chain_iso_Psi_F/mF_Equation_} is a composition
\[
   F_x/\goth m F_x \stackrel{p_I^*}\arrow F_1/ \goth m_{l_x} F_1
   \stackrel {\Psi_{1,2}} \arrow F_2/ \goth m_{l_x} F_2
   \stackrel {\Psi_{1,2}^{-1}} \arrow F_1/ \goth m_{l_x} F_1
   \stackrel{e_I^*}\arrow F_x/ \goth m F_x.
\]
This map is clearly equivalent to a composition
\[  F_x/\goth m F_x \stackrel{p_I^*}\arrow F_1/ \goth m_{l_x} F_1
    \stackrel{e_I^*}\arrow F_x/ \goth m F_x, 
\]
which is identity according to the same reasoning which proved that
\eqref{_restri_F/mF_F_2_and_back_Equation_} is identity. 
We proved \ref{_Psi^F_identi_on_F/mF_Claim_}
and \ref{_Psi^F_on_m^iF/m^i+1F_Lemma_}. \endproof

\hfill

Consider the $\lambda^n$-eigenspaces $F_{\lambda^n}$ of $F_x$.
Consider the $\oplus V_{\lambda^n}$-submodule 
$\oplus F_{\lambda^n}\subset F_x$, where 
$\oplus V_{\lambda^n}\subset \calo_{x, I}$ is
the ring defined in Subsection 
\ref{_C^*_action_on_loca_ring_Subsection_}.
From \ref{_Psi^F_identi_on_F/mF_Claim_} and
\eqref{_Psi_acts_on_cota_Equation_} it follows that
$\oplus F_{\lambda^n}$ is dense in 
$F_x$, with respect to the adic topology on
$F_x$. For $z\in \C^*$, let
$\Psi_{I,J}^F(z):\; \oplus F_{\lambda^n}\arrow \oplus F_{\lambda^n}$
act on $F_{\lambda^n}$ as a multiplication by $z^n$.
As in \ref{_Psi(z)_Definition_}, we extend 
$\Psi_{I,J}^F(z)$ to $F_x = \widehat{\oplus F_{\lambda^n}}$.
This automorphism makes $F_x$ into a 
\index{terms}{$\C^*$-equivariant sheaves} 
$\C^*$-equivariant
module over $\calo_{x, I}$

\hfill

\definition\label{_C^*_stru_on_sge_Definition_}
The constructed above $\C^*$-equivariant
\index{terms}{$\C^*$-equivariant sheaves!construction of}  
structure
on $F_x$ is called {\bf the canonical $\C^*$-equivariant structure
on $F_x$ associated with $J$}.

\subsection{Desingularization of $\C^*$-equivariant sheaves}
\label{_desingu_she_Subsection_}

Let $M$ be a hyperk\"ahler manifold, $I$ an induced complex structure
and $F$ a \index{terms}{coherent sheaves!reflexive}
reflexive sheaf with isolated singularities
over $(M,I)$, equipped with a
\index{terms}{connections!hyperholomorphic} hyperholomorphic connection. We have shown that the sheaf $F$ admits
a $\C^*$-equivariant
\index{terms}{$\C^*$-equivariant structure}  
structure compatible
with the canonical $\C^*$-action on the local ring of $(M,I)$. 
Therefore, \ref{_desingu_hyperho_Theorem_}
is implied by the following proposition.

\hfill

\proposition \label{_desingu_C^*_equi_Proposition_}
Let $B$ be a complex manifold, $x\in B$ a point. 
Assume that there is an action $\Psi(z)$ of $\C^*$ on $B$
which fixes $x$ and acts on $T_x B$ be dilatations. Let
$F$ be a \index{terms}{coherent sheaves!reflexive}
reflexive coherent sheaf on $B$, which is 
locally trivial outside of $x$.
 Assume that the germ
$F_x$ of $F$ in $x$ is equipped with a $\C^*$-equivariant
\index{terms}{$\C^*$-equivariant sheaves!desingularization of}   
structure,
compatible with $\Psi(z)$.
Let $\tilde B$ be a blow-up of $B$
in $x$, and $\pi:\; \tilde B\arrow B$ the
standard projection. Then the pullback
sheaf $\tilde F:= \pi^* F$ is locally trivial on $\tilde B$.

\hfill

{\bf Proof:} Let $C:= \pi^{-1}(x)$ be the singular
locus of $\pi$. The sheaf $F$ is locally trivial
outside of $x$. Let $d$ be the rank of
$F\restrict{B\backslash x}$. To prove that
$\tilde F$ is locally trivial, we need to show that
for all points $y\in \tilde B$, the fiber
$\pi^* F\restrict y$ is $d$-dimensional. Therefore,
to prove \ref{_desingu_C^*_equi_Proposition_}
it suffices to show that $\pi^* F\restrict C$
is a vector bundle of dimension $d$.

 The variety $C$ is naturally identified
with the projectivization ${\Bbb P} T_x B$ of the
tangent space $T_xB$. The total space of
$T_x B$ is equipped with a natural action of $\C^*$,
acting by dilatations. Clearly, coherent sheaves on 
${\Bbb P} T_x B$ are in one-to-one 
correspondence with $\C^*$-equivariant
\index{terms}{$\C^*$-equivariant sheaves}   
coherent sheaves on $T_x B$. 
Consider a local isomorphism 
$\phi:\; T_x B \arrow  B$ which is compatible
with \index{terms}{$\C^*$-action on hyperk\"ahler manifolds} 
$\C^*$-action, maps $0\in T_x B$ to $x$ and acts as  
identity on the tangent space $T_0 (T_x B)=T_x B$. 
The sheaf $\phi^* F$ is $\C^*$-equivariant.
Clearly, the corresponding sheaf on
${\Bbb P} T_x B$ is canonically isomorphic with
$\pi^*F \restrict C$. Let $l\in T_x B$ be a line passing
through $0$, and $l\backslash 0$ its complement to $0$.
Denote the corresponding point of ${\Bbb P} T_x B$ by $y$.
The restriction $\phi^* F\restrict {l\backslash 0}$
is a $\C^*$-equivariant vector bundle.
The $\C^*$-equivariant structure
identifies all the fibers of the bundle 
$\phi^* F\restrict {l\backslash 0}$.
Let $F_l$ be one of these fibers.
Clearly, the fiber of $\pi^*F \restrict C$ in $y$
is canonically isomorphic to $F_l$. Therefore,
the fiber of $\pi^*F \restrict C$ in $y$ is $d$-dimensional.
We proved that $\pi^*F$ is a bundle. \endproof


\section{Twistor transform and quaternionic-K\"ahler geometry}
\label{_twisto_tra_Section_}


This Section is a compilation of results known from the literature.
Subsection \ref{_dire_inve_twi_Subsection_} is based on
\cite{_NHYM_} and the results of 
Subsection \ref{_twi_tra_Hermi_Subsection_}
are implicit in  \cite{_NHYM_}.
Subsection \ref{_B_2_bundles_Subsection_} is based
on \cite{_Salamon_}, \cite{_Nitta:bundles_}
and \cite{_Nitta:Y-M_}, and Subsection
\ref{_specia_and_q-K-Subsection_} is a recapitulation
of the results of A. \index{names}{Swann, A.} Swann (\cite{_Swann_}).

\subsection{Direct and inverse twistor transform}
\label{_dire_inve_twi_Subsection_}

In this Subsection, we recall the definition and the main properties
of the direct and inverse twistor \index{terms}{twistor transform} transform for bundles over
hyperk\"ahler manifolds (\cite{_NHYM_}). 

\hfill

The following definition is a non-Hermitian analogue of the
notion of a \index{terms}{connections!hyperholomorphic} hyperholomorphic connection.

\hfill

\definition
Let $M$ be a hyperk\"ahler manifold,
not necessarily compact, and
$(B, \nabla)$ be a vector bundle with a connection
over $M$, not necessarily Hermitian. Assume that
the curvature of $\nabla$ is contained in the space
$\Lambda^2_{inv}(M, \End(B))$ of $SU(2)$-invariant
2-forms with coefficients in $\End(B)$. Then
$(B, \nabla)$ is called {\bf an autodual bundle},
\index{terms}{autodual bundles}
\index{terms}{connections!autodual}
and $\nabla$ {\bf an autodual connection}. 

\hfill

Let
$\Tw(M)$ be the \index{terms}{twistor space} twistor space of $M$, equipped with the standard
maps $\pi:\; \Tw(M) \arrow \C P^1$, $\sigma:\; \Tw(M) \arrow M$.

\hfill

\index{terms}{autodual bundles!direct and inverse twistor transform for}
We introduce the direct and inverse twistor 
\index{terms}{twistor transform} transforms which relate autodual
bundles on the hyperk\"ahler manifold $M$ and 
holomorphic bundles on its twistor space $\Tw(M)$.

\hfill

Let $B$ be a complex vector bundle on $M$ equipped with a connection
$\nabla$. The pullback $\sigma^*B$ of $B$ to $\Tw(M)$ is equipped with a
pullback connection $\sigma^*\nabla$.

\hfill

\lemma \label{_autodua_(1,1)-on-twi_Lemma_}
(\cite{_NHYM_}, Lemma 5.1)
The connection $\nabla$ is autodual if and only if the connection
$\sigma^*\nabla$ has curvature of Hodge type $(1,1)$. 

{\bf Proof:} Follows from \ref{_SU(2)_inva_type_p,p_Lemma_}. \endproof

\hfill

 In assumptions of
\ref{_autodua_(1,1)-on-twi_Lemma_},
consider the $(0,1)$-part $(\sigma^*\nabla)^{0,1}$ 
of the connection $\sigma^*\nabla$. Since
$\sigma^*\nabla$ has curvature of Hodge type $(1,1)$,
we have \[ \left((\sigma^*\nabla)^{0,1}\right)^2=0, \] and by 
\ref{_Newle_Nie_for_NH_bu_Proposition_}, this connection
is integrable. Consider  $(\sigma^*\nabla)^{0,1}$
as a holomorphic structure operator on $\sigma^* B$.

\hfill

\index{terms}{autodual bundles!direct and inverse twistor transform for}
Let $\c A$ be the category of autodual bundles on
$M$, and $\c C$ the category of holomorphic vector
bundles on $\Tw(M)$. We have constructed a functor
\[ (\sigma^* \bullet)^{0,1}:\; \c A \arrow \c C, \]
$\nabla \arrow (\sigma^*\nabla)^{0,1}$.
Let $s\in \Hor\subset \Tw(M)$ be a horizontal twistor 
line (Subsection \ref{_twi_lines_C^*_Subsection_}).
For any $(B, \nabla)\in \c A$, 
consider corresponding holomorphic vector 
bundle $(\sigma^* B, (\sigma^*\nabla)^{0,1})$.
The restriction of $(\sigma^* B, (\sigma^*\nabla)^{0,1})$
to $s\cong \C P^1$ is a trivial vector bundle.
A converse statement is also true. Denote by
$\c C_0$ the category of holomorphic vector bundles
$C$ on $\Tw(M)$, such that the restriction of
$C$ to any horizontal \index{terms}{twistor lines} twistor line is trivial.

\hfill

\theorem \label{_dire_inve_twisto_Theorem_}
Consider the functor
\[ (\sigma^* \bullet)^{0,1}:\; \c A \arrow \c C_0\]
constructed above. Then it is an equivalence of categories.

{\bf Proof:} \cite{_NHYM_}, Theorem 5.12. \endproof

\hfill

\definition
Let $M$ be a hyperk\"ahler manifold,
$\Tw(M)$ its twistor space and $\c F$ a holomorphic vector
bundle. We say that $\c F$ is {\bf compatible with twistor
\index{terms}{vector bundles!compatible with twistor transform} 
\index{terms}{twistor transform!vector bundles compatible with}
transform} if the restriction of
$C$ to any horizontal \index{terms}{twistor lines} twistor line 
$s\in \Tw(M)$ is a trivial bundle on $s\cong \C P^1$.

\hfill

Recall that a connection $\nabla$ in a vector bundle over
a complex manifold is called $(1,1)$-connection if its
curvature is of Hodge type $(1,1)$.

\hfill

\remark \label{_cano_conne_Remark_}
Let $\c F$ be a holomorphic bundle over $\Tw(M)$
which is compatible with twistor transform. Then $\c F$
is equipped with a natural $(1,1)$-connection 
$\nabla_{\c F}= \sigma^* \nabla$,
where $(B, \nabla)$ is the corresponding autodual bundle over $M$.
The connection $\nabla_{\c F}$ is not, generally speaking,
Hermitian, or compatible with a Hermitian structure.
\index{terms}{autodual bundles!Hermitian structures on}
\index{terms}{vector bundles!compatible with twistor transform} 

\subsection{Twistor transform and Hermitian structures
on vector bundles}
\label{_twi_tra_Hermi_Subsection_}

Results of this Subsection were implicit in \cite{_NHYM_},
but in this presentation, they are new.

Let $M$ be a hyperk\"ahler manifold, not necessarily compact,
and $\Tw(M)$ its \index{terms}{twistor space} twistor space.
In Subsection \ref{_dire_inve_twi_Subsection_}, we have shown that certain
holomorphic vector bundles over $\Tw(M)$ admit a canonical
(1,1)-connection $\nabla{_\c F}$ (\ref{_cano_conne_Remark_}).
This connection can be non-Hermitian. Here we study
the Hermitian structures on $(\c F, \nabla{_\c F})$
in terms of holomorphic properties of $\c F$.

\hfill

\definition 
Let $F$ be a real analytic complex vector bundle over a 
real analytic manifold $X_\R$, and $h:\; F\times F \arrow \C$ a
$\calo_{X_\R}$-linear pairing on $F$. Then $h$ is called
{\bf semilinear} if for all $\alpha \in \calo_{X_\R}\otimes_\R \C$,
we have
\[ h(\alpha x, y) = \alpha\cdot h(x, y), \text{\ \ and\ \ }
   h(x,\alpha y) = \bar \alpha\cdot h(x, y). 
\]
For $X$ a complex manifold and $F$ a holomorphic vector bundle,
by a semilinear pairing on $F$ we understand a semilinear pairing
on the underlying real analytic bundle.
Clearly, a real analytic Hermitian metric is always semilinear.

\hfill

Let $X$ be a complex manifold, $I:\; TX \arrow TX$
the complex structure operator, and $i:\; X \arrow X$ a
real analytic map. We say that $X$ is anticomplex if
the induced morphism of tangent spaces satisfies
$i\circ I = - I \circ I$. For a complex
vector bundle $F$ on $X$, consider the complex
adjoint vector bundle $\bar F$, which coincides with
$F$ as a real vector bundle, with $\C$-action which is
conjugate to that defined on $F$. Clearly, for every
holomorphic vector bundle $F$, and any anticomplex map
$i:\; X \arrow X$, the bundle
$i^* \bar F$ is equipped with a natural
holomorphic structure.

\hfill

Let $M$ be a hyperk\"ahler manifold, and $\Tw(M)$ its \index{terms}{twistor space} twistor space.
Recall that $\Tw(M)=\C P^1\times M$ 
is equipped with a canonical anticomplex involution $\iota$,
which acts as identity on $M$ and as central symmetry
$I\arrow -I$ on $\C P^1= S^2$.
For any holomorphic bundle $\c F$ on $\Tw(M)$, consider the corresponding
holomorphic bundle $\iota^* \bar{\c F}$.

\hfill

\index{terms}{autodual bundles!direct and inverse twistor transform for}
Let $M$ be a hyperk\"ahler manifold, $I$ an induced complex
structure, $F$ a vector bundle over $M$, equipped with an autodual
connection $\nabla$, and $\c F$ the corresponding holomorphic
vector bundle over $\Tw(M)$, equipped with a canonical connection
$\nabla_{\c F}$. 
As usually, we identify $(M, I)$ and the fiber $\pi^{-1}(I)$ 
of the \index{terms}{twistor projection} 
twistor projection $\pi:\; \Tw(M) \arrow \C P^1$.
Let $\nabla_{\c F} = \nabla^{1,0}_I + \nabla^{0,1}_I$
be the Hodge decomposition of $\nabla$ with respect to $I$.
\begin{equation} \label{_nabla^1,0_as_holo_Equation_}
\begin{minipage}[m]{0.8\linewidth}
Clearly, the operator $\nabla^{1,0}_I$ can be considered
as a holomorphic structure operator on $F$, considered as
a complex vector bundle over $(M, -I)$.
\end{minipage}
\end{equation}
Then the holomorphic structure operator on $\c F\restrict{(M,I)}$
is equal to $\nabla^{0,1}_I$, and the 
holomorphic structure operator on $\c F\restrict{(M,-I)}$
is equal to $\nabla^{1,0}_I$. 

Assume that 
the bundle $(\c F, \nabla_{\c F})$ is equipped with a
non-degenerate semilinear pairing $h$ which is compatible with
the connection. Consider the natural connection $\nabla_{\c F^*}$
on the dual bundle to $\c F$, and its Hodge decomposition
(with respect to $I$)
\[ \nabla_{\c F^*} = \nabla^{1,0}_{\c F^*} + \nabla^{0,1}_{\c F^*}. \]
 Clearly, the pairing $h$ gives a $C^\infty$-isomorphism
of $\c F$ and the complex conjugate of its dual bundle, denoted as
$\bar{\c F}^*$. Since $h$ is semilinear and
compatible with the connection, it maps the holomorphic structure
operator $\nabla^{0,1}_I$ to the complex conjugate of 
$\nabla^{1,0}_{\c F^*}$. On the other hand, 
the operator $\nabla^{1,0}_{\c F^*}$ is a 
holomorphic structure operator in 
$\c F^*\restrict{(M,-I)}$, as 
\eqref{_nabla^1,0_as_holo_Equation_} claims.
We obtain that the map $h$ can be considered as
an isomorphism of holomorphic vector bundles
\[ h:\; \c F \arrow (\iota^* \bar{\c F})^*.\]
This correspondence should be thought of as a 
(direct) twistor
\index{terms}{twistor transform!for for bundles with a semilinear pairing} 
transform for bundles with a semilinear pairing.

\hfill

\proposition\label{_twi_tra_for_semili_Proposition_}
(direct and inverse twistor
transform for bundles with semilinear pairing)
Let $M$ be a hyperk\"ahler manifold, and $\c C_{sl}$
the category of autodual bundles over $M$ equipped with
\index{terms}{autodual bundles!Hermitian structures on}
a non-degenerate semilinear pairing. Consider the category
$\c C_{hol,sl}$ of holomorphic vector bundles $\c F$
on $\Tw(M)$, compatible with twistor transform
and equipped with an isomorphism
\[ {\goth h}:\; \c F \arrow (\iota^* \bar{\c F})^*.\]
Let $\c T:\; \c C_{sl}\arrow \c C_{hol,sl}$ be the functor
constructed above. Then $\c T$ is an isomorphism of categories.

\hfill

{\bf Proof:} Given a pair 
$\c F, {\goth h}:\; \c F \arrow (\iota^* \bar{\c F})^*$,
we need to construct a non-degenerate semilinear pairing
$h$ on $\c F\restrict{(M,I)}$, compatible with a connection. 
Since $\c F$ is compatible with twistor transform, it
\index{terms}{vector bundles!compatible with twistor transform} 
is a pullback of a bundle $(F, \nabla)$ on $M$.
This identifies the real analytic bundles $\c F\restrict{(M,I')}$,
for all induced complex structures $I'$. Taking $I' = \pm I$,
we obtain an identification of the $C^{\infty}$-bundles
$\c F\restrict{(M,I)}$, $\c F\restrict{(M,-I)}$. Thus,
$\goth h$ can be considered as an isomorphism of
$F= \c F\restrict{(M,I)}$ and $(\bar {\c F})^*\restrict{(M,I)}$.
This allows one to consider $\goth h$ as a semilinear form $h$ on 
$F$. We need only to show that $h$ is compatible with the connection
$\nabla$. Since $\nabla_{\c F}$ is an invariant of holomorphic
structure, the map ${\goth h}:\; \c F \arrow (\iota^* \bar{\c F})^*$
is compatible with the connection $\nabla_{\c F}$. Thus, the obtained
above form $h$ is compatible with 
the connection $\nabla_{\c F}\restrict{(M,I)}=\nabla$.
This proves \ref{_twi_tra_for_semili_Proposition_}.
\endproof

\subsection{ $B_2$-bundles on quaternionic-K\"ahler manifolds}
\label{_B_2_bundles_Subsection_}

\definition \label{_q-K_Definition_}
(\cite{_Salamon_}, \cite{_Besse:Einst_Manifo_})
Let $M$ be a Riemannian manifold. Consider a
bundle of algebras $\End(TM)$, where $TM$ is the tangent
bundle to $M$. Assume that $\End(TM)$ contains a 
4-dimensional bundle of subalgebras $W\subset \End(TM)$,
with fibers isomorphic to a quaternion algebra 
${\Bbb H}$. Assume, moreover, that
$W$ is closed under the transposition map
$\bot:\; \End(TM)\arrow \End(TM)$ and
is preserved by the Levi-Civita connection.
Then $M$ is called {\bf quaternionic-K\"ahler}.
\index{terms}{quaternionic-K\"ahler manifolds!definition of} 

\hfill

\example
Consider the quaternionic \index{terms}{quaternionic projective space}
projective space
\[ {\Bbb H} P^n= ({\Bbb H}^n \backslash 0) / {\Bbb H}^*. \]
It is easy to see that ${\Bbb H} P^n$
is a quaternionic-K\"ahler manifold.
For more examples of quaternionic-K\"ahler manifolds,
see \cite{_Besse:Einst_Manifo_}.

\hfill

\index{terms}{quaternionic-K\"ahler manifolds!is Einstein} 
A quaternionic-K\"ahler manifold is Einstein
(\cite{_Besse:Einst_Manifo_}),
i. e. its Ricci tensor is proportional to the
metric: $Ric(M) = c \cdot g$, with $c\in \R$.
When $c=0$, the manifold $M$ is hyperk\"ahler,
and its restricted \index{terms}{holonomy group} 
holonomy group is $Sp(n)$;
otherwise, the restricted holonomy is
$Sp(n)\cdot Sp(1)$. The number $c$ is called
{\bf the scalar curvature} of $M$. Further on, we 
shall use the term {\it quaternionic-K\"ahler manifold}
for manifolds with non-zero scalar curvature.

The quaternionic \index{terms}{quaternionic projective space}
projective space
${\Bbb H} P^n$ has positive scalar curvature.

\hfill

The quaternionic \index{terms}{quaternionic projective space}
projective space is the only example of
\index{terms}{quaternionic-K\"ahler manifolds} 
quaternionic-K\"ahler manifold which we need, in the course of this
paper. However, the formalism of quaternionic-K\"ahler manifolds 
is very beautiful and significantly 
simplifies the arguments, so we state the definitions and
results for a general quaternionic-K\"ahler manifold 
whenever possible.

\hfill

Let $M$ be a quaternionic-K\"ahler manifold,
and $W\subset \End(TM)$ the corresponding 4-dimensional
bundle. For $x\in M$, consider the set $\c R_x\subset W\restrict x$,
consisting of all $I\in W\restrict x$ satisfying
$I^2=-1$. Consider $\c R_x$ as a Riemannian 
submanifold of the total space of $W\restrict x$. Clearly,
$\c R_x$ is isomorphic to a 2-dimensional sphere. 
Let $\c R= \cup_x \c R_x$ be the corresponding spherical fibration
over $M$, and $\Tw(M)$ its total space. The manifold
$\Tw(M)$ is equipped with an almost complex structure,
which is defined in the same way as the almost complex
structure for the twistor space of a hyperk\"ahler manifold.
This almost complex structure is known to be integrable
(see \cite{_Salamon_}).

\hfill

\index{terms}{quaternionic-K\"ahler manifolds!twistors of} 
\definition\label{_twi_q-K_Definition_}
(\cite{_Salamon_}, \cite{_Besse:Einst_Manifo_})
Let $M$ be a quaternionic-K\"ahler manifold.
Consider the complex manifold $\Tw(M)$ constructed
above. Then $\Tw(M)$ is called {\bf the 
\index{terms}{twistor space!of a quaternionic-K\"ahler manifold} 
twistor space of 
$M$.}

\hfill

Note that (unlike in the hyperk\"ahler case) the space
$\Tw(M)$ is K\"ahler. For 
\index{terms}{quaternionic-K\"ahler manifolds!of positive scalar curvature} 
quaternionic-K\"ahler manifolds
with positive scalar curvature, the anticanonical bundle
of $\Tw(M)$ is ample, so $\Tw(M)$ is a Fano manifold.

\hfill

Quaternionic-K\"ahler analogue of a twistor transform
was studied by T. \index{names}{Nitta, T.} Nitta in a serie of papers
(\cite{_Nitta:bundles_}, \cite{_Nitta:Y-M_} etc.)
It turns out that the picture given in
\cite{_NHYM_} for K\"ahler manifolds is
very similar to that observed by T. \index{names}{Nitta, T.} Nitta.

\hfill

A role of $SU(2)$-invariant 2-forms is played
by the so-called $B_2$-forms.

\hfill

\definition
Let $SO(TM)\subset \End(TM)$ be a group bundle of all 
orthogonal automorphisms of $TM$, and
$G_M:= W\cap SO(TM)$. Clearly, the fibers of $G_M$ 
are isomorphic to $SU(2)$. Consider the action of
$G_M$ on the bundle of 2-forms $\Lambda^2(M)$.
Let $\Lambda^2_{inv}(M)\subset \Lambda^2(M)$
be the bundle of $G_M$-invariants. The bundle 
$\Lambda^2_{inv}(M)$ is called {\bf the bundle 
of \index{terms}{$B_2$-forms} $B_2$-forms}. In a similar fashion we define
$B_2$-forms with coefficients in a bundle.

\hfill

\definition\label{_B_2_bu_Definition_}
In the above assumptions, let $(B, \nabla)$ be a bundle
with connection over $M$. The bundle $B$ is called
{\bf a \index{terms}{$B_2$-bundles!definition of} $B_2$-bundle}, and 
$\nabla$ is called {\bf a $B_2$-connection},
 if its curvature is a $B_2$-form.

\hfill

Consider the natural projection $\sigma:\; \Tw(M)\arrow M$.
The proof of the following claim is completely analogous
to the proof of \ref{_SU(2)_inva_type_p,p_Lemma_} and
\ref{_autodua_(1,1)-on-twi_Lemma_}.

\hfill

\claim\label{_B_2_=_holo_on_Tw_Claim_}
\begin{description}
\item[(i)] Let $\omega$ be a $2$-form on $M$.
The pullback $\sigma^* \omega$ is of type
$(1,1)$ on $\Tw(M)$ if and only if $\omega$
is a \index{terms}{$B_2$-forms} \index{terms}{$B_2$-bundles}
$B_2$-form on $M$.

\item[(ii)]
Let $B$ be a complex vector bundle on $M$ equipped with a 
connection $\nabla$, not necessarily Hermitian. 
The pullback $\sigma^*B$ of $B$ to $\Tw(M)$ is equipped 
with a pullback connection $\sigma^*\nabla$.
Then, $\nabla$ is a \index{terms}{$B_2$-connections} 
$B_2$-connection if and only if 
$\sigma^*\nabla$ has curvature of Hodge type $(1,1)$. 
\end{description}

\endproof

There exists an analogue of direct and inverse twistor
transform as well.
\index{terms}{twistor transform!for $B_2$-bundles}

\hfill

\theorem \label{_dire_inve_q-K_Theorem_}
For any \index{terms}{$B_2$-bundles!twistor transform for} 
$B_2$-connection $(B, \nabla)$, 
consider the corresponding holomorphic vector 
bundle \[ (\sigma^* B, (\sigma^*\nabla)^{0,1}). \]
The restriction of $(\sigma^* B, (\sigma^*\nabla)^{0,1})$
to a line $\sigma^{-1}(m)\cong \C P^1$ is a trivial vector bundle,
for any point $m\in M$.
Denote by $\c C_0$ the category of holomorphic vector bundles
$C$ on $\Tw(M)$, such that the restriction of
$C$ to $\sigma^{-1}(m)$  is trivial, for all $m\in M$,
and by $\c A$ the category of $B_2$-bundles (not necessarily
Hermitian). Consider the functor
\[ (\sigma^* \bullet)^{0,1}:\; \c A \arrow \c C_0\]
constructed above. Then it is an equivalence of categories.

{\bf Proof:} It is easy to modify
the proof of the direct and inverse twistor transform
theorem from \cite{_NHYM_} to work in 
\index{terms}{quaternionic-K\"ahler manifolds!twistor transform for} 
quaternionic-K\"ahler 
situation. \endproof

\hfill

We will not use \ref{_dire_inve_q-K_Theorem_},
except for its consequence, which was proven 
in \cite{_Nitta:bundles_}.

\hfill

\corollary \label{_twi_tra_q-K_Hermi_Corollary_}
Consider the functor
\[ (\sigma^* \bullet)^{0,1}:\; \c A \arrow \c C_0\]
constructed in \ref{_dire_inve_q-K_Theorem_}.
Then $(\sigma^* \bullet)^{0,1}$ gives an injection $\kappa$
from the set of equivalence classes of Hermitian 
$B_2$-connections over $M$ to
the set of equivalence classes of holomorphic connections
over $\Tw(M)$.

\endproof

\hfill

Let $M$ be a \index{terms}{quaternionic-K\"ahler manifolds} quaternionic-K\"ahler manifold. The space $\Tw(M)$
has a natural K\"ahler metric $g$,
such that the standard map $\sigma:\; \Tw(M)\arrow M$ is a Riemannian
submersion, and the restriction of $g$ to the fibers $\sigma^{-1}(m)$ 
of $\sigma$ is a metric of constant curvature on $\sigma^{-1}(m)= \C P^1$ 
 (\cite{_Salamon_}, \cite{_Besse:Einst_Manifo_}). 

\hfill

\example\label{_Tw_HP^n_Fu-St_Example_}
In the case
\index{terms}{twistor space!of the quaternionic projective space} 
$M= {\Bbb H}P^n$, we have $\Tw(M) = \C P^{2n+1}$, and the K\"ahler
metric $g$ is proportional to the Fubini-Study metric on
$\C P^{2n+1}$.

\hfill

\theorem \label{_twi_tra_YM_q-K_Theorem_}
(T. \index{names}{Nitta, T.} Nitta)
Let $M$ be a 
\index{terms}{quaternionic-K\"ahler manifolds!of positive scalar curvature} 
\index{terms}{twistor space!of a quaternionic-K\"ahler manifold} 
quaternionic-K\"ahler manifold of positive scalar curvature,
$\Tw(M)$ its  twistor space, equipped with
a natural K\"ahler structure,
and $B$ a Hermitian 
\index{terms}{$B_2$-bundles!and Yang--Mills connections} 
$B_2$-bundle on $M$.
Consider the pullback $\sigma^* B$, equipped with
a Hermitian connection. Then $\sigma^* B$ 
is a Yang-Mills bundle on $\Tw(M)$, and
$\deg c_1(\sigma^* B) =0$.
\index{terms}{degree!of a bundle obtained as a twistor transform}

{\bf Proof:} \cite{_Nitta:Y-M_}. \endproof

\hfill

Let $\kappa$ be the map considered in
\ref{_twi_tra_q-K_Hermi_Corollary_}. Assume that
$M$ is a compact manifold. In \cite{_Nitta:Y-M_}, T. 
\index{names}{Nitta, T.} Nitta
defined the moduli space of Hermitian $B_2$-bundles. By
\index{terms}{$B_2$-bundles!moduli of}  
\index{terms}{moduli space!of Hermitian $B_2$-bundles} 
\index{names}{Uhlenbeck, K.}
\index{names}{Yau, S.-T.}
\index{terms}{vector bundles!Yang--Mills}
\index{terms}{vector bundles!polystable}
 Uhlenbeck-Yau theorem, Yang-Mills bundles
are polystable. Then the map $\kappa$ provides an 
embedding from the moduli of non-decomposable
Hermitian $B_2$-bundles to the moduli $\c M$
of stable bundles on $\Tw(M)$. The image of
$\kappa$ is a totally real subvariety in $\c M$
(\cite{_Nitta:Y-M_}).

\subsection{Hyperk\"ahler manifolds with special
${\Bbb H}^*$-action and qua\-ter\-ni\-o\-nic-\-K\"ah\-ler manifolds
of positive scalar curvature}
\label{_specia_and_q-K-Subsection_}

Further on, we shall need the following definition.

\hfill

\definition
{\bf An almost hypercomplex manifold}
is a smooth manifold $M$ with an action of quaternion
algebra in its tangent bundle
For each $L\in \Bbb H$, $L^2 = -1$, $L$ gives 
an almost complex structure on $M$. 
The manifold $M$ is caled {\bf hypercomplex}
if the almost complex structure $L$ is integrable,
for all possible choices $L\in \Bbb H$.

\hfill

The \index{terms}{twistor space!of a hypercomplex manifold} 
twistor space for a hypercomplex manifold is defined in the 
same way as for hyperk\"ahler manifolds. It is also a
complex manifold (\cite{_Kaledin_}). The formalism of direct
and inverse twistor 
\index{terms}{twistor transform!for hypercomplex manifolds} 
transform can be repeated for
hypercomplex manifolds verbatim.

\hfill

Let ${\Bbb H}^*$ be the group of non-zero quaternions.
Consider an embedding $SU(2)\hookrightarrow {\Bbb H^*}$.
Clearly, every quaternion $h\in {\Bbb H}^*$
can be uniquely represented as $h= |h| \cdot g_h$,
where $g_h\in SU(2)\subset {\Bbb H}^*$. This gives a 
natural decomposition ${\Bbb H}^*= SU(2)\times {\Bbb R}^{>0}$.
Recall that $SU(2)$ acts naturally on the set of induced complex
structures on a hyperk\"ahler manifold.

\hfill

\definition\label{_H^*_specia_Definition_}
Let $M$ be a hyperk\"ahler manifold equipped with a free 
smooth action $\rho$ of the group ${\Bbb H}^*= SU(2)\times \R^{>0}$. 
The action $\rho$ is called {\bf special}
if the following conditions hold.

\begin{description}
\item[(i)] The subgroup $SU(2)\subset {\Bbb H}^*$
acts on $M$ by isometries.

\item[(ii)] For $\lambda\in {\Bbb R}^{>0}$, the corresponding
action $\rho(\lambda):\; M\arrow M$ is compatible
with the hyperholomorphic structure (which is a fancy way
of saying that $\rho(\lambda)$ is holomorphic with respect to any
of induced complex structures).

\item[(iii)] Consider the smooth
${\Bbb H^*}$-action $\rho_e:\; {\Bbb H^*}\times \End(TM) \arrow \End(TM)$ 
induced on $\End(TM)$ by $\rho$. For any $x\in M$ and
any induced complex structure $I$, the restriction
$I\restrict x$ can be considered as a point in the total
space of $\End(TM)$.
Then, for all induced complex structures $I$,
all $g\in  SU(2)\subset {\Bbb H^*}$, and all $x\in M$,
the map $\rho_e(g)$ maps $I\restrict{x}$ to 
$g(I)\restrict{\rho_e(g)(x)}$.

Speaking informally, this can be stated as
``${\Bbb H}^*$-action interchanges the induced complex
structures''.

\item[(iv)]
Consider the 
automorphism of $S^2 T^*M$ induced by $\rho(\lambda)$, where
$\lambda\in {\Bbb R}^{>0}$. Then $\rho(\lambda)$
maps the Riemannian metric tensor
$s\in S^2 T^*M$ to $\lambda^2 s$. 
\end{description}

\hfill

\example
Consider the flat hyperk\"ahler manifold 
$M_\fl= {\Bbb H}^n \backslash 0$. 
There is a natural action of ${\Bbb H^*}$
on ${\Bbb H}^n\backslash 0$. This gives a special action of
${\Bbb H^*}$ on $M_\fl$.

\hfill

The case of a flat manifold $M_\fl= {\Bbb H}^n \backslash 0$
is the only case where we apply the results of this
section. However, the general statements are just as 
difficult to prove, and much easier to comprehend.

\hfill

\definition\label{_speci_equi_Definition_}
Let $M$ be a hyperk\"ahler manifold with a special action $\rho$ of
${\Bbb H^*}$. Assume that $\rho(-1)$ acts non-trivially on $M$.
Then $M/\rho(\pm 1)$ is also a hyperk\"ahler manifold with
a special action of ${\Bbb H^*}$. We say that the manifolds
$(M, \rho)$ and $(M/\rho(\pm 1), \rho)$ are 
 {\bf hyperk\"ahler manifolds with special action of
${\Bbb H^*}$ which are special equivalent}.
Denote by $H_{sp}$ the category
of hyperk\"ahler manifolds with a special action of
${\Bbb H^*}$ defined up to special equivalence.

\hfill

A. \index{names}{Swann, A.} Swann (\cite{_Swann_}) developed an equivalence between the category of
qua\-ter\-ni\-o\-nic-\-K\"ah\-ler manifolds of positive scalar curvature and
the category $H_{sp}$. The purpose of this Subsection is to
give an exposition of \index{terms}{Swann's formalism} Swann's formalism.

\hfill

Let $Q$ be a 
\index{terms}{quaternionic-K\"ahler manifolds!holonomy group of} 
quaternionic-K\"ahler manifold. The restricted holonomy
\index{terms}{holonomy group} 
group of $Q$ is $Sp(n)\cdot Sp(1)$, that is, 
$(Sp(n)\times Sp(1))/\{\pm 1\}$. Consider the principal bundle $\c G$
with the fiber $Sp(1)/\{\pm 1\}$,
corresponding to the subgroup 
\[ Sp(1)/\{\pm 1\}\subset (Sp(n)\times Sp(1))/\{\pm 1\}.\]
of the holonomy. There is a natural $Sp(1)/\{\pm 1\}$-action
on the space \\ ${\Bbb H^*}/\{\pm 1\}$. Let 
\[ \c U(Q):= \c G\times_{Sp(1)/\{\pm 1\}}{\Bbb H^*}/\{\pm 1\}.\]
Clearly, $\c U(Q)$ is fibered over $Q$, with fibers
which are isomorphic to \\ ${\Bbb H^*}/\{\pm 1\}$.
We are going to show that the manifold $\c U(Q)$ is equipped
with a natural hypercomplex structure.

\hfill

There is a natural smooth decomposition
$\c U(Q)\cong \c G \times {\Bbb R}^{>0}$
which comes from the isomorphism 
${\Bbb H^*}\cong Sp(1)\times {\Bbb R}^{>0}$.

\hfill

Consider the standard 4-dimensional bundle
$W$ on $Q$. Let $x\in Q$ be a point. The fiber
$W\restrict q$ is isomorphic to $\Bbb H$, in a non-canonical way.
The choices of isomorphism $W\restrict q\cong \Bbb H$
are called {\bf quaternion frames in $q$}.
The set of quaternion frames gives a fibration
over $Q$, with a fiber $\Aut({\Bbb H})\cong Sp(1)/\{\pm 1\}$.
Clearly, this fibration coincides with the principal
bundle $\c G$ constructed above. Since 
$\c U(Q)\cong \c G \times {\Bbb R}^{>0}$,
a choice of $u\in \c U(Q)\restrict q$
determines an isomorphism $W\restrict q\cong \Bbb H$.

\hfill

Let $(q, u)$ be the point of $\c U(Q)$,
with $q\in Q$, $u\in \c U(Q)\restrict q$.
The natural connection in $\c U(Q)$
gives a decomposition 
\[ T_{(q, u)}U(Q) = T_u \bigg(\c U(Q)\restrict q\bigg) \oplus T_q Q. \]
The space $\c U(Q)\restrict q \cong {\Bbb H^*}/\{\pm 1\}$
is equipped with a natural hypercomplex structure.
This gives a quaternion action on $T_u \bigg(\c U(Q)\restrict q\bigg)$
The choice of $u\in \c U(Q)\restrict q$ determines a
quaternion action on $T_q Q$, as we have seen above.
We obtain that the total space of $\c U(Q)$ is 
an almost hypercomplex manifold.

\hfill

\proposition\label{_U(Q)_hypercomple_Proposition_}
(A. \index{names}{Swann, A.} Swann) 
Let $Q$ be a 
\index{terms}{quaternionic-K\"ahler manifolds!hyperk\"ahler manifolds associated with} 
quaternionic-K\"ahler manifold.
Consider the manifold $\c U(Q)$ 
constructed as above, and equipped with a
quaternion algebra action in its tangent space. 
Then $\c U(Q)$ is a hypercomplex manifold.

\hfill

{\bf Proof:} Clearly, the manifold $\c U(Q)$
is equipped with a ${\Bbb H}^*$-action,
which is related with the almost hypercomplex
structure as prescribed by \ref{_H^*_specia_Definition_} 
(ii)-(iii). Pick an induced complex structure 
$I\in {\Bbb H}$. This gives an
algebra embedding $\C \arrow {\Bbb H}$.
Consider the corresponding
\index{terms}{$\C^*$-action on twistor spaces} 
$\C^*$-action $\rho_I$ on an almost complex manifold
$(\c U(Q), I)$. This $\C^*$-action is compatible
with the almost complex structure. The quotient
$\c U(Q)/\rho(I)$ is an almost complex manifold, which
is naturally isomorphic to the twistor space $\Tw(Q)$.
Let $L^*$ be a complex vector bundle of all $(1,0)$-vectors
$v\in T (\Tw(Q))$ tangent to the fibers of the standard
projection $\sigma:\; \Tw(Q)\arrow Q$, and $L$
be the dual vector bundle. Denote by 
$Tot_{\neq 0}(L)$ the complement
$\Tot(L)\backslash N$, where $N=\Tw(Q)\subset \Tot(L)$
is the zero section of $L$. 
Using the natural connection in $L$, we obtain
an almost complex structure on $\Tot(L)$.

Consider the natural
projection $\phi:\; Tot_{\neq 0}(L)\arrow Q$.
The fibers $\phi^{-1}(q)$ of $\phi$ are identified with the
space of non-zero vectors in the total space
of the cotangent bundle
$T^* \sigma^{-1}(q)\cong T^*(\C P^1)$. This space is naturally
isomorphic to 
\[ \c G\restrict q \times \R^{>0}= 
   \c U(Q)\restrict q\cong {\Bbb H}^*/\{\pm 1\}.
\]
This gives a canonical isomorphism of almost complex manifolds
\[ (\c U(Q), I)\arrow Tot_{\neq 0}(L).\] Therefore, to prove
that $(\c U(Q), I)$ is a complex manifold, it suffices to
show that the natural almost complex structure on
$Tot_{\neq 0}(L)\subset \Tot(L)$ is integrable. 
Consider the natural connection
$\nabla_L$ on $L$. To prove that $\Tot(L)$
is a complex manifold, it suffices to show that
$\nabla_L$ is a holomorphic connection. 
The bunlde $L$ is known under the name of
{\bf holomorphic contact bundle}, and
it is known to be holomorphic
(\cite{_Salamon_}, \cite{_Besse:Einst_Manifo_}).
\endproof

\hfill

\remark
The result of \ref{_U(Q)_hypercomple_Proposition_}
is well known. We have given its proof because
we shall need the natural identification 
$Tot_{\neq 0}(L)\cong \c U(Q)$ further on
in this paper. 

\hfill

\theorem\label{_U(Q)_hyperk_Theorem_}
Let $Q$ be a 
\index{terms}{quaternionic-K\"ahler manifolds!of positive scalar curvature} 
quaternionic-K\"ahler manifold of positive scalar curvature,
and $\c U(Q)$ the hypercomplex manifold constructed above.
Then $\c U(Q)$ admits a unique (up to a scaling)
hyperk\"ahler metric compatible with the 
hypercomplex structure.

{\bf Proof:} \cite{_Swann_}. \endproof

\hfill

Consider the action of ${\Bbb H}^*$ on
 $\c U(M)$ defined in the 
proof of \ref{_U(Q)_hypercomple_Proposition_}. 
This action satisfies the conditions
(ii) and (iii) of \ref{_H^*_specia_Definition_}.
The conditions (i) and (iv) of \ref{_H^*_specia_Definition_}
are easy to check (see \cite{_Swann_} for details).
This gives a functor from the category $\c C$ of 
\index{terms}{quaternionic-K\"ahler manifolds!hyperk\"ahler manifolds associated with} 
quaternionic-K\"ahler
manifolds of positive scalar curvature to
the category $H_{sp}$ of \ref{_speci_equi_Definition_}.

\hfill

\theorem \label{_U(Q)_equiva_cate_Theorem_}
The functor $Q\arrow \c U(Q)$ from
$\c C$ to $H_{sp}$ is an equivalence of categories. 

{\bf Proof:} \cite{_Swann_}. \endproof

\hfill

The inverse functor from $H_{sp}$ to $C$
is constructed by taking a quotient of $M$ by the action
of ${\Bbb H}^*$. Using the technique of quaternionic-K\"ahler
reduction anf hyperk\"ahler potentials 
(\cite{_Swann_}), one can equip the quotient
$M/{\Bbb H}^*$ with a natural quaternionic-K\"ahler structure.
\index{terms}{quaternionic-K\"ahler manifolds!hyperk\"ahler manifolds associated with}


\section{$\C^*$-equivariant twistor spaces}
\label{_C_equiv_twi_spa_Section_}


\index{terms}{twistor transform}
\index{terms}{$B_2$-bundles}
\index{names}{Swann, A.}
\index{terms}{Swann's formalism!for vector bundles}
\index{terms}{Swann's formalism}

In Section \ref{_twisto_tra_Section_}, we gave an exposition of 
the twistor transform,  $B_2$-bundles and  Swann's formalism.
In the present Section, we give a synthesis of these theories,
obtaining a construction with should be thought
of as Swann's formalism for vector bundles.

Consider the equivalence of categories $Q\arrow \c U(Q)$ 
constructed in \ref{_U(Q)_equiva_cate_Theorem_} 
(we call this equivalence ``Swann's formalism'').
We show that 
\index{terms}{$B_2$-bundles!and $\C^*$-equivariant holomorphic bundles} 
$B_2$-bundles on $Q$ are in functorial
bijective correspondence with
\index{terms}{$\C^*$-equivariant holomorphic bundles on twistor spaces} 
$\C^*$-equivariant holomorphic bundles on $\Tw(\c U(Q))$
(\ref{_B_2_to_C^*_equiva_Theorem_}). 

In Subsection
\ref{_hyperho_shea_C^*_equiv_Y-M_on_blow-up_Subsection_}, 
this equivalence is applied 
to the vector bundle $\pi^*(F)$ of
\ref{_desingu_hyperho_Theorem_}. We use it
to construct a canonical $SU(2)$-equivariant connection
on $\pi^*(F)\restrict C$, where $C$ is 
a special fiber of $\pi:\; \tilde M\arrow (M, I)$
(see \ref{_desingu_hyperho_Theorem_} for details 
and notation). This implies that the holomorphic
bundle $\pi^*(F)\restrict C$ is a direct sum of stable bundles.

\subsection[$B_2$-bundles on quaternionic-K\"ahler manifolds
and $\C^*$-equi\-va\-ri\-ant holomorphic bundles over twistor
spaces]{$B_2$-bundles on quaternionic-K\"ahler manifolds
and \\$\C^*$-equi\-va\-ri\-ant holomorphic bundles over twistor
spaces}
\label{_B_2_to_C^*-invaholo_ove_twi_Subsection_}

For the duration of this Subsection, we fix a hyperk\"ahler manifold
$M$, equipped with a special ${\Bbb H}^*$-action $\rho$, and the corresponding
\index{terms}{quaternionic-K\"ahler manifolds} 
quaternionic-K\"ahler manifold $Q= M/{\Bbb H}^*$. Denote the natural
quotient map by $\phi:\; M\arrow Q$.

\hfill

\lemma \label{_phi^*_B_2-forms_1,1_Lemma_}
Let $\omega$ be a 2-form
over $Q$, and $\phi^* \omega$ its pullback to $M$.
Then the following conditions are equivalent
\begin{description}
\item[(i)] $\omega$ is a \index{terms}{$B_2$-forms} $B_2$-form
\item[(ii)] $\phi^* \omega$ is of Hodge type $(1,1)$
with respect to some induced complex structure $I$ on $M$
\item[(iii)] $\phi^* \omega$ is $SU(2)$-invariant.

\end{description}

{\bf Proof:} 
Let $I$ be an induced complex structure on $M$.
As we have shown in the proof of \ref{_U(Q)_hypercomple_Proposition_},
the complex manifold $(M, I)$ is identified with an open subset 
of the total space $\Tot(L)$ of a holomorphic line bundle $L$
over $\Tw(Q)$. The map $\phi$ is represented as a composition
of the projections $h:\; \Tot(L)\arrow \Tw(Q)$ and
$\sigma_Q:\; \Tw(Q) \arrow Q$. Since the map $h$ is
smooth and holomorphic, the form  $\phi^* \omega$ is of Hodge type $(1,1)$
if and only if $\sigma_Q^*\omega$ is of type $(1,1)$.
By \ref{_B_2_=_holo_on_Tw_Claim_} (i), this happens
if and only if $\omega$ is a $B_2$-form. This proves
an equivalence (i) $\Leftrightarrow$ (ii). Since 
the choice of $I$ is arbitrary, the pullback
$\phi^* \omega$ of a $B_2$-form is 
of Hodge type $(1,1)$
with respect to all induced complex structures.
By \ref{_SU(2)_inva_type_p,p_Lemma_}, this proves the implication
(i) $\Rightarrow$ (iii). The implication
(iii) $\Rightarrow$ (ii) is clear.
\endproof

\hfill

\proposition
Let $(B, \nabla)$ be a complex vector bundle with connection
over $Q$, and $(\phi^*B, \phi^*\nabla)$ its pullback to $M$.
Then the following conditions are equivalent
\begin{description}
\item[(i)]  $(B, \nabla)$ is a \index{terms}{$B_2$-bundles} $B_2$-form
\item[(ii)] The curvature of $(\phi^*B, \phi^*\nabla)$ is of 
Hodge type (1,1)
with respect to some induced complex structure $I$ on $M$
\item[(iii)] The bundle $(\phi^*B, \phi^*\nabla)$ is autodual
\index{terms}{Swann's formalism!for vector bundles}
\end{description}

{\bf Proof:} Follows from \ref{_phi^*_B_2-forms_1,1_Lemma_}
applied to $\omega = \nabla^2$. \endproof

\hfill

For any point $I\in \C P^1$,
consider the corresponding algebra embedding 
$\C \stackrel {c_I} \hookrightarrow {\Bbb H}$. Let 
$\rho_I$ be the action of $\C^*$ on $(M, I)$ obtained
as a restriction of $\rho$ to $c_I(\C^*)\subset {\Bbb H}^*$.
Clearly from \ref{_H^*_specia_Definition_} (ii),
$\rho_I$ acts on $(M, I)$ by holomorphic automorphisms.

Consider $\Tw(M)$ as a union
\[ \Tw(M)  = \bigcup_{I\in \C P^1} 
   \pi^{-1}(I), \ \ \pi^{-1}(I)= (M, I) 
\]
Gluing $\rho(I)$ together, we obtain a 
\index{terms}{twistor space!$\C^*$-equivariant}
smooth \index{terms}{$\C^*$-action on twistor spaces} 
$\C^*$-action $\rho_\C$ on $\Tw(M)$.

\hfill

\claim \label{_C^*_acti_on_Tw_holo_Claim_}
Consider the action $\rho_\C:\; \C^*\times \Tw(M) \arrow \Tw(M)$
constructed above. Then $\rho_\C$ is holomorphic.

{\bf Proof:} It is obvious from construction that
$\rho_\C$ is compatible with the complex structure
on $\Tw(M)$. \endproof

\hfill

\example
Let $M = {\Bbb H}^n \backslash 0$. Since $\Tw({\Bbb H}^n)$ is 
canonically isomorphic to a total space
of the bundle $\calo(1)^n$ over $\C P^1$, the twistor
\index{terms}{twistor space!of a flat manifold} 
space $\Tw(M)$ is $\Tot(\calo(1)^n)$ without zero section.
The group $\C^*$ acts on $\Tot(\calo(1)^n)$ by dilatation,
and the restriction of this action to $\Tw(M)$
coincides with $\rho_\C$.

\hfill

Consider the map $\sigma:\; \Tw(M) \arrow M$.
Let $(B, \nabla)$ be a $B_2$-bundle over $Q$.
\index{terms}{$B_2$-bundles!and $\C^*$-equivariant holomorphic bundles} 
Since the bundle $(\phi^*B, \phi^*\nabla)$
is autodual, the curvature of $\sigma^*\phi^*\nabla$
has type $(1,1)$. Let
$(\sigma^* \phi^*B, (\sigma^*\phi^*\nabla)^{0,1})$ 
be the holomorphic bundle obtained from $(\phi^*B, \phi^*\nabla)$
by twistor \index{terms}{twistor transform} transform. Clearly, this bundle 
is $\C^*$-equivariant, with respect to the natural
\index{terms}{$\C^*$-equivariant holomorphic bundles on twistor spaces} 
\index{terms}{$\C^*$-action on twistor spaces} 
$\C^*$-action on $\Tw(M)$. It turns out that any
$\C^*$-equivariant bundle $\c F$ on $\Tw(M)$
can be obtained this way, assuming that
$\c F$ is compatible with twistor transform.

\hfill

\theorem\label{_B_2_to_C^*_equiva_Theorem_}
In the above assumptions, let $\c C_{B_2}$ 
be the category of of $B_2$-bundles on $Q$,
\index{terms}{$B_2$-bundles!and $\C^*$-equivariant holomorphic bundles} 
and $\c C_{\Tw, \C^*}$ the category of $\C^*$-equivariant
\index{terms}{$\C^*$-equivariant holomorphic bundles on twistor spaces} 
holomorphic bundles on $\Tw(M)$ which are compatible
with the twistor \index{terms}{twistor transform} transform. Consider 
the functor 
\[ (\sigma^* \phi^*)^{0,1}:  \c C_{B_2}\arrow \c C_{\Tw, \C^*},
\]
$(B, \nabla)\arrow (\sigma^* \phi^*B, (\sigma^*\phi^*\nabla)^{0,1})$,
constructed above. Then $(\sigma^* \phi^*)^{0,1}$
establishes an equivalence of categories.

\hfill

We prove \ref{_B_2_to_C^*_equiva_Theorem_}
in Subsection \ref{_twi_tra_H^*_Subsection_}.

\hfill

\remark
Let $Q$ be an arbitrary quaternionic-K\"ahler manifold,
\index{terms}{quaternionic-K\"ahler manifolds!hyperk\"ahler manifolds associated with} 
and $M= \c U(Q)$ the corresponding fibration. Then
$M$ is hypercomplex, and its twistor space is equipped with
a natural holomorphic action of $\C^*$. This gives necessary
ingredients needed to state \ref{_B_2_to_C^*_equiva_Theorem_}
for $Q$ with negative scalar curvature. 
The proof which we give for $Q$ with positive scalar curvature
will in fact work for all quaternionic-K\"ahler manifolds.

\hfill

\question
What happens with this construction when $Q$ is a hyperk\"ahler manifold?

\hfill

In this paper, we need \ref{_B_2_to_C^*_equiva_Theorem_}
only in the case $Q={\Bbb H}P^n$, $M = {\Bbb H}^n\backslash 0$,
but the general proof is just as difficult.

\hfill

\subsection{$\C^*$-equivariant bundles and twistor transform}
\label{_C^*_equiva_and_twistor_Subsection_}

Let $M$ be a hyperk\"ahler manifold, and $\Tw(M)$ its twistor space.
Recall that $\Tw(M)=\C P^1\times M$ 
is equipped with a canonical anticomplex involution $\iota$,
which acts as identity on $M$ and as central symmetry
$I\arrow -I$ on $\C P^1= S^2$.

\hfill

\proposition \label{_conne_flat_along_leave_C^*_Proposition_}
Let $M$ be a hyperk\"ahler manifold, and $\Tw(M)$ its twistor space.
Assume that $\Tw(M)$ is equipped with a free holomorphic action 
$\rho(z):\; \Tw(M)\arrow \Tw(M)$ of $\C^*$, acting along the
fibers of $\pi:\; \Tw(M)\arrow \C P^1$. Assume, moreover, that
$\iota\circ \rho(z) = \rho(\bar z) \circ \iota$,
\index{terms}{$\C^*$-action on twistor spaces}
\index{terms}{twistor space!$\C^*$-equivariant} 
where $\iota$ is the natural anticomplex involution of $\Tw(M)$.%
\footnote{These assumptions are automatically satisfied 
when $M$ is equipped with a special ${\Bbb H}^*$-action,
and $\rho(z)$ is the corresponding  $\C^*$-action on
$\Tw(M)$.}
Let $\c F$ be a $\C^*$-equivariant holomorphic vector
\index{terms}{$\C^*$-equivariant holomorphic bundles on twistor spaces} 
bundle on $\Tw(M)$. Assume that $\c F$ is compatible
with the twistor \index{terms}{twistor transform} transform. Let $\nabla_{\c F}$ be the
natural connection on $\c F$ (\ref{_cano_conne_Remark_}). 
Then $\nabla_{\c F}$ is flat along the leaves of
$\rho$.

\hfill

{\bf Proof:} First of all, let us 
recall the construction of
the natural connection $\nabla_{\c F}$. 
Let $\c F$ be an arbitrary bundle compatible
with the twistor transform. We construct 
$\nabla_{\c F}$ in terms of the
isomorphism $\Psi_{1,2}$ defined 
in \ref{_F_1_=F_2_Lemma_}.

\hfill

Consider an induced complex structure $I$.
Let $F_I$ be the restriction of $\c F$ to
$(M, I)= \pi^{-1}(I)\subset \Tw(M)$. 
Consider the evaluation map \[ p_I:\; \Lin(M) \arrow (M, I)\]
(Subsection \ref{_twi_lines_C^*_Subsection_}). 
In a similar way we define the holomorphic vector bundle
$F_{-I}$ on $(M, -I)$ and the map 
$p_{-I}:\; \Lin(M) \arrow (M, -I)$.
Denote by $F_1$, $F_{-1}$ the sheaves
$p_I^*(F_I)$, $p_{-I}^*(F_{-I})$. 
In \ref{_F_1_=F_2_Lemma_}, we constructed
an isomorphism $\Psi_{1,-1}:\; F_1 \arrow F_{-1}$.

Let us identify $\Lin(M)$ with $(M, I)\times (M, I)$
(this identification is naturally defined in a neighbourhood
of $\Hor\subset \Lin(M)$  -- see 
\ref{_Lin_is_MxM_Proposition_}). Then
the maps $p_I$, $p_{-I}$ became projections
to the relevant components. Let 
\begin{equation*}
\begin{split}
 \bar\6:\; F_1 &\arrow F_1 \otimes p_I^* \Omega^1(M, -I), \\
 \6:\; F_{-1} &\arrow F_{-1} \otimes p_{-I}^* \Omega^1(M, I), 
\end{split}
\end{equation*}
be the sheaf maps obtained as pullbacks of de Rham differentials
(the tensor product is taken in the category of coherent sheaves
over $\Lin(M)$).
Twisting $\6$ by an isomorphism $\Psi_{1,-1}:\; F_1 \arrow F_{-1}$,
we obtain a map
\[ \6^\Psi:\; F_1 \arrow F_1 \otimes p_I^* \Omega^1(M, I).
\]
Adding $\bar \6$ and $\6^\Psi$, we obtain 
\[ 
\nabla:\; F_1 \arrow F_1 \otimes 
  \bigg (p_I^* \Omega^1(M, I) \oplus p_I^* \Omega^1(M, -I)\bigg).
\]
Clearly, $\nabla$ satisfies the Leibniz rule.
Moreover, the sheaf 
$p_I^* \Omega^1(M, I) \oplus p_I^* \Omega^1(M, -I)$
is naturally isomorphic to the sheaf of differentials over
\[ \Lin(M) = (M,I)\times (M, -I).\] 
Therefore, $\nabla$ can be considered as a connection
in $F_1$, or as a real analytic connection in a 
real analytic complex vector bundle underlying $F_I$.
From the definition of $\nabla_{\c F}$ (\cite{_NHYM_}),
it is clear that $\nabla_{\c F}\restrict{(M, I)}$
equals $\nabla$.

\hfill

Return to the proof of \ref{_conne_flat_along_leave_C^*_Proposition_}.
Consider a $\C^*$-action 
\index{terms}{$\C^*$-action on hyperk\"ahler manifolds} 
$\rho_I(z)$ on $(M, I)$, $(M, -I)$
induced from the natural embeddings 
$(M,I)\hookrightarrow \Tw(M)$, $(M,-I)\hookrightarrow \Tw(M)$.
Then $F_I$ is a $\C^*$-equivariant bundle. 
\index{terms}{$\C^*$-equivariant holomorphic bundles on twistor spaces} 
Since $\iota\circ \rho(z) = \rho(\bar z) \circ \iota$,
the identification  $\Lin(M)=(M, I)\times (M, I)$
is compatible with $\C^*$-action.
Let ${\mathbf r}= \frac{d}{dr}$ be the holomorphic
vector field on $(M, I)$ corresponding to the $\C^*$-action.
To prove \ref{_conne_flat_along_leave_C^*_Proposition_},
we have to show that the operator
\[ [\nabla_{\mathbf r},\nabla_{\bar {\mathbf r}}]:\; 
   F_I \arrow F_I\otimes \Lambda^{1,1}(M, I)
\]
vanishes.

Consider the equivariant structure operator
\[ \rho(z)^F: \rho_I(z)^* F_I \arrow F_I. \]
Let $U$ be a $\C^*$-invariant Stein subset of $(M, I)$.
Consider $\rho(z)^F$ an an endomorphism of the
space of global holomorphic sections $\Gamma_U(F_I)$.
Let
\[ D_r(f):= \lim\limits_{\epsilon\rightarrow 0}
   \frac{\rho_I(1+\epsilon)}{\epsilon},
\]
for $f\in \Gamma_U(F_I)$.
Clearly, $D_r$ is a well defined sheaf endomorphism of $F_I$,
satisfying
\[ D_r(\alpha \cdot f) = 
   \frac{d}{dr}\alpha \cdot f +\alpha \cdot D_r(f),
\]
for all $\alpha\in \calo_{(M,I)}$.
We say that a holomorphic section $f$ of $F_I$ is {\bf $\C^*$-invariant}
if $D_r(f)=0$. 
Clearly, the $\calo_{(M,I)}$-sheaf $F_I$
is generated by $\C^*$-invariant sections. Therefore,
it suffices to check the equality
\[ [\nabla_{\mathbf r},\nabla_{\bar {\mathbf r}}] (f)=0
\]
for holomorphic $\C^*$-invariant $f\in F_I$.

Since $f$ is holomorphic, we have $\nabla_{\bar {\mathbf r}} f =0$.
Thus,
\[ [\nabla_{\mathbf r},\nabla_{\bar {\mathbf r}}] (f)
   = \nabla_{\bar {\mathbf r}} \nabla_{\mathbf r}(f).
\]
We obtain that \ref{_conne_flat_along_leave_C^*_Proposition_}
is implied by the following lemma.

\hfill

\lemma \label{_conne_on_C^*_inva_Lemma_}
In the above assumptions, let $f$ be a $\C^*$-invariant
section of $F_I$. Then $\nabla_{\mathbf r}(f)=0$.

\hfill

{\bf Proof:} Return to the notation we used in the beginning
of the proof of \ref{_conne_flat_along_leave_C^*_Proposition_}.
Then, $\nabla(f) = \bar \6(f) +\6^\Psi(f)$. Since 
$f$ is holomorphic, $\bar \6(f)=0$, so we 
need to show that $\6^\Psi(f)({\mathbf r})=0$.
By definition of $\6^\Psi$, this is equivalent to
proving that
\[ \6\Psi_{1,-1}(f)({\mathbf r})=0.\]
Consider the 
\index{terms}{$\C^*$-action on hyperk\"ahler manifolds} 
$\C^*$-action on $\Lin(M)$ which is induced by
the $\C^*$-action on $\Tw(M)$. Since the maps $p_I$,
$p_{-I}$ are compatible with the $\C^*$-action, the sheaves
$F_1$, $F_{-1}$ are $\C^*$-equivariant. We can repeat the construction
of the operator $D_r$ for the sheaf $F_{-I}$. This allows
one to speak of holomorphic $\C^*$-invariant sections
of $F_{-I}$. Pick a $\C^*$-invariant Stein subset $U\subset (M, -I)$.
Since the statement of \ref{_conne_on_C^*_inva_Lemma_} is local,
we may assume that $M=U$. Let $g_1, ... , g_n$ be a set of
$\C^*$-invariant sections of $F_I$ which generated $F_I$.
Then, the sections $p_{-I}^*(g_1), ..., p_{-I}^*(g_n)$ generate
$F_{-1}$. Consider the section $\Psi_{1,-1}(f)$ of $F_{-1}$.
Clearly, $\Psi_{1,-1}$ commutes with the natural
$\C^*$-action. Therefore, the section $\Psi_{1,-1}(f)$
is $\C^*$-invariant, and can be written as
\[ \Psi_{1,-1}(f)= \sum \alpha_i p_{-I}^* (g_i), \]
where the functions $\alpha_i$ are $\C^*$-invariant.
By definition of $\6$ we have 
\[ \6\left(\sum \alpha_i p_{-I}^* (g_i)\right) = 
   \sum \6(\alpha_i p_{-I}^* (g_i))+\sum \alpha_i \6 (p_{-I}^* (g_i)).
\]
On the other hand, $g_i$ is a holomorphic section of $F_{-I}$,
so $\6 p_{-I}^* (g_i)=0$. We obtain
\[ \6\left(\sum \alpha_i p_{-I}^* \cdot (g_i)\right) = 
   \sum \6\alpha_i p_{-I}^* (g_i). 
\]
Thus,
\[ \6\Psi_{1,-1}(f)({\mathbf r}) = \sum \frac{\6\alpha_i}{\6 r} 
    p_{-I}^* (g_i),
\]
but since the functions $\alpha_i$ are $\C^*$-invariant,
their derivatives along ${\mathbf r}$ vanish. 
We obtain $\6\Psi_{1,-1}(f)({\mathbf r})=0$. This proves 
\ref{_conne_on_C^*_inva_Lemma_}.
\ref{_conne_flat_along_leave_C^*_Proposition_} 
is proven.
\endproof

\subsection{Twistor transform and the ${\Bbb H}^*$-action}
\label{_twi_tra_H^*_Subsection_}

For the duration of this Subsection, we fix a hyperk\"ahler manifold
$M$, equipped with a special ${\Bbb H}^*$-action $\rho$, and the 
corresponding \index{terms}{quaternionic-K\"ahler manifolds} 
quaternionic-K\"ahler manifold $Q= M/{\Bbb H}^*$.
Denote the natural quotient map by $\phi:\; M\arrow Q$.
Clearly, \ref{_B_2_to_C^*_equiva_Theorem_}
is an immediate consequence of the following theorem.

\hfill

\theorem \label{_C^*_equi_cano_conne_Theorem_}
Let $\c F$ be a $\C^*$-equivariant holomorphic bundle
\index{terms}{$\C^*$-equivariant holomorphic bundles on twistor spaces} 
over $\Tw(M)$, which is compatible with the twistor 
\index{terms}{twistor transform} transform.
Consider the natural connection
$\nabla_{\c F}$ on $\c F$. Then $\nabla_{\c F}$
is flat along the leaves of ${\Bbb H}^*$-action.

\hfill

{\bf Proof:} The leaves of ${\Bbb H}^*$-action are parametrized by
the points of $q\in Q$. Consider such a leaf 
$M_q:=\phi^{-1}(q)\subset M$.
Clearly, $M_q$ is a hyperk\"ahler submanifold in
$M$, equipped with a special action of ${\Bbb H}^*$.
Moreover, the restriction of $\c F$ to $\Tw(M_q)\subset \Tw(M)$
satisfies assumptions of \ref{_C^*_equi_cano_conne_Theorem_}.
To prove that $\nabla_{\c F}$
is flat along the leaves of ${\Bbb H}^*$-action,
we have to show that ${\c F}\restrict{\Tw(M_q)}$
is flat, for all $q$. Therefore, it suffices to prove
\ref{_C^*_equi_cano_conne_Theorem_} for 
$\dim_{\Bbb H} M=1$. 

\hfill

\lemma\label{_cano_conne_flat_on_4-dim_Lemma_}
We work in notation and assumptions of 
\ref{_C^*_equi_cano_conne_Theorem_}. Assume that
$\dim_{\Bbb H} M=1$. Then the connection
$\nabla_{\c F}$ is flat.

\hfill

{\bf Proof:} 
Let $I$ be an induced complex structure, and 
$F_I:= F\restrict{(M,I)}$ the corresponding holomorphic
bundle on $(M, I)$. Denote by $z_I$ the vector field
corresponding to the \index{terms}{$\C^*$-action on twistor spaces} 
$\C^*$-action $\rho_I$ on $(M, I)$.
By definition, the  connection $\nabla\restrict{F_I}$
has $SU(2)$-invariant curvature $\Theta_I$. 
On the other hand, $\Theta_I(z_I, \bar z_I)=0$ by
\ref{_conne_flat_along_leave_C^*_Proposition_}.
Since $\nabla_{\c F}=\sigma^*\nabla$ is a pullback of an autodual connection
$\nabla$ on $M$, its curvature is a pullback of 
$\Theta_I$. In particular, $\Theta= \Theta_I$ is independent from
the choice of induced complex structure $I$.
We obtain that $\Theta(z_I, \bar z_I)=0$ for all
induced complex structures $I$ on $M$. 

Now \ref{_cano_conne_flat_on_4-dim_Lemma_} is implied
by the following linear-algebraic claim.

\hfill

\claim \label{_SU_2_inva_2-form_z_bar_z=0_is_zero_Claim_}
Let $M$ be a hyperk\"ahler manifold equipped with a special
${\Bbb H}^*$-action, $\dim_{\Bbb H}M=1$. Consider the vectors 
$z_I$, $\bar z_I$ defined above. Let $\Theta$ be a smooth 
$SU(2)$-invariant 2-form, such that for all induced complex structures,
$I$, we have $\Theta(z_I, \bar z_I)=0$. Then $\Theta=0$.

\hfill

{\bf Proof:} The proof of
\ref{_SU_2_inva_2-form_z_bar_z=0_is_zero_Claim_}
is an elementary calculation. 
Fix a point $m_0\in M$.
Consider the flat hyperk\"ahler manifold ${\Bbb H}\backslash 0$,
equipped with a natural special action of ${\Bbb H}^*$.
From the definition of a special action, it is clear
that the map $\rho$ defines a covering ${\Bbb H}\backslash 0\arrow M$,
$h\arrow \rho(h)m_0$
of hyperk\"ahler manifolds, and this covering is compatible
with the special action. Therefore, the hyperk\"ahler manifold $M$ is flat,
and the ${\Bbb H}^*$-action is linear in the
flat coordinates.

 Let 
\[ \Lambda^2(M) = \Lambda^+(M) \oplus\Lambda^-(M)
\]
be the standard decomposition of $\Lambda^2(M)$ according
to the eigenvalues of the Hodge $*$ operator. Consider the
natural Hermitian metric on $\Lambda^2(M)$. Then
$\Lambda^-(M)$ is the bundle of $SU(2)$-invariant 2-forms
(see, e. g., \cite{_Verbitsky:Hyperholo_bundles_}), and
$\Lambda^+(M)$ is its orthogonal complement. 
Consider the corresponding orthogonal projection
 $\Pi:\; \Lambda^2(M) \arrow \Lambda^-(M)$. Denote by
$dz_I\wedge d \bar z_I$ the differential form which is dual to
the bivector $z_I\wedge \bar z_I$.
Let $R\subset \Lambda^-(M)$ be the $C^\infty(M)$-subsheaf
of $\Lambda^-(M)$ generated by $\Pi(dz_I\wedge d \bar z_I)$,
for all induced complex structures $I$ on $M$. 
Clearly, $\Theta\in \Lambda^-(M)$ and $\Theta$ is orthogonal
to $R\subset \Lambda^-(M)$. Therefore, to prove that
$\Theta=0$ it suffices to show that $R= \Lambda^-(M)$.
Since $M$ is covered by ${\Bbb H}\backslash 0$,
we may prove $R= \Lambda^-(M)$ in assumption
$M={\Bbb H}\backslash 0$.

\hfill

Let $\gamma$ be the real vector field corresponding to dilatations
of $M={\Bbb H}\backslash 0$, and $d\gamma$ the dual 1-form.
Clearly,
\[ dz_I\wedge d \bar z_I = 2\1 d\gamma\wedge I(d\gamma).
\]
Averaging $d\gamma\wedge I(d\gamma)$ by $SU(2)$, we obtain
\[ \Pi(dz_I\wedge d \bar z_I) =  \1\bigg(
   d\gamma\wedge I(d\gamma) - J(d\gamma) \wedge K(d\gamma)\bigg)
\]
where $I$, $J$, $K$ is the standard triple of generators for
quaternion algebra. Similarly,
\[ \Pi(dz_J\wedge d \bar z_J) =  \1\bigg(
   d\gamma\wedge J(d\gamma) + K(d\gamma) \wedge I(d\gamma)\bigg)
\]
and 
\[ \Pi(dz_K\wedge d \bar z_K) =  \1\bigg(
   d\gamma\wedge K(d\gamma) + I(d\gamma) \wedge J(d\gamma)\bigg)
\]
Thus, $\Pi(R)$ is a 3-dimensional sub-bundle of $\Lambda^-(M)$.
Since $\dim \Lambda^-(M) =3$, we have $\Pi(R)= \Lambda^-(M)$.
This proves \ref{_SU_2_inva_2-form_z_bar_z=0_is_zero_Claim_}.
\ref{_cano_conne_flat_on_4-dim_Lemma_} and
\ref{_C^*_equi_cano_conne_Theorem_} is proven. \endproof

\subsection{Hyperholomorphic sheaves and $\C^*$-equivariant bundles
over $M_\fl$}
\label{_hyperho_shea_C^*_equiv_Y-M_on_blow-up_Subsection_}

Let $M$ be a hyperk\"ahler manifold, $I$ an induced complex structure
and $F$ a \index{terms}{coherent sheaves!reflexive}
reflexive sheaf over $(M, I)$, equipped with a hyperholomorphic
connection. Assume that $F$ has an isolated singularity
in $x\in M$. Consider the sheaf $\c F$ on $\Tw(M)$ 
corresponding to $\c F$ as in the proof of
\ref{_conne_=>_hyperho_Proposition_}. 
Let $s_x\subset \Tw(M)$ be the horizontal \index{terms}{twistor lines} twistor line
corresponding to $x$, and $\goth m$ its ideal.
Consider the associated graded sheaf of $\goth m$.
Denote by $\Tw^{gr}$ the spectre of this associated graded
sheaf. 
 Clearly, $\Tw^{gr}$ is naturally isomorphic to
$\Tw(T_x M)$, where $T_xM$ is the flat hyperk\"ahler manifold
corresponding to the space $T_xM$ with induced quaternion action.
Consider the natural ${\Bbb H}^*$-action on 
$T_xM$. This provides the hyperk\"ahler manifold
$T_xM\backslash 0$ with a special ${\Bbb H}^*$-action.

Let $s_0\subset \Tw^{gr}$ be the horizontal twistor
line corresponding to $s_x$. The space $\Tw^{gr}\backslash s_0$
is equipped with a holomorphic 
\index{terms}{$\C^*$-action on twistor spaces} 
$\C^*$-action (\ref{_C^*_acti_on_Tw_holo_Claim_}). 
Denote by $\c F^{gr}$ the sheaf on $\Tw^{gr}$ associated with 
$\c F$. Clearly, $\c F^{gr}$ is $\C^*$-equivariant.
\index{terms}{$\C^*$-equivariant holomorphic bundles on twistor spaces} 
In order to be able to apply \ref{_B_2_to_C^*_equiva_Theorem_} 
and \ref{_C^*_equi_cano_conne_Theorem_} to
$\c F^{gr}\restrict{\Tw^{gr}\backslash s_0}$,
we need only to show that $\c F^{gr}$
is compatible with twistor transform.

\hfill

\proposition \label{_F^gr_compa_twi_tra_Proposition_}
Let $M$ be a hyperk\"ahler manifold, $I$ an induced complex structure
and $F$ a \index{terms}{coherent sheaves!reflexive}
reflexive sheaf over $(M, I)$, equipped with a hyperholomorphic
connection. Assume that $\c F$ has an isolated singularity
in $x\in M$. Let $\c F^{gr}$ be the $\C^*$-equivariant bundle
on $\Tw^{gr}\backslash s_0$ constructed above. Then 
\begin{description}
\item[(i)] the bundle $\c F^{gr}$
is compatible with twistor \index{terms}{twistor transform} transform. 
\item[(ii)]
Moreover, the natural
connection $\nabla_{\c F^{gr}}$ (\ref{_cano_conne_Remark_})
is Hermitian.
\end{description}

{\bf Proof:} The argument is clear, but cumbersome, and
essentially hinges on taking associate graded quotients everywhere
and checking that all equations remain true. We give a simplified
version of the proof, which omits some details and notation.

Consider the bundle $\c F\restrict {M\backslash s_x}$.
This bundle is compatible with
twistor \index{terms}{twistor transform} transform, and therefore, is equipped with a natural
connection $\nabla_{\c F}$. This connection is constructed
using the isomorphism $\Psi_{1,-1}:\; F_1\arrow F_{-1}$
(see the proof of \ref{_conne_flat_along_leave_C^*_Proposition_}).
We apply the same consideration to $\c F^{gr}\restrict{(T_x M, I)}$, 
and show that the resulting connection $\nabla_{\c F^{gr}}$
is hyperholomorphic. This implies that $\c F^{gr}$
admits a $(1,1)$-connection which is a pullback of
some connection on $\c F^{gr}\restrict{(T_x M, I)}$.
This argument is used to prove 
that $\c F^{gr}$ is compatible with the twistor
\index{terms}{twistor transform} transform.

We use the notation introduced in 
the proof of \ref{_conne_flat_along_leave_C^*_Proposition_}.
Let $\Lin^{gr}$ be the space of twistor maps in $\Tw^{gr}$.
Consider the maps
$p_{\pm I}^{gr}:\; \Lin^{gr}\arrow (T_x M, \pm I)$   
and the sheaves $F^{gr}_{\pm1}:=(p_{\pm I}^{gr})^*\c F^{gr}_{\pm I}$
obtained in the same way as the maps $p_{\pm I}$
and the sheaves $F_{\pm1}$ from the corresponding associated graded 
objects. Taking the associated graded of $\Psi_{1,-1}$
gives an isomorphism $\Psi^{gr}_{1,-1}:\; F^{gr}_1\arrow F^{gr}_{-1}$.
Using the same construction as in the proof of 
\ref{_conne_flat_along_leave_C^*_Proposition_},
we obtain a connection operator
\[ 
\bar \6^{gr}+\6^{\Psi^{gr}}=\nabla^{gr}_I:\; F_1^{gr} \arrow F_1^{gr} \otimes 
  \bigg ((p^{gr}_{-I})^* \Omega^1(T_x M, I) \oplus 
   (p^{gr}_I)^* \Omega^1(T_x M, -I)\bigg).
\]
Since $(\bar \6^{gr})^2 = (\6^{\Psi^{gr}})^2=0$, the
curvature of $\nabla^{gr}_I$ has Hodge type $(1,1)$ with respect to $I$. 
To prove that $\nabla^{gr}_I$ is hyperholomorphic,
we need to show that the curvature of $\nabla^{gr}_I$  
has type $(1,1)$ with respect to every induced
complex structure. Starting from another induced complex
structure $J$, we obtain a connection $\nabla^{gr}_J$,
with the curvature of type $(1,1)$ with respect to $J$.
To prove that $\nabla^{gr}_J$ is hyperholomorphic
it remains to show that $\nabla^{gr}_J=\nabla^{gr}_I$.

Let $\nabla_I$, $\nabla_J$ be the corresponding operators on
$F_1$. From the construction, it is clear that
$\nabla^{gr}_I$, $\nabla^{gr}_J$ are obtained from 
$\nabla_I$, $\nabla_J$ by taking the associated graded quotients.
On the other hand, $\nabla_I = \nabla_J$.
Therefore, the connections $\nabla^{gr}_I$ and $\nabla^{gr}_J$
are equal. We proved that the bundle
$\c F^{gr}\restrict{\Tw^{gr}\backslash s_0}$ is compatible
with the twistor transform.
To prove \ref{_F^gr_compa_twi_tra_Proposition_},
it remains to show that the natural connection
on $\c F^{gr}$ is Hermitian. 

The bundle ${\c F} \restrict{\Tw(M\backslash{x_0})}$ is by definition
Hermitian. Consider the corresponding isomorphism
$\c F \arrow (\iota^*\bar{\c F})^*$ 
(\ref{_twi_tra_for_semili_Proposition_}). Taking an associate
graded map, we obtain an isomorphism 
\[ \c F^{gr}\oldtilde\rightarrow (\iota^*\bar{\c F}^{gr})^*.\]
This gives a non-degenerate semilinear form $h^{gr}$ on $\c F^{gr}$.
It remains only to show that
$h^{gr}$ is pseudo-Hermitian (i. e. satisfies
$h(x, y) = \overline{h(y,x)}$) and positive definite.

\hfill

Let $M^{gr}_\C$ be a \index{terms}{complexification} 
complexification of $M^{gr}=T_x M$,
$M^{gr}_\C = \Lin(M^{gr})$.
Consider the corresponding complex vector bundle 
$\c F^{gr}_\C$ over $M^{gr}_\C$ underlying $\c F^{gr}$. The metric
$h^{gr}$ can be considered as a semilinear form
$\c F^{gr}_\C\times \c F^{gr}_\C\arrow \calo_{M^{gr}_\C}$.
This semilinear form is obtained from the corresponding form
$h$ on $\c F$ by taking the associate graded quotients.
Since $h$ is Hermitian, the form $h^{gr}$ is pseudo-Hermitian.
To prove that $h^{gr}$ is positive semidefinite, we need to show that
for all $f\in \c F^{gr}_\C$, the function
$h^{gr}(f, \bar f)$ belongs to
$\calo_{M^{gr}_\C}\cdot\bar\calo_{M^{gr}_\C}$, 
where $\calo_{M^{gr}_\C}\cdot\bar\calo_{M^{gr}_\C}$
denotes the $\R^{>0}$-semigroup of $\calo_{M^{gr}_\C}$
generated by $x\cdot\bar x$, for all $x\in \calo_{M^{gr}_\C}$.
A similar property for $h$ holds, because $h$ is positive definite.
Clearly, taking associated graded quotient of the
semigroup $\calo_{M_\C}\cdot\bar\calo_{M_\C}$, we obtain
$\calo_{M^{gr}_\C}\cdot\bar\calo_{M^{gr}_\C}$.
Thus,
\[ h^{gr}(f, \bar f)\in 
   \left(\calo_{M_\C}\cdot\bar\calo_{M_\C}\right)^{gr} =
   \calo_{M^{gr}_\C}\cdot\bar\calo_{M^{gr}_\C}
\]
This proves that $h^{gr}$ is positive semidefinite.
Since $h^{gr}$ is non-degenerate, this form
in positive definite. \ref{_F^gr_compa_twi_tra_Proposition_} is proven.
\endproof

\hfill

\remark\label{_exte_conne_conje_Remark_}
Return to the notations of \ref{_desingu_hyperho_Theorem_}.
Consider the bundle $\pi^* F\restrict C$, where $C= {\Bbb P}T_xM$ 
is the blow-up divisor. Clearly, this bundle corresponds to
the graded sheaf $F_I^{gr}= \c F^{gr}\restrict{(M,I)}$
on $(T_xM, I)$. By \ref{_F^gr_compa_twi_tra_Proposition_}
(see also \ref{_C^*_equi_cano_conne_Theorem_}), 
the bundle $\pi^* F\restrict C$ is equipped with a natural
${\Bbb H}^*$-invariant connection and Hermitian structure.%
\footnote{As usually, coherent sheaves over projective variety $X$
correspond to finitely generated graded modules over
the graded ring $\oplus \Gamma(\calo_X(i))$.}
The sheaf $\pi^*F \restrict{\tilde M \backslash C}$
is a hyperholomorphic bundle over $\tilde M \backslash C\cong M\backslash x_0$.
Therefore, $\pi^*F \restrict{\tilde M \backslash C}$
is equipped with a natural metric and a hyperholomorphic 
connection. It is expected that the natural connection and metric
on $\pi^* F\restrict{\tilde M\backslash C}$ can be extended to 
$\pi^* F$, and the rectriction of the resulting connection and metric
to $\pi^* F\restrict C$ coincides with that given by
\ref{_F^gr_compa_twi_tra_Proposition_} and 
\ref{_C^*_equi_cano_conne_Theorem_}. This will give
an alternative proof of \ref{_F^gr_compa_twi_tra_Proposition_} (ii),
because a continuous extension of a positive definite
Hermitian metric is a positive semidefinite Hermitian metric.

\subsection{Hyperholomorphic sheaves and stable bundles
on $\C P^{2n+1}$}

The purpose of the current Section was to prove the following
result, which is a consequence of 
\ref{_F^gr_compa_twi_tra_Proposition_} and
\ref{_C^*_equi_cano_conne_Theorem_}.

\hfill

\theorem \label{_hyperho_blow-up_stable_Theorem_}
Let $M$ be a hyperk\"ahler manifold, $I$ an induced complex
structure and $F$ a \index{terms}{coherent sheaves!reflexive}
reflexive sheaf on $(M,I)$ admitting
a hyperholomorphic connection. Assume that $F$ has an isolated
singularity in $x\in M$, and is locally trivial outside of $x$.
Let $\pi:\; \tilde M\arrow (M,I)$ be the blow-up of $(M,I)$ in $x$.
Consider the holomorphic vector bundle $\pi^* F$ on $\tilde M$
(\ref{_desingu_hyperho_Theorem_}). Let $C\subset (M,I)$ be the 
blow-up divisor, $C= {\Bbb P}T_xM$. Then the holomorphic
bundle $\pi^* F\restrict{ C}$ admits a 
natural Hermitian connection $\nabla$ which is equivariant with
respect to the natural $SU(2)$-action. Moreover,
the holomorphic vector bundle $\pi^* F\restrict{C}$
is a direct sum of stable bundles, with degree which is 
expressed in terms of the equivariant structure.

\hfill

{\bf Proof:} By definition, coherent sheaves on 
$C= {\Bbb P}T_xM$ correspond bijectively to $\C^*$-equivariant
\index{terms}{$\C^*$-equivariant holomorphic bundles on twistor spaces} 
sheaves on $T_xM\backslash 0$. Let $F^{gr}$ be the 
associated graded sheaf of $F$ (Subsection 
\ref{_hyperho_shea_C^*_equiv_Y-M_on_blow-up_Subsection_}).
Consider $F^{gr}$ as a bundle on $T_xM\backslash 0$.
In the notation of \ref{_F^gr_compa_twi_tra_Proposition_},
$F^{gr}= \c F^{gr}\restrict{(M, I)}$. By 
\ref{_F^gr_compa_twi_tra_Proposition_}, the sheaf
$\c F^{gr}$ is $\C^*$-equivariant. 

The bundle $F^{gr}$ is obtained as a pullback
of $\pi^*F\restrict C$, where \[\pi:\; T_xM \backslash 0 \arrow C\]
denotes the standard projection. The vector bundles on $\C P T_xM$
is the same thing as $\C^*$-equivariant vector bundles
on $T_xM\backslash 0$. We obtained two $\C^*$-equivariant
structures on $F^{gr}$: one is given by
\ref{_F^gr_compa_twi_tra_Proposition_},
another is given from the identification
$F^{gr}=\pi^*F\restrict C$.

The quotient of these two equivariant structures
gives a $\C^*$-action $\rho$ on $F^{gr}$. This gives a
decomposition
$F^{gr}=\bigoplus F^{gr}_i(i)$, with 
$\rho(t)\restrict{F^{gr}_i}=t^i$. 

Consider now the
quaternionic projective space ${\Bbb H}P T_xM$. For each 
$p\in{\Bbb H}P T_xM$,
denote by $H_p$ the fiber of the natural projection
\[ T_xM\backslash 0 \arrow {\Bbb H}P T_xM,\] isomorphic
to ${\Bbb H}\backslash 0$. The bundles $F^{gr}_i$ are trivial
on each of the fibers $H_p$. Therefore, these bundles
are compatible with the twistor transform. From
\ref{_C^*_equi_cano_conne_Theorem_}
it follows that each of $F^{gr}_i$ is Yang-Mills and polystable.
\endproof


\section{Moduli spaces of hyperholomorphic sheaves and bundles}
\label{_modu_hyperho_Section_}
\index{terms}{coherent sheaves!hyperholomorphic!moduli of}
\index{terms}{vector bundles!hyperholomorphic!moduli of}
 

\subsection{Deformation of hyperholomorphic sheaves with 
isolated singularities}

The following theorem is an elementary
consequence of \ref{_hyperho_blow-up_stable_Theorem_}.
The proof uses well known results on stability and reflexization
(see, for instance, \cite{_OSS_}).
The main idea of the proof is the following.
Given a family of 
\index{terms}{coherent sheaves!hyperholomorphic!a family of} 
hyperholomorphic sheaves with an
isolated singularity, we blow-up this singularity
and restrict the obtained family to a blow-up
divisor. We obtain a family of coherent sheaves
$\goth V_s$, $s\in S$ over $\C P^{2n+1}$, with fibers 
\index{terms}{coherent sheaves!semistable!a family of} 
semistable
\index{terms}{slope!of a coherent sheaf}  
of slope zero. Assume that for all $s\in S$,
$s\neq s_0$, the sheaf $\goth V_s$ is trivial.
Then the family  $\goth V$ is also trivial,
up to a \index{terms}{reflexization}
reflexization.

\hfill

We use the following property of \index{terms}{coherent sheaves!reflexive}
reflexive sheaves.

\hfill

\definition 
Let $X$ be a complex manifold, and $F$ a torsion-free coherent sheaf.
We say that $F$ is {\bf normal} if 
for all open subvarieties $U\subset X$,
and all closed subvarieties $Y\subset U$ of codimension 2,
the restriction \[ \Gamma_U(F) \arrow \Gamma_{U\backslash Y}(F)\]
is an isomorphism.

\hfill

\lemma\label{_normal_refle_Lemma_}
Let $X$ be a complex manifold, and $F$ a torsion-free coherent sheaf.
Then $F$ is \index{terms}{coherent sheaves!reflexive!are normal}
\index{terms}{coherent sheaves!normal}
reflexive if and only if $F$ is normal.

{\bf Proof:} \cite{_OSS_}, Lemma 1.1.12. \endproof

\hfill

\theorem \label{_reflexi_defo_loca_trivi_Theorem_}
Let $M$ be a hyperk\"ahler manifol, $I$ an induced 
complex structure, $S$ a complex variety 
and $\goth F$ a family of coherent sheaves over
$(M, I)\times S$.  Consider the sheaf 
$F_{s_0}:= {\goth F}\restrict{(M, I)\times \{s_0\}}$.
Assume that the sheaf $F_{s_0}$ is equipped with a filtration
$\xi$.
Let $F_i$, $i= 1, ..., m$ denote the associated graded components of 
$\xi$, and $F_i^{**}$ denote their reflexizations.
Assume that $\goth F$ is locally trivial outside
of $(x_0, s_0)\in (M, I)\times S$. 
Assume, moreover, that all sheaves $F_i^{**}$,
$i= 1, ..., m$ admit a hyperholomorphic 
connection. Then the \index{terms}{reflexization}
reflexization 
${\goth F}^{**}$ is locally trivial.

\hfill

{\bf Proof:} 
Clearly, it suffices to prove 
\ref{_reflexi_defo_loca_trivi_Theorem_} for
${\goth F}$ \index{terms}{coherent sheaves!reflexive}
reflexive.
Let $\tilde X$ be the blow-up of
$(M, I)\times S$ in $\{x_0\} \times S$,
and $\tilde {\goth F}$ the pullback of
$\goth F$ to $\tilde X$. Clearly,
$\tilde X= \tilde M\times S$, where $\tilde M$ is a 
blow-up of $(M,I)$ in $x_0$. Denote by $C\subset \tilde M$ the blow-up
divisor of $\tilde M$. Taking $S$ sufficiently small,
we may assume that the bundle
$\goth F\restrict{\{x_0\}\times(S\backslash\{s_0\})}$
is trivial. Thus, the bundle
$\tilde{\goth F}\restrict {(C\times S)\backslash (C\times \{s_0\})}$,
which is a pullback of 
$\goth F\restrict{\{x_0\}\times(S\backslash\{s_0\})}$
under the natural projection
$(C\times S)\backslash (C\times \{s_0\})
\arrow\{x_0\}\times(S\backslash\{s_0\}$
is trivial.
To prove that $\goth F$ is locally trivial, we have to show that
$\tilde {\goth F}$ is locally trivial,
and that the restriction of
$\tilde {\goth F}$ to $C\times S$ 
is trivial along the fibers of the natural
projection $C \times S\arrow S$. Clearly, to show that
$\tilde {\goth F}$ is locally trivial we need only to prove
that the fiber $\tilde {\goth F}\restrict z$ has constant
dimension for all $z\in C \times S$. 
Thus, $\tilde {\goth F}$ is locally trivial 
if and only if $\tilde {\goth F}\restrict {C\times S}$
is locally trivial. This sheaf is \index{terms}{coherent sheaves!reflexive}
reflexive, since 
it corresponds to an associate graded sheaf of a reflexive sheaf,
in the sense of Footnote to \ref{_exte_conne_conje_Remark_}.
It is non-singular in codimension 2, because all
reflexive sheaves are non-singular in codimension 2
(\cite{_OSS_}, Ch. II, Lemma 1.1.10).
By \ref{_hyperho_blow-up_stable_Theorem_}, the sheaf
$\tilde {\goth F}\restrict {C\times\{s\}}$
is \index{terms}{coherent sheaves!semistable} a direct sum of
stable bundles, with the slopes depending on the
equivariant action. Since the equivariant
action is deformed continuously with the sheaf,
and the corresponding numbers are all 0 when
the sheaf is smooth, they are all 0 in the
singular points too. Then
\ref{_hyperho_blow-up_stable_Theorem_}
implies that it is polystable.
\ref{_reflexi_defo_loca_trivi_Theorem_} is implied
by the following lemma, applied to the
sheaf $\tilde {\goth F}\restrict {C\times S}$.

\hfill

\lemma \label{_F_to_blow-up_stable=>loc_triv_Lemma_}
Let $C$ be a complex projective space, $S$ a complex variety
and $\goth F$ a torsion-free sheaf over $C\times S$.
Consider an open set $U\stackrel j\hookrightarrow C\times S$,
which is a complement of $C\times\{s_0\}\subset C\times S$.
Assume that the sheaf ${\goth F}\restrict U$ is trivial:
${\goth F}\restrict U\cong \calo_U^n$.
Assume, moreover, that $\goth F$ is non-singular in codimension 2,
and the sheaf $\left({\goth F}\restrict{C \times\{s_0\}}\right)^{**}$
is \index{terms}{coherent sheaves!polystable} 
polystable.
\index{terms}{slope!of a coherent sheaf} 
\[ \rk {\goth F}=\rk {\goth F}\restrict {C\times \{s_0\}}. \]
 Then the \index{terms}{reflexization}
reflexization  
${\goth F}^{**}$ of $\goth F$ is a trivial bundle.

\hfill

{\bf Proof:} 
Using induction, 
it suffices to prove \ref{_F_to_blow-up_stable=>loc_triv_Lemma_}
assuming that it is proven for all $\goth F'$ with
$\rk {\goth F}'<\rk {\goth F}$.
 We may also assume that
$S$ is Stein, smooth and 1-dimensional.

\hfill

{\bf Step 1:} {\it We construct an exact sequence 
\[ 0\arrow\goth F_2\arrow\goth F\arrow \im p_{O_1} \arrow 0 \]
of sheaves of positive rank, which, as we prove in Step 3, satisfy assumptions
of \ref{_F_to_blow-up_stable=>loc_triv_Lemma_}.}

\hfill

Consider the pushforward sheaf $j_* \calo_U^n$. 
From the definition of $j_*$, we obtain a 
canonical map
\begin{equation}\label{_embe_tilde_F_to_j_*_Lemma_}
{\goth F}\arrow j_* \calo_U^n,
\end{equation}
and the kernel of this map is a torsion 
subsheaf in ${\goth F}$.

Let $f$ be a coordinate
function on $S$, which vanishes in $s_0\in S$. Clearly,
\[
   j_* \calo_U^n\cong \calo_{C\times S}^n\left[\frac{1}{f}\right].
\]
On the  other hand, the sheaf $\calo_{C\times S}\left[\frac{1}{f}\right]$
is a direct limit of the following diagram:
\[ \calo_{C\times S}^n\stackrel{\cdot f}\arrow
   \calo_{C\times S}^n\stackrel{\cdot f}\arrow
   \calo_{C\times S}^n\stackrel{\cdot f}\arrow ...,
\]
where $\cdot f$ is the injection given by the multiplication by $f$.
Thus, the map \eqref{_embe_tilde_F_to_j_*_Lemma_}
gives an embedding
\[ {\goth F}\stackrel p \hookrightarrow \calo_{C\times S}^n,\]
which is identity outside of $(x_0, s_0)$. 
Multiplying
$p$ by $\frac{1}{f}$ if necessary, we may assume that
the restriction $p\restrict {C\times \{s_0\}}$ is non-trivial.
Thus, $p$ gives a map
\begin{equation} \label{_goth_F_to_calo_Equation_}
   {\goth F}\restrict {C\times \{s_0\}}\arrow
   \calo^n_{C\times \{s_0\}}.
\end{equation}
with image of positive rank.
Since both sides of \eqref{_goth_F_to_calo_Equation_}
are \index{terms}{coherent sheaves!semistable} semistable of slope zero, 
and $\calo^n_{C\times \{s_0\}}$ is po\-ly\-stable,
the map \eqref{_goth_F_to_calo_Equation_}
satisfies the following conditions.
(see \cite{_OSS_}, Ch. II, Lemma 1.2.8 
for details). 

\hfill

\begin{minipage}[m]{0.8\linewidth}
Let $F_1:= \im p\restrict {C\times \{s_0\}}$, 
and $F_2:= \ker p\restrict {C\times \{s_0\}}$. 
Then the \index{terms}{reflexization}
reflexization
of $F_1$ is a trivial bundle $\calo^k_{C\times \{s_0\}}$, 
and $p$ maps $F_1$ to the direct summand 
$O_1'= \calo^k_{C\times \{s_0\}}\subset \calo^n_{C\times \{s_0\}}$.
\end{minipage}

\hfill

Let $O_1 = \calo^k_{C\times S}\subset \calo^n_{C\times S}$
be the corresponding free subsheaf of $\calo^n_{C\times S}$.
Consider the natural projection $\pi_{O_1}$ of
$\calo^n_{C\times S}$ to $O_1$. 
Let $p_{O_1}$
be the composition of $p$ and $\pi_{O_1}$,
$\goth F_1$ the image of $p_{O_1}$, 
and $\goth F_2$ the kernel of $p_{O_1}$.

\hfill

{\bf Step 2:} {\it We show that the sheaves
$\goth F_2$ and $\goth F_1$ and non-singular
in codimension 2.}

\hfill

Consider the exact sequence
\[ Tor^1(\calo_{C\times \{s_0\}}, {\goth F}_1)\arrow
   \goth F_2\restrict {C\times \{s_0\}}\arrow 
   \goth F\restrict {C\times \{s_0\}} \arrow
   {\goth F}_1\restrict {C\times \{s_0\}} \arrow 0
\]
obtained by tensoring the 
sequence \[ 0\arrow\goth F_2\arrow\goth F\arrow {\goth F}_1 \arrow 0\]
with $\calo_{C\times \{s_0\}}$. From this sequence,
we obtain an isomorphism
${\goth F}_1\restrict {C\times \{s_0\}}\cong F_1$.

A torsion-free coherent sheaf over a smooth manifold
is non-singular in codimension 1
(\cite{_OSS_}, Ch. II, Corollary 1.1.8).

Since $\goth F$ is non-singular in codimension 2, 
the restriction $\goth F\restrict {C\times \{s_0\}}$
is non-singular in codimension 1. Therefore, the torsion of
$\goth F\restrict {C\times \{s_0\}}$ has support of
codimension at least 2 in $C\times \{s_0\}$.
Since the sheaf $F_2$  is a subsheaf of
$\goth F\restrict {C\times \{s_0\}}$, its torsion 
has support of codimension at least 2.
Therefore, the singular set of $F_2$ 
has codimension at least 2 in $C\times \{s_0\}$.
The rank of  $F_2$ is by definition equal to $n-k$.

Since $F_1$ has rank $k$, the singular
set of $\goth F_1$ coincides with the singular set of
$F_1$. Since the restriction
${\goth F}_1\restrict {C\times \{s_0\}}=F_1$,
is a subsheaf of a trivial
bundle of dimension $k$ on $C\times \{s_0\}$, it
is torsion-free. Therefore, the singularities of
${\goth F}_1$ have codimension at least 
2 in $C\times \{s_0\}$.

We obtain that the support of 
$Tor^1(\calo_{C\times \{s_0\}}, {\goth F}_1)$
has codimension at least 2 in $C\times \{s_0\}$. Since
the quotient sheaf
\begin{equation}\label{_F_2_quoti_Equation_}
  \goth F_2\restrict {C\times \{s_0\}}\bigg /
   Tor^1(\calo_{C\times \{s_0\}}, {\goth F}_1)\cong F_2
\end{equation}
is isomorphic to the sheaf $F_2$,
this quotient is 
non-singular in codimension 1. Since we proved that 
$F_2$ is non-singular in codimension 1, the sheaf 
$\goth F_2 \restrict {C\times \{s_0\}}$ is also non-singular
in codimension 1, and its rank is equal to the rank 
of $F_2$.

Let $R$ be the union of singular sets
of the sheaves $\goth F_2$, $\goth F$, ${\goth F}_1$.
Clearly, $R$ is contained in $C\times \{s_0\}$,
and $R$ coincides with the set of all
$x\in C\times \{s_0\}$ where the dimension of the 
fiber of the sheaves 
$\goth F_2$, $\goth F$, ${\goth F}_1$ is
not equal to $n-k$, $n$, $k$. We have seen that
the restrictions of $\goth F_2$, ${\goth F}_1$
to ${C\times \{s_0\}}$ have ranks $n-k$, $k$.
Therefore, the singular sets of $\goth F_2$, ${\goth F}_1$ 
coincide with the singular sets of 
$\goth F_2\restrict{C\times \{s_0\}}$, 
${\goth F}_1\restrict{C\times \{s_0\}}$. We have shown that these
singular sets have codimension at least 2 in $C\times \{s_0\}$.
On the other hand, $\goth F$ is non-singular in codimension 2,
by the conditions of \ref{_F_to_blow-up_stable=>loc_triv_Lemma_}.
Therefore, $R$ has codimension at least 3 in $C\times S$.

\hfill

{\bf Step 3:} {\it We check the assumptions of
\ref{_F_to_blow-up_stable=>loc_triv_Lemma_} applied to the sheaves
$\goth F_2$, ${\goth F}_1$.}

\hfill

Since
the singular set of ${\goth F}_1$ has codimension 2 in $C\times \{s_0\}$,
the $\calo_{C\times \{s_0\}}$-sheaf 
$Tor^1(\calo_{C\times \{s_0\}},{\goth F}_1)$ is a torsion sheaf
with support of codimension 2 in $C\times \{s_0\}$. 
By \eqref{_F_2_quoti_Equation_}, the
\index{terms}{reflexization}
reflexization of $\goth F_2\restrict {C\times \{s_0\}}$
coincides with the \index{terms}{reflexization}
reflexization of $F_2$. Thus, the sheaf
$\left(\goth F_2\restrict {C\times \{s_0\}}\right)^{**}$
is \index{terms}{coherent sheaves!semistable} semistable. On the other hand, outside of ${C\times \{s_0\}}$,
the sheaf $\goth F_2$ is a trivial bundle. 
Thus, $\goth F_2$ satisfies assumptions of
\ref{_F_to_blow-up_stable=>loc_triv_Lemma_}.
Similarly, the sheaf $\goth F_1$ is non-singular in
codimension 2, its restriction to
$C\times \{s_0\}$ has trivial reflexization, and
it is free outside of $C\times \{s_0\}$.

\hfill

{\bf Step 4:} {\it We apply induction
and prove \ref{_F_to_blow-up_stable=>loc_triv_Lemma_}.}

\hfill

By induction assumption, the
\index{terms}{reflexization}
reflexization of $\goth F_2$ is isomorphic to 
a trivial bundle $\calo^{n-k}_{C\times S}$.
and reflexization of $\goth F_1$ is
$\calo^k_{C\times S}$. We obtain an exact sequence
\begin{equation}\label{_F_2_to_F_exa_Equation_} 
  0 \arrow \goth F_2\arrow \goth F \arrow {\goth F}_1 \arrow 0,
\end{equation}
where the sheaves $\goth F_2$ and ${\goth F}_1$ have trivial
reflexizations.

Let $V:= C\times S\backslash R$.
Restricting the exact sequence \eqref{_F_2_to_F_exa_Equation_} to $V$,
we obtain an exact sequence
\begin{equation}\label{_F_2_to_F_restri_V_exa_Equation_}
   0 \arrow \calo^{n-k}_{V} \stackrel a \arrow {\goth F}\restrict V
   \stackrel b\arrow \calo^k_{V} \arrow 0.
\end{equation}
Since $V$ is a complement of a codimension-3 complex subvariety in
a smooth Stein domain, the first cohomology of a trivial sheaf on
$V$ vanish. Therefore, the sequence
\eqref{_F_2_to_F_restri_V_exa_Equation_}
splits, and the sheaf ${\goth F}\restrict V$
is a trivial bundle. Consider the pushforward
$\zeta_*{\goth F}\restrict V$, where $\zeta:\; V\arrow C\times S$
is the standard map. Then $\zeta_*{\goth F}\restrict V$
is a \index{terms}{reflexization}
reflexization of ${\goth F}$
(a pushforward of a \index{terms}{coherent sheaves!reflexive!push-forward of}
reflexive sheaf over a
subvariety of codimension 2 or more is
reflexive -- see \ref{_normal_refle_Lemma_}).
On the other hand, since the sheaf ${\goth F}\restrict V$
is a trivial bundle, its push-forward over a
subvariety of codimension at least 2 is also
a trivial bundle over $C\times S$. We proved that
the sheaf ${\goth F}^{**}=\zeta_*{\goth F}\restrict V$ is a trivial bundle 
over $C\times S$. The push-forward
$\zeta_*{\goth F}\restrict V$ coincides with reflexization
of $\goth F$, by \ref{_normal_refle_Lemma_}.
This proves
\ref{_F_to_blow-up_stable=>loc_triv_Lemma_} and  
\ref{_reflexi_defo_loca_trivi_Theorem_}.
\endproof

\subsection{The Maruyama moduli space of coherent sheaves}

This Subsection is a compilation of results of \index{names}{Gieseker, D.} Gieseker and \index{names}{Maruyama, M.} Maruyama
on the moduli of coherent sheaves over projective manifolds.
We follow \cite{_OSS_}, \cite{_Maruyama:Si_}.

To study the \index{terms}{moduli space!of Maruyama} 
moduli spaces of holomorphic bundles and coherent sheaves,
we consider the following definition of stability.

\hfill

\definition
(Gieseker--Maruyama stability)
(\cite{_Gieseker_}, \cite{_OSS_})
Let $X$ be a projective variety, $\calo(1)$ the standard line bundle
and $F$ a torsion-free coherent sheaf. The sheaf $F$ is called
{\bf \index{names}{Gieseker, D.} Gieseker--Maruyama 
\index{terms}{coherent sheaves!stable!Gieseker--Maruyama stable} stable} 
(resp. Gieseker--Maruyama
\index{terms}{coherent sheaves!semistable!Gieseker--Maruyama semistable} 
semistable) if for all coherent subsheaves $E\subset F$
with $0<\rk E<\rk F$, we have
\[ p_F(k) < p_E(k) \ \ (\text{resp.}, \ \ p_F(k) \leq p_E(k))
\]
for all sufficiently large numbers $k\in \Z$. Here
\[ p_F(k) = \frac{\dim \Gamma_X(F\otimes \calo(k))}{\rk F}. 
\]

\hfill

Clearly, Gieseker--Maruyama stability is weaker than the
Mum\-ford-\-Ta\-ke\-mo\-to stability. Every Gieseker--Maruyama
\index{terms}{coherent sheaves!semistable!Jordan-H\"older filtration on} 
\index{terms}{Jordan-H\"older filtration} 
semistable sheaf $F$ has a so-called Jordan-H\"older filtration
$F_0\subset F_1\subset ...\subset F$ with Gieseker--Maruyama
stable successive quotients $F_i/F_{i-1}$. The corresponding 
associated graded sheaf
\[ \oplus F_i/F_{i-1} \]
is independent from a choice of a filtration.
It is called {\bf the associate graded quotient of
the Jordan-H\"older filtration on $F$}.

\hfill

\definition
Let $F$, $G$ be Gieseker--Maruyama
\index{terms}{coherent sheaves!semistable!$S$-equivalent}
 \index{terms}{$S$-equivalence of coherent sheaves}
semistable sheaves on $X$. Then 
$F$, $G$ are called {\bf $S$-equivalent}
if the corresponding 
associate graded quotients
$\oplus F_i/F_{i-1}$, $\oplus G_i/G_{i-1}$
are isomorphic.

\hfill

\definition
Let $X$ be a complex manifold, $F$ a torsion-free sheaf on $X$,
and $Y$ a complex variety. Consider a sheaf $\c F$ on
$X\times Y$ which is flat over $Y$. Assume that for 
some point $s_0\in Y$, the sheaf $\c F\restrict {X\times \{s_0\}}$
is isomorphic to $F$. Then $\c F$ is called {\bf a deformation of
$F$ parametrized by $Y$}. 
We say that a sheaf $F'$ on $X$ is {\bf deformationally
equivalent} to $F$ if for some $s\in Y$, the restriction 
$\c F\restrict {X\times \{s\}}$ is isomorphic to $F'$.
Slightly less formally, such sheaves are called {\bf deformations of
$F$}. If $F'$ is a (semi-)stable bundle, it is called
{\bf a (semi-) stable bundle deformation of $F$.}
\index{terms}{vector bundles!stable!deformations of}
\index{terms}{vector bundles!semistable!deformations of}

\hfill

\remark\label{_Chern_deforma_equal_}
Clearly, the Chern classes of deformationally equivalent
sheaves are equal.

\hfill

\definition \label{_coarse_modu_Definition_}
Let $X$ be a complex manifold, and $F$ a torsion-free sheaf on $X$,
and $\c M_{mar}$ a complex variety. We say that $\c M_{mar}$ is a 
{\bf coarse \index{terms}{moduli space!coarse} 
moduli space of deformations of $F$}
if the following conditions hold.

\begin{description}
\item[(i)] The 
points of $s\in \c M_{mar}$ are in bijective correspondence
with $S$-equi\-va\-lence classes of coherent sheaves $F_s$
which are deformationally
equivalent to $F$.
\item[(ii)] For any flat deformation
$\c F$ of $F$ parametrized by $Y$, there exists a unique
morphism $\phi:\; Y\arrow \c M_{mar}$ such that
for all $s\in Y$, the restriction
$\c F\restrict {X\times \{s\}}$ is $S$-equivalent
to the sheaf $F_{\phi(s)}$
corresponding to $\phi(s)\in \c M_{mar}$.
\end{description}

Clearly, the coarse moduli space is unique.
By \ref{_Chern_deforma_equal_},
the Chern classes of $F_s$ are equal for
all $s\in \c M_{mar}$.

\hfill

It is clear how to define other kinds of 
\index{terms}{moduli space!of semistable bundles} 
moduli spaces.
For instance, replacing the word {\it sheaf} by the word
{\it bundle} throughout \ref{_coarse_modu_Definition_},
we obtain a definition of {\bf the coarse moduli space
of semistable bundle deformations of $F$}. Further on,
we shall usually omit the word ``coarse'' and say
``moduli space'' instead.
\index{terms}{vector bundles!semistable!moduli of} 

\hfill

\theorem\label{_Maruya_exists_}
(\index{names}{Maruyama, M.}Maruyama) Let $X$ be a projective manifold and
$F$ a coherent sheaf over $X$. Then 
the Maruyama \index{terms}{moduli space!of Maruyama} 
moduli space $\c M_{mar}$ of deformations of $F$ exists
and is compact.

{\bf Proof:} See, e. g., \cite{_Maruyama:Si_}. \endproof

\hfill

\subsection{Moduli of hyperholomorphic sheaves and $C$-restricted
comples structures}

\index{terms}{vector bundles!semistable!moduli of} 
Usually, the moduli space of semistable bundle
deformations of a bundle $F$ is not compact. To compactify 
this moduli space, \index{names}{Maruyama, M.} 
Maruyama adds points corresponding to
the deformations of $F$ which are singular 
(these deformations can be non-reflexive
and can have singular \index{terms}{reflexization}
reflexizations). Using the
\index{terms}{hyperk\"ahler desingularization} desingularization theorems for hyperholomorphic sheaves,
we were able to obtain \ref{_reflexi_defo_loca_trivi_Theorem_},
which states (roughly speaking) that a deformation of a 
semistable \index{terms}{vector bundles!hyperholomorphic} 
hyperholomorphic bundle is again a semistable 
bunlde, assuming that all its singularities are isolated. In 
Section \ref{_C_restri_Section_}, we 
showed that under certain conditions,
a deformation of a hyperholomorphic sheaf is
again \index{terms}{coherent sheaves!hyperholomorphic!deformations of} 
hyperholomorphic (\ref{_sheaf_on_C_restr_hyperho_Theorem_}). 
This makes it possible to prove that a deformation of a semistable 
hyperholomorphic bundle is locally trivial.

\hfill

In \cite{_Verbitsky:Hilbert_}, we have shown that
a Hilbert scheme of a K3 surface has no non-trivial \index{terms}{trianalytic subvarieties} trianalytic
subvarieties, for a general hyperk\"ahler 
\index{terms}{hyperk\"ahler structures!generic} structure.

\hfill

\theorem\label{_space_semista_bu_compa_Theorem_}
Let $M$ be a compact hyperk\"ahler manifold
without non-trivial \index{terms}{trianalytic subvarieties} trianalytic subvarieties, $\dim_{\Bbb H}\geq 2$,
and $I$ an induced complex structure.
Consider a 
\index{terms}{vector bundles!hyperholomorphic!moduli of} 
hyperholomorphic bundle $F$ on $(M, I)$
(\ref{_hyperho_shea_Definition_}). Assume that $I$ is a 
\index{terms}{$C$-restricted complex structures} 
$C$-restricted complex structure, $C= \deg_I c_2(F)$. Let
$\c M$ be the moduli space of semistable bundle
\index{terms}{vector bundles!semistable!moduli of} 
deformations of $F$ over $(M, I)$. Then
$\c M$ is compact.

\hfill

{\bf Proof:} The complex structure $I$ is by definition
algebraic, with unique polarization. This makes it possible
to speak of Gieseker--Maruyama stability on $(M, I)$. Denote 
by $\c M_{mar}$ the Maruyama \index{terms}{moduli space!of Maruyama} 
moduli of deformations of $F$.
Then $\c M$ is naturally an open subset of $\c M_{mar}$. 
Let $s\in \c M_{mar}$ be an arbitrary point
and $F_s$ the corresponding coherent sheaf
on $(M, I)$, defined up to $S$-equivalence.
According to \ref{_Chern_deforma_equal_},
the Chern classes of $F$ and $F_s$ are equal.
Thus, by \ref{_sheaf_on_C_restr_hyperho_Theorem_}, 
the sheaf $F_s$ is hyperholomorphic. Therefore,
$F_s$ admits a filtration with hyperholomorphic stable
quotient sheaves $F_i$, $i= 1, ..., m$. By \ref{_singu_triana_Claim_},
the singular set $S$ of $F_s$ is \index{terms}{trianalytic subvarieties} trianalytic. 
Since $M$ has no proper \index{terms}{trianalytic subvarieties} trianalytic subvarieties, 
$S$ is a collection of points. We obtain that
$F_s$ has isolated singularities. Let $\goth F$ be a family
of deformations of $F$, parametrized by $Y$. 
The points $y\in Y$ correspond to deformations
$F_y$ of $F_s$. Assume that for all $y\in Y$,
$y\neq s$, the sheaf $F_y$ is a bundle. 
Since $\c M$ is open in $\c M_{mar}$, such a 
deformation always exists.

The sheaf $F_s$ has isolated singularities and
admits a filtration with hyperholomorphic stable
quotient sheaves. This implies that the family $\goth F$ satisfies the
conditions of \ref{_reflexi_defo_loca_trivi_Theorem_}.
By \ref{_reflexi_defo_loca_trivi_Theorem_},
the \index{terms}{reflexization}
reflexization $\goth F^{**}$ is locally trivial.
To prove that $\c M= \c M_{mar}$, we have to show that 
for all $s\in \c M_{mar}$, the corresponding 
coherent sheaf $F_s$ is locally trivial.
Therefore, to finish the proof of
\ref{_space_semista_bu_compa_Theorem_}, it
remains to prove the following 
algebro-geometric claim.

\hfill

\claim \label{_defo_shea-w-holes-has-holes_Claim_}
Let $X$ be a compact complex manifold, $\dim_\C X>2$,
and $\goth F$ a torsion-free coherent sheaf over
$X\times Y$ which is flat over $Y$.
Assume that the \index{terms}{reflexization}
reflexization of $\goth F$ is locally trivial,
$\goth F$ has isolated singularities, and for
some point $s\in Y$, the restriction of
$\goth F$ to the complement $(X\times Y)\backslash (X\times \{s\})$
is locally trivial. Then the reflexization
$\left(\goth F\restrict X\times\{s\}\right)^{**}$ 
is locally trivial.

\hfill

\remark
We say that a kernel of a map from a bundle
to an Artinian sheaf is {\bf a bundle with holes}.
In slightly more intuitive terms,
\ref{_defo_shea-w-holes-has-holes_Claim_}
states that a flat deformation of a bundle with holes 
is again a bundle with holes, and cannot be smooth,
assuming that $\dim_\C X>2$.

\hfill

{\bf Proof of \ref{_defo_shea-w-holes-has-holes_Claim_}:} 
\ref{_defo_shea-w-holes-has-holes_Claim_} is well known.
Here we give a sketch of a proof.
Consider a coherent sheaf $F_s= \goth F \restrict{X\times \{s\}}$,
and an exact sequence
\[ 0\arrow F_s\arrow F_s^{**} \arrow k\arrow 0, \]
where $k$ is an Artinian sheaf. By definition,
the sheaf $F_s^{**}$ is locally trivial. The flat deformations
of $F_s$ are infinitesimally classified by $Ext^1(F_s, F_s)$.
Replacing $F_s$ by a quasi-isomorphic complex of sheaves
$F_s^{**} \arrow k$, we obtain a spectral sequence converging
to $Ext^\bullet(F_s, F_s)$. In the $E_2$-term of this sequence,
we observe the group
\[ Ext^1(F_s^{**}, F_s^{**})\oplus Ext^1(k, k)
   \oplus  Ext^2(k, F_s^{**})\oplus  Ext^0(F_s^{**}, k).
\]
which is responsible for $Ext^1(F_s, F_s)$.
The term $Ext^1(F_s^{**}, F_s^{**})$ is responsible
for deformations of the bundle $F_s^{**}$,
the term $Ext^0(F_s^{**}, k)$ for the deformations
of the map $F_s^{**} \arrow k$, and
the term $Ext^1(k, k)$ for the deformations
of the Artinian sheaf $k$. Thus,
the term $Ext^2(k, F_s^{**})$ is responsible for
the deformations of $F_s$ which change the dimension
of the cokernel of the embedding $F_s\arrow F_s^{**}$.
We obtain that whenever $Ext^2(k, F_s^{**})=0$,
all deformations of $F_s$ are singular. 
On the other hand, $Ext^2(k, F_s^{**})=0$,
because the $i$-th $Ext$ from the skyscraper to a 
free sheaf on a manifold of dimension
more than $i$ vanishes (this is a basic result of 
\index{names}{Grothendieck, A.} Grothendieck's duality,  \cite{_Hartshorne:Grothendieck's_}).
\endproof


\section{New examples of hyperk\"ahler manifolds}
\label{_new_exa_Section_}


\subsection{Twistor paths}
\label{_twi_paths_Subsection_}

This Subsection contains an exposition and
further elaboration of the results of
\cite{_coho_announce_} concerning
the twistor curves in the 
\index{terms}{moduli space!of complex structures!twistor curves in}
\index{terms}{moduli space!of complex structures} 
moduli space of complex structures on a complex manifold of hyperk\"ahler type.

Let $M$ be a compact manifold admitting a hyperk\"ahler structure.
In \ref{_moduli_hyperka_Definition_}, we 
defined the coarse, marked  moduli space of
complex structures on $M$, denoted by $Comp$. 
For the duration of this section, we fix a 
compact simple hyperk\"ahler manifold $M$, and
its moduli $Comp$.

\hfill

Further on, we shall need the following fact.

\hfill

\claim \label{_simple=hyperho_Claim_}
Let $M$ be a hyperk\"ahler manifold,
$I$ an induced complex structure of general type,
and $B$ a holomorphic vector bundle over $(M, I)$.
Then $B$ is 
\index{terms}{coherent sheaves!stable!on generic hyperk\"ahler manifolds} 
stable if an only if $B$ is 
simple.\footnote{Simple sheaves are coherent sheaves which have
no proper subsheaves}

{\bf Proof:} 
By \ref{_Lambda_of_inva_forms_zero_Lemma_},
for all $\omega \in Pic(M, I)$, we have
$\deg_I(\omega)=0$. Therefore, every subsheaf
of $B$ is destabilising. 
\endproof

\hfill

\remark
In assumptions of \ref{_simple=hyperho_Claim_},
all stable bundles are 
\index{terms}{vector bundles!hyperholomorphic} hyperholomorphic
(\ref{_inva_then_hyperho_Theorem_}).
Therefore, \ref{_simple=hyperho_Claim_}
implies that $B$ is hyperholomorphic
if it is simple.

\hfill

In Subsection \ref{_modu_and_C-restri_Subsection_},
we have shown that every hyperk\"ahler 
\index{terms}{hyperk\"ahler structures!twistor curves corresponding to!definition of} 
structure $\c H$
corresponds to a holomorphic 
embedding \[ \kappa(\c H):\; \C P^1 \arrow Comp, \ \ 
L \arrow (M, L).\]

\definition
A projective line $C \subset  Comp$ is called
{\bf a \index{terms}{twistor curves!definition of} 
twistor curve} if $C= \kappa(\c H)$ for some 
hyperk\"ahler structure $\c H$ on $M$.
\index{terms}{twistor curves}

\hfill

The following theorem was proven in \cite{_coho_announce_}.

\hfill

\theorem \label{_twistor_connect_Theorem_} 
(\cite{_coho_announce_}, Theorem 3.1)
Let $I_1, I_2\in Comp$. Then there exist a sequence of 
intersecting \index{terms}{twistor curves} twistor curves which connect $I_1$ with $I_2$.

\endproof

\hfill

\definition 
Let  $P_0$, ..., $P_n\subset Comp$ be
a sequence of \index{terms}{twistor curves} twistor curves, supplied with an intersection point
$x_{i+1}\in P_i\cap P_{i+1}$ for each $i$. We say that
$\gamma= P_0, ..., P_n, x_1, ..., x_n$ is 
a {\bf 
\index{terms}{twistor paths!definition of} 
twistor path}. Let $I$, $I'\in Comp$.
We say that $\gamma$ is {\bf a twistor path
connecting $I$ to $I'$} if $I\in P_0$ and $I'\in P_n$.
The lines $P_i$ are called {\bf the edges},
and the points $x_i$ {\bf the vertices}
of a 
\index{terms}{twistor paths!vertices of} 
\index{terms}{twistor paths!edges of} 
twistor path.

\hfill

Recall that in \ref{_generic_manifolds_Definition_},
we defined induced complex structures which are 
generic with respect to a hyperk\"ahler structure.

\hfill

Given a \index{terms}{twistor curves!hyperk\"ahler structures associated to} 
twistor curve $P$, the corresponding hyperk\"ahler
\index{terms}{hyperk\"ahler structures} structure $\c H$ is unique (\ref{_hyperka_etale_Theorem_}).
We say that a point $x\in P$ is {\bf of general type},
or {\bf generic with respect
to $P$} if the corresponding complex structure is generic
with respect to $\c H$. 

\hfill

\definition \label{_admi_twi_path_Definition_}
Let $I$, $J\in Comp$ and $\gamma= P_0, ..., P_n$ be a 
\index{terms}{twistor paths!admissible!definition of} 
twistor path from $I$ to $J$, which corresponds to the 
hyperk\"ahler structures 
$\c H_0$, ..., $\c H_n$. We say that $\gamma$ is {\bf admissible}
if all vertices of $\gamma$ are of general type with respect 
to the corresponding edges.

\hfill

\remark 
In \cite{_coho_announce_},  \index{terms}{twistor paths!admissible}
admissible twistor paths were defined slightly differently. 
In addition to the conditions above, we required that
$I$, $J$ are of general type with respect to $\c H_0$, $\c H_n$.

\hfill

\ref{_twistor_connect_Theorem_} proves that 
every two points $I$, $I'$ in $Comp$ are connected
with a \index{terms}{twistor paths} twistor path. Clearly, each twistor path
induces a diffeomorphism $\mu_\gamma:\; (M,I)\arrow (M,I')$.
In \cite{_coho_announce_}, Subsection 5.2,
we studied algebro-geometrical properties of this
diffeomorphism. 

\hfill

\theorem \label{_admi_twi_impli_Theorem_} 
Let $I$, $J\in Comp$, and $\gamma$ be an admissible 
\index{terms}{twistor paths!admissible}
twistor path from $I$ to $J$. Then

\begin{description}

\item[(i)]  There exists a natural 
isomorphism of tensor cetegories
\[ \Phi_{\gamma}:\; Bun_I(\c H_0)\arrow Bun_J(\c H_n),\] 
where $Bun_I(\c H_0)$, $Bun_J(\c H_n)$ are the categories of 
polystable \index{terms}{vector bundles!hyperholomorphic!category of} 
hyperholomorphic vector bundles on $(M, I)$,
$(M, J)$, taken with respect to $\c H_0$,
$\c H_n$ respectively.

\item[(ii)] Let $B\in Bun_I(\c H_0)$ be a stable hyperholomorphic bundle,
and \[ \c M_{I, \c H_0}(B)\] the 
\index{terms}{moduli space!of stable hyperholomorphic bundles} 
moduli of stable deformations of $B$,
where stability is taken with respect to the K\"ahler metric induced
by $\c H_0$. Then $\Phi_{\gamma}$ maps \index{terms}{vector bundles!stable} stable bundles 
which are deformationally equivalent to $B$ to the stable bundles
which are deformationally equivalent to $\Phi_\gamma(B)$.
Moreover, obtained this way bijection 
\[ \Phi_\gamma:\; \c M_{I, \c H_0}(B)\arrow \c M_{J, \c H_n}(\Phi_\gamma(B))\]
induces a real analytic isomorphism of deformation spaces.

\end{description}

{\bf Proof:} \ref{_admi_twi_impli_Theorem_} (i) is 
a consequence of \cite{_coho_announce_}, 
Corollary 5.1. Here we give a sketch of its proof.

Let $I$ be an induced complex structure
of general type. By \ref{_simple=hyperho_Claim_}, 
a bundle $B$ over $(M, I)$
is \index{terms}{vector bundles!stable} stable if and only if it is simple. Thus,
the category $Bun_I(\c H)$ is independent from
the choice of $\c H$ (\ref{_simple=hyperho_Claim_}). 

In \ref{_equi_cate_Theorem_},
we constructed the equivalence of categories 
$\Phi_{I, J}$, which gives the functor
$\Phi_\gamma$ for 
\index{terms}{twistor paths!consisting of a single twistor curve} 
twistor path which
consists of a single \index{terms}{twistor curves} twistor curve. 
This proves \ref{_admi_twi_impli_Theorem_} (i)
for $n=1$. A composition of isomorphisms
$\Phi_{I, J}\circ \Phi_{J, J'}$ is well
defined, because the category $Bun_I(\c H)$
is independent from the choice of $\c H$.
Taking successive compositions of the maps
$\Phi_{I, J}$, we obtain an isomorphism $\Phi_\gamma$.
This proves \ref{_admi_twi_impli_Theorem_} (i).

 The
variety $\c M_{I, \c H}(B)$ is singular hyperk\"ahler
(\cite{_Verbitsky:Hyperholo_bundles_}),
and the variety $\c M_{J, \c H}(B)$ is the same
singular hyperk\"ahler variety, taken with another
induced complex structure. By definition of
singular hyperk\"ahler varieties, this implies
that $\c M_{I, \c H}(B)$, $\c M_{J, \c H}(B)$
are real analytic equivalent, with equivalence 
provided by $\Phi_{I, J}$. This proves
\ref{_admi_twi_impli_Theorem_} (ii).
\endproof

\hfill

For $I\in Comp$, denote by $Pic(M, I)$ 
the group $H^{1,1}(M, I)\cap H^2(M,\Z)$, and by
$Pic(I, \Q)$ the space 
$H^{1,1}(M, I)\cap H^2(M, \Q)\subset H^2(M)$.
Let $Q\subset H^2(M, \Q)$ be a subspace of $H^2(M, \Q)$,
and
\[ Comp_Q:= \{ I\in Comp \;\; | \;\; Pic(I, \Q) =Q\}. \]

\theorem\label{_admi_exi_Theorem_} 
Let $\c H$, $\c H'$
be hyperk\"ahler \index{terms}{hyperk\"ahler structures} structures, and 
$I$, $I'$ be complex structures of general type to and induced by
$\c H$, $\c H'$. Assume that $Pic(I, \Q) = Pic(I', \Q) =Q$,
and $I$, $I'$ lie in the same connected component of $Comp_Q$.
Then $I$, $I'$ can be connected by an admissible path.

{\bf Proof:}  This is \cite{_coho_announce_}, Theorem 5.2.
$\;\;\blacksquare$

\hfill

For a general $Q$, we have no control over the number of connected 
components of $Comp_Q$ (unless global Torelli theorem is proven), 
and therefore we cannot directly apply 
\ref{_admi_exi_Theorem_} to obtain results from algebraic
geometry.\footnote{Exception is a K3 surface, where Torelli holds.
For K3, $Comp_Q$ is connected for all $Q\subset H^2(M, \Q)$.}
However, when $Q=0$, $Comp_Q$ is clearly
connected and dense in $Comp$. This is used to prove 
the following corollary.

\hfill

\corollary \label{_I_conne_w_admi_Corollary_}
Let $I$, $I'\in Comp_0$. Then $I$ can be connected to $I'$
by an admissible twistor path. 
\index{terms}{twistor paths!admissible!existence}

{\bf Proof} This is \cite{_coho_announce_}, Corollary 5.2.
$\blacksquare$

\hfill

\definition\label{_gene_pola_Definition_}
Let $I\in Comp$ be a complex structure, $\omega$ be a K\"ahler form on 
$(M, I)$, and $\c H$ the corresponding hyperk\"ahler
metric, which exists by \index{terms}{Calabi--Yau Theorem} Calabi-Yau theorem. 
Then $\omega$ is called {\bf a generic polarization}
\index{terms}{generic polarization}
if any of the following conditions hold
\begin{description}
\item[(i)] For all $a\in Pic(M, I)$, the degree 
$\deg_\omega(a)\neq 0$, unless $a=0$.
\item[(ii)] For all $SU(2)$-invariant integer classes
$a\in H^2(M, \Z)$, we have $a=0$.
\end{description}
The conditions (i) and (ii) are equivalent
by \ref{_Lambda_of_inva_forms_zero_Lemma_}.

\hfill

\claim\label{_omega_gene_otho_Claim_}
Let $I\in Comp$ be a complex structure, $\omega$ be a K\"ahler form on 
$(M, I)$, and $\c H$ the corresponding hyperk\"ahler structure,
which exists by \index{terms}{Calabi--Yau Theorem} Calabi-\index{names}{Yau, S.-T.}Yau theorem. Then $\omega$ is generic
if and only if for all integer classes $a\in H^{1,1}(M, I)$,
the class $a$ 
is not orthogonal to $\omega$ with respect to the \index{terms}{Bogomolov--Beauville pairing}
Bogomolov-Beauville
pairing.

{\bf Proof:} Clearly, the map $\deg_\omega:\; H^2(M) \arrow \R$
is equal (up to a scalar multiplier) to 
the orthogonal projection onto the line $\R\cdot\omega$. 
Then, \ref{_omega_gene_otho_Claim_}
is equivalent to \ref{_gene_pola_Definition_}, (i).
\endproof

\hfill

From \ref{_omega_gene_otho_Claim_}
it is clear that the set of generic polarizations is a complement
to a countable union of hyperplanes. Thus, generic
polarizations are dense in the K\"ahler cone of $(M, I)$,
for all $I$.

\hfill

\claim \label{_admi_pa_exist_for_gene_pol_Claim_}
Let $I, J\in Comp$, and $a$, $b$ be generic polarizations 
on $(M, I)$. Consider the corresponding hyperk\"ahler 
\index{terms}{hyperk\"ahler structures} structures
$\c H_0$ and $\c H_n$ inducing $I$ and $J$. Then there exists an
admissible twistor path starting from $I, \c H_0$
\index{terms}{twistor paths!admissible!existence}
and ending with $\c H_n, J$.

\hfill

{\bf Proof:} Consider the \index{terms}{twistor curves} twistor curves $P_0$,
$P_n$ corresponding to $\c H_0$, $\c H_n$. Since $a$,
$b$ are generic, the curves $P_0$,
$P_n$ intersect with $Comp_0$.
Applying \ref{_I_conne_w_admi_Corollary_},
we connect the curves $P_0$ and $P_n$ by an admissible path.
\endproof

\hfill

Putting together \ref{_admi_pa_exist_for_gene_pol_Claim_}
and \ref{_admi_twi_impli_Theorem_}, we obtain the following
result.

\hfill

\theorem\label{_iso_Bun_exists_gene_pola_Theorem_}
Let $I$, $J\in Comp$ be complex structures, and $a, b$ be generic
polarizations on $(M, I)$, $(M, J)$. 
Then
\begin{description}

\item[(i)]  There exist an
isomorphism of tensor cetegories
\[ \Phi_{\gamma}:\; Bun_I(a)\arrow Bun_J(a),\] 
where $Bun_I(a)$, $Bun_J(b)$ are the categories of 
polystable 
\index{terms}{vector bundles!hyperholomorphic!category of} 
hyperholomorphic vector bundles on $(M, I)$,
$(M, J)$, taken with respect to 
the hyperk\"ahler structures defined by
the K\"ahler classes $a$, $b$ as in
\ref{_symplectic_=>_hyperkahler_Proposition_}.

\item[(ii)] Let $B\in Bun_I(a)$ be a stable hyperholomorphic bundle,
and \[ \c M_{I, a}(B)\] 
the moduli of stable deformations of $B$,
\index{terms}{moduli space!of stable hyperholomorphic bundles} 
where stability is taken with respect to the 
polarization $a$. Then $\Phi_{\gamma}$ maps stable bundles 
which are deformationally equivalent to $B$ to the stable bundles
which are deformationally equivalent to $\Phi_\gamma(B)$.
Moreover, obtained this way bijection 
\[ \Phi_\gamma:\; \c M_{I, a}(B)\arrow \c M_{J, b}(\Phi_\gamma(B))\]
induces a real analytic isomorphism of deformation spaces.
\end{description}
\endproof

\hfill

\lemma \label{_Phi_of_tange_Lemma_}
In assumptions of \ref{_admi_twi_impli_Theorem_},
let $B$ be a holomorphic tangent bundle of $(M, I)$. Then
$\Phi_\gamma(B)$ is a holomorphic tangent bundle of $(M, J)$.

{\bf Proof:} Clear. \endproof

\hfill

\corollary \label{_mod_of_tange_compa_Corollary_}
Let $I, J\in Comp$ be complex structures, and
$a$, $b$ generic polarizations on $(M, I)$, $(M, J)$. Assume that
the moduli of stable deformations $\c M_{I, a}(T(M, I))$ of 
\index{terms}{moduli space!of stable hyperholomorphic bundles} 
\index{terms}{vector bundles!stable!moduli of} 
the holomorphic tangent bundle $T^{1,0}(M, I)$
is compact. Then the space $\c M_{J, b}(T(M, J))$
is also compact.

{\bf Proof:} 
Let $\gamma$ be the twistor path
of \ref{_admi_pa_exist_for_gene_pol_Claim_}.
By \ref{_Phi_of_tange_Lemma_},
$\Phi_\gamma(T(M, I)) = T(M, J)$. Applying
\ref{_admi_twi_impli_Theorem_}, we 
obtain a real analytic equivalence
from $\c M_{I, a}(T(M, I))$ to
$\c M_{J, b}(T(M, J))$.
\endproof

\subsection{New examples of hyperk\"ahler manifolds}

\theorem\label{_space_sta_bu_compa_hyperka_Theorem_}
Let $M$ be a compact hyperk\"ahler manifold
without non-trivial \index{terms}{trianalytic subvarieties} trianalytic subvarieties, $\dim_{\Bbb H}M\geq 2$,
and $I$ an induced complex structure.
Consider a \index{terms}{vector bundles!hyperholomorphic!moduli of} 
\index{terms}{moduli space!of stable hyperholomorphic bundles} 
hyperholomorphic bundle $F$ on $M$
(\ref{_hyperho_shea_on_M_Definition_}). Let
$F_I$ be the corresponding holomorphic bundle
over $(M, I)$. Assume that $I$ is a 
\index{terms}{$C$-restricted complex structures} $C$-restricted complex structure, $C= \deg_I c_2(F)$. 
Assume, moreover, that all semistable bundle deformations of
\index{terms}{vector bundles!semistable!deformations of} 
$F_I$ are stable.\footnote{This may happen, for instance,
when $\rk F= \dim_\C M=n$, and the number $c_n(F)$ is prime.}
Denote by $\c M_F^I$ the moduli of stable bundle deformations of 
$F_I$ over $(M, I)$. Then 
\begin{description} 
\item[(i)]
the normalization
$\tilde{\c M}_F^I$ is a compact and 
smooth complex manifold equipped with a natural hyperk\"ahler
structure. 
\item[(ii)]
Moreover, for all induced complex  structures $J$ on $M$,
the the variety $\c M_F^J$ is compact, and has a smooth 
normalization $\tilde{\c M}_F^J$, 
which is also equipped with a natural hyperk\"ahler
structure.
\item[(iii)] Finally, the hyperk\"ahler manifolds
$\tilde{\c M}_F^J$, $\tilde{\c M}_F^I$ are naturally isomorphic.
\end{description}
\index{terms}{hyperk\"ahler structures}

{\bf Proof:} The variety $\c M_F^I$ is compact by
\ref{_space_semista_bu_compa_Theorem_}.
In \cite{_Verbitsky:Hyperholo_bundles_}, 
it was proven that the space $\c M_F^I$ of stable
\index{terms}{vector bundles!stable!moduli of} 
\index{terms}{moduli space!of stable hyperholomorphic bundles} 
deformations of $F$ is a singular hyperk\"ahler 
variety (see also \cite{_NHYM_} for an explicit construction
of the twistor space of $\c M_F^I$).
Then \ref{_space_sta_bu_compa_hyperka_Theorem_}
is a consequence of the Desingularization Theorem
for singular hyperk\"ahler varietiess (\ref{_desingu_Theorem_}).
\endproof

\hfill

The assumptions of \ref{_space_sta_bu_compa_hyperka_Theorem_}
are quite restrictive. Using the technique of
\index{terms}{twistor paths} twistor paths, developed in 
Subsection \ref{_twi_paths_Subsection_},
it is possible to prove a more accessible form
of \ref{_space_sta_bu_compa_hyperka_Theorem_}.

\hfill

Let $M$ be a hyperk\"ahler manifold, and $I$, $J$ induced complex structures.
Given an admissible twistor
path from $I$ to $J$, we obtain an equivalence 
$\Phi_\gamma$ between the category of 
\index{terms}{vector bundles!hyperholomorphic!category of} 
hyperholomorphic bundles
on $(M, I)$ and $(M, J)$.

\hfill

\index{terms}{hyperk\"ahler structures!corresponding to generic polarizations} 
\theorem \label{_twi_pa_space_sta_compa_Theorem_}
Let $M$ be a compact simple hyperk\"ahler manifold, $\dim_{\Bbb H}M >1$,
and $I$ a complex structure
on $M$. Consider a generic polarization $a$ on $(M, I)$.
Let $\c H$ be the corresponding hyperk\"ahler structure, and $F$
a hyperholomorphic bundle on $(M, I)$. Fix
 a hyperk\"ahler structure $\c H'$ on $M$ admitting \index{terms}{$C$-restricted complex structures} $C$-restricted
complex structures, such that $M$ has no \index{terms}{trianalytic subvarieties} trianalytic subvarieties
with respect to $\c H'$. Assume that for some
 \index{terms}{$C$-restricted complex structures} $C$-restricted complex structure $J$
induced by $\c H'$, $C=\deg_I c_2(F)$, all admissible twistor 
paths $\gamma$ from $I$ to $J$, and all 
semistable bundles $F'$ which are deformationally
equivalent to $\Phi_\gamma(F)$, the bundle $F'$ is stable.
Then the space of stable deformations of $F$ is compact.
\index{terms}{vector bundles!stable!moduli of} 
\index{terms}{moduli space!of stable hyperholomorphic bundles} 

\hfill

\remark
The space of stable deformations of $F$ is singular hyperk\"ahler
(\cite{_Verbitsky:Hyperholo_bundles_}) and its normalization
is smooth and hyperk\"ahler (\ref{_desingu_Theorem_}).

\hfill

{\bf Proof of \ref{_twi_pa_space_sta_compa_Theorem_}:}
Clearly, $F'$ satisfies assumptions
of \ref{_space_sta_bu_compa_hyperka_Theorem_},
and the moduli space of its stable deformations is compact.
\index{terms}{moduli space!of stable hyperholomorphic bundles} 
\index{terms}{vector bundles!stable!moduli of} 
Since $\Phi_\gamma$ induces a homeomorphism of moduli spaces
(\ref{_admi_twi_impli_Theorem_}),
the space of stable deformations of $F$ is also compact.
\endproof

\hfill

Applying \ref{_twi_pa_space_sta_compa_Theorem_}
to the holomorphic tangent bundle $T(M, I)$,
we obtain the following corollary.

\hfill

\theorem \label{_defo_tange_compact_Theorem_}
Let $M$ be a compact simple hyperk\"ahler manifold,
$\dim_{\Bbb H}(M)>1$. Assume that for a generic
hyperk\"ahler \index{terms}{hyperk\"ahler structures!generic} 
structure $\c H$  on $M$, this manifold admits
no \index{terms}{trianalytic subvarieties} trianalytic subvarieties.\footnote{This assumption holds for
a Hilbert scheme of points on a K3 surface.} Assume, moreover,
that for some \index{terms}{$C$-restricted complex structures} $C$-restricted induced complex structure $I$,
all semistable bundle deformations of
\index{terms}{vector bundles!semistable!deformations of} 
$T(M, I)$ are stable, for $C> \deg_I c_2(M)$.
Then, for all complex structures $J$ on $M$ and 
all generic polarizations $\omega$
on $(M, J)$, the deformation space 
$\c M_{J, \omega}(T(M, J))$ is compact.

{\bf Proof:} Follows from
\ref{_twi_pa_space_sta_compa_Theorem_}
and \ref{_mod_of_tange_compa_Corollary_}.
\endproof

\subsection{How to check that we obtained new examples 
of hyperk\"ahler manifolds?}
\label{_new_exa_F-M_checking_Subsection_}

A. \index{names}{Beauville, A.} Beauville \cite{_Beauville_} 
described two families of compact hyperk\"ahler manifolds, one
obtained as the Hilbert scheme of points on a K3-surface,
another obtained as the Hilbert scheme of a
2-dimensional torus factorized by the free torus action.

\hfill

\conjecture\label{_anti_Beauville_Conjecture}
There exist compact simple hyperk\"ahler manifolds 
which are not isomorphic to deformations of these two fundamental examples. 

\hfill

Here we explain our strategy of a proof of
\ref{_anti_Beauville_Conjecture} using results on 
compactness of the moduli space of 
\index{terms}{moduli space!of stable hyperholomorphic bundles} 
\index{terms}{vector bundles!hyperholomorphic!moduli of} 
hyperholomorphic bundles.

\hfill

The results of this subsection are still in writing, 
so all statements below this line should be considered as conjectures.
We give an idea of a proof for each result and label it 
as ``proof'', but these ``proofs'' are merely sketches.

\hfill

First of all, it is possible to prove the following theorem.

\hfill

\theorem \label{_no_rk-2-bu_on_Hilb_Theorem_}
Let $M$ be a complex K3 surface without automorphisms. Assume that
$M$ is Mumford-Tate generic with respect to some hyperk\"ahler structure. 
Consider the Hilbert scheme $M^{[n]}$ of points on $M$, $n>1$.
Pick a hyperk\"ahler structure $\c H$ on $M^{[n]}$ which is compatible with
the complex structure. Let $B$ be a \index{terms}{vector bundles!hyperholomorphic} hyperholomorphic
bundle on $(M^{[n]}, \c H)$, $\rk B=2$. Then
$B$ is a trivial bundle.

{\bf Proof:} The proof of \ref{_no_rk-2-bu_on_Hilb_Theorem_}
is based on the same ideas as the proof of
\ref{_no_triana_subva_of_Hilb_Theorem_}.
\endproof

\hfill

For a compact complex manifold $X$ of hyperk\"ahler type, denote its
coarse, marked moduli space (\ref{_moduli_hyperka_Definition_})
\index{terms}{moduli space!of complex structures} 
by $Comp(X)$.

\hfill

\corollary\label{_no_rk-2-bu_on_def_Theorem_}
Let $M$ be a K3 surface, 
$I\in Comp(X)$ an arbitrary complex structure
on $X = M^{[n]}$, $n>1$,
and $a$ a generic polarization on $(X, J)$.
Consider the hyperk\"ahler 
\index{terms}{hyperk\"ahler structures!corresponding to generic polarizations} 
structure $\c H$ which
corresponds to ($I$, $a$) 
as in \ref{_symplectic_=>_hyperkahler_Proposition_}.
Let $B$, $\rk B=2$ be a \index{terms}{vector bundles!hyperholomorphic} hyperholomorphic bundle over $(X, \c H)$.
Then $B$ is trivial.

{\bf Proof:} Follows from
\ref{_no_rk-2-bu_on_Hilb_Theorem_}
and \ref{_iso_Bun_exists_gene_pola_Theorem_}.
\endproof

\hfill

\corollary\label{_defo_4-dim_bu_Corollary_}
Let $M$ be a K3 surface, 
$I\in Comp(X)$ an arbitrary complex structure
on $X = M^{[n]}$, $n>1$, and $a$ a generic polarization on $(X, I)$.
Consider the hyperk\"ahler structure $\c H$ which
\index{terms}{hyperk\"ahler structures!corresponding to generic polarizations} 
corresponds to ($J$, $a$)
(\ref{_symplectic_=>_hyperkahler_Proposition_}).
Let $B$, $\rk B\leq 6$ be a stable 
\index{terms}{vector bundles!hyperholomorphic} hyperholomorphic bundle on 
$(X, \c H)$. Assume that the Chern class $c_{\rk B}(B)$ is non-zero.
Assume, moreover, that $I$ is \index{terms}{$C$-restricted complex structures} $C$-restricted, $C = \deg_I(c_2(B))$.
Let $B'$ be a semistable deformation of $B$ over $(X, I)$.
Then $B'$ is stable.

\hfill

{\bf Proof:} Consider the Jordan--H\"older serie for $B'$.
Let $Q_1 \oplus Q_2 \oplus ...$ be the associated graded 
sheaf. By \ref{_sheaf_on_C_restr_hyperho_Theorem_}, 
the stable bundles $Q_i$ are hyperholomorphic.
Since $c_{\rk B}(B)\neq 0$, we have
$c_{\rk Q_i}(Q_i)\neq 0$. Therefore,
the bundles $Q_i$ are non-trivial.
By \ref{_no_rk-2-bu_on_def_Theorem_},
$\rk Q_i >2$. Since all the Chern classes
of the bundles $Q_i$ are $SU(2)$-invariant,
the odd Chern classes of $Q_i$ vanish 
(\ref{_SU(2)_inva_type_p,p_Lemma_}).
Therefore, $\rk Q_i\geq 4$ for all $i$. Since
 $\rk B\leq 6$, we have $i=1$ and 
the bundle $B'$ is \index{terms}{vector bundles!stable} stable.
\endproof

\hfill

Let $M$ be a K3 surface, $X= M^{[i]}$, $i= 2$, $3$ be its second 
or third Hilbert
scheme of points, $I\in Comp(X)$ arbitrary complex structure
on $X$, and $a$ a generic polarization on 
$(X, I)$.  Consider the hyperk\"ahler structure $\c H$ which
\index{terms}{hyperk\"ahler structures!corresponding to generic polarizations} 
corresponds to $J$ and $a$ by
\index{terms}{Calabi--Yau Theorem} Calabi-Yau theorem (\ref{_symplectic_=>_hyperkahler_Proposition_}). 
Denote by $TX$ the tangent bundle of $X$,
considered as a 
\index{terms}{vector bundles!hyperholomorphic!tangent bundle considered as} 
hyperholomorphic bundle. Let $\Def(TX)$ denote the
hyperk\"ahler \index{terms}{hyperk\"ahler desingularization} desingularization of the moduli of stable deformations of $TX$.
\index{terms}{moduli space!of stable hyperholomorphic bundles} 
\index{terms}{vector bundles!stable!moduli of} 
By \ref{_iso_Bun_exists_gene_pola_Theorem_},
the real analytic subvariety underlying
$\Def(TX)$ is independent from the choice of $I$.
Therefore, its dimension is also independent from the choice of $I$.
The dimension of the
deformation space $\Def(TX)$ can be estimated by a direct
computation, for $X$ a Hilbert scheme. We obtain that
$\dim  \Def(TX)> 40$. 

\hfill

\claim \label{_Def_TX_compa_}
In these assumptions, the space $\Def(TX)$
is a compact hyperk\"ahler manifold.

{\bf Proof:} By \ref{_defo_4-dim_bu_Corollary_},
all semistable bundle deformations of $TX$ are stable.
Then \ref{_Def_TX_compa_} is implied by
\ref{_defo_tange_compact_Theorem_}.
\endproof

\hfill

Clearly, deforming the complex structure on $X$, we obtain a
deformation of complex structures on $\Def(TX)$.
This gives a map
\begin{equation} \label{_Comp_X_to_Comp_Def_Equation_}
Comp(X) \arrow Comp(\Def(TX)).
\end{equation}
It is easy to check that the map
\eqref{_Comp_X_to_Comp_Def_Equation_} is complex analytic,
and maps \index{terms}{twistor curves} twistor curves to twistor curves.

\hfill

\claim \label{_maps_pre_tw_curves_Claim_}
Let $X$, $Y$ be hyperk\"ahler manifolds,
and \[ \phi:\; Comp(X) \arrow Comp(Y)\] be a holomorphic
map of corresponding moduli spaces 
\index{terms}{moduli space!of complex structures} 
which maps twistor curves to twistor curves.
Then $\phi$ is locally an embedding.

{\bf Proof:} An elementary argument using the period maps, in
the spirit of Subsection \ref{_modu_and_C-restri_Subsection_}.
\endproof

\hfill

The following result, along with 
\ref{_no_rk-2-bu_on_Hilb_Theorem_}, is the major stumbling block
on the way to proving \ref{_anti_Beauville_Conjecture}.
The other results of this Subsection are elementary or
routinely proven, but the complete proof of 
\ref{_no_rk-2-bu_on_Hilb_Theorem_} and 
\ref{_defo_simple_Theorem_} seems to be difficult.

\hfill

\theorem \label{_defo_simple_Theorem_}
Let $X$ be a simple hyperk\"ahler manifold
without proper \index{terms}{trianalytic subvarieties} trianalytic subvarieties, $B$ a 
\index{terms}{vector bundles!hyperholomorphic!moduli of} 
hyperholomorphic
bundle over $X$, and $I$ an induced complex structure. 
Denote the corresponding holomorphic bundle over $(X, I)$ by $B_I$.
Assume that the space $\c M$ of stable bundle deformations of $B$ is compact.
\index{terms}{vector bundles!stable!moduli of} 
Let $\Def(B)$ be the hyperk\"ahler \index{terms}{hyperk\"ahler desingularization} desingularization of $\c M$.
Then $\Def(M)$ is a simple hyperk\"ahler manifold.

\hfill

{\bf Proof:} Given a decomposition $\Def(M) = M_1\times M_2$,
we obtain a parallel 2-form on $\Omega_1$ on $\Def(B)$,
which is a pullback of the holomorphic \index{terms}{holomorphic symplectic structure} symplectic form on $M_1$.
Consider the space $\c A$ of connections on $B$, which is an
infinitely-dimensional complex analytic Banach manifold. 
Then $\Omega_1$ corresponds to a holomorphic 
2-form $\tilde \Omega_1$ on $\c A$. Since $\Omega_1$
is parallel with respect to the natural connection
on $\Def(B)$, the form $\tilde \Omega_1$ is also a
parallel 2-form on the tangent space to $\c A$,
which is identified with $\Omega^1(X, \End(B))$.
It is possible to prove that this 2-form is obtained
as 
\[ A, B \arrow \int_{Y} \Theta\left(A\restrict Y, 
   B\restrict Y\right)\Vol(Y)
\]
where 
\[ \Theta:\;\Omega^1(Y, \End(B)) \times \Omega^1(Y, \End(B))\arrow
    \calo_Y
\]
is a certain holomorphic 
pairing on the bundle $\Omega^1(Y, \End(B))$,
and $Y$ is a \index{terms}{trianalytic subvarieties} trianalytic subvariety of $X$.
Since $X$ has no \index{terms}{trianalytic subvarieties} trianalytic subvarieties, 
$\tilde \Omega_1$ is obtained from a $\calo_X$-linear
pairing
\[ \Omega^1(X, \End(B))\times \Omega^1(X, \End(B)) \arrow \calo_X.
\]
Using stability of $B$, it is possible to show that
such a pairing is unique, and thus, $\Omega_1$ coincides
with the holomorphic \index{terms}{holomorphic symplectic structure} symplectic form on $\Def(B)$.
Therefore, $\Def(B) = M_1$, and this manifold is simple.
\endproof

\hfill

Return to the deformations of tangent bundles
on $X= M^{[i]}$, $i=2,3$. Recall that the
second Betti number of a Hilbert scheme of 
points on a K3 surface is equal to $23$,
and that of the generalized Kummer variety is
7 (\cite{_Beauville_}).
Consider the map \eqref{_Comp_X_to_Comp_Def_Equation_}.
By \ref{_defo_simple_Theorem_}, the manifold
$\Def(TX)$ is simple. By 
\index{names}{Bogomolov, F. A.} \index{terms}{Bogomolov's theorem!on period mapping}
Bogomolov's theorem (\ref{_Bogomo_etale_Theorem_}),
we have \[ \dim Comp(\Def(TX)) = \dim H^2(\Def(TX)) -2.\]
Therefore, either $\dim H^2(\Def(TX))> \dim H^2(X)=23$,
or the map \eqref{_Comp_X_to_Comp_Def_Equation_} is etale.
In the first case, the second Betti number of 
$\Def(TX)$ is bigger than that of known simple
hyperk\"ahler manifolds, and thus, $\Def(TX)$
is a new example of a simple hyperk\"ahler
manifold; this proves \ref{_anti_Beauville_Conjecture}.
Therefore, to prove \ref{_anti_Beauville_Conjecture},
we may assume that $\dim H^2(\Def(TX))=23$, 
the map \eqref{_Comp_X_to_Comp_Def_Equation_} is etale,
and $\Def(TX)$ is a deformation of a Hilbert scheme of
points on a K3 surface.

\hfill

Consider the universal bundle $\tilde B$ over 
$X\times \Def(TX)$. Restricting $\tilde B$ to
$\{x\} \times \Def(TX)$, we obtain a bundle $B$ on
$\Def(TX)$. Let $\Def(B)$ be the hyperk\"ahler \index{terms}{hyperk\"ahler desingularization} desingularization
of the moduli space of stable deformations of $B$.
\index{terms}{moduli space!of stable hyperholomorphic bundles} 
\index{terms}{vector bundles!stable!moduli of} 
Clearly, the manifold $\Def(B)$ is independent from the
choice of $x\in X$. Taking the generic hyperk\"ahler structure
\index{terms}{hyperk\"ahler structures!generic} 
on $X$, we may assume that the hyperk\"ahler structure $\c H$ on $\Def(TX)$ 
is also generic. Thus, $(\Def(TX),\c H)$ 
admits \index{terms}{$C$-restricted complex structures} $C$-restricted complex structures and has no \index{terms}{trianalytic subvarieties} trianalytic
subvarieties. In this situation,
\ref{_defo_4-dim_bu_Corollary_} implies
that the hyperk\"ahler manifold $\Def(B)$ is compact.
Applying \ref{_maps_pre_tw_curves_Claim_} again,
we obtain a sequence of maps
\[ Comp(X) \arrow Comp(\Def(TX))\arrow Comp(\Def(B)) \]
which are locally closed embeddings. By the same argument
as above, we may assume that the composition
$Comp(X) \arrow Comp(\Def(B))$ is etale, and
the manifold $\Def(B)$ is a deformation of a Hilbert
scheme of points on K3. Using \index{names}{Mukai, S.} 
Mukai's version of Fourier
\index{terms}{Fourier--Mukai transform} 
transform (\cite{_Orlov:K3_}, \cite{_BBH-R_}), 
we obtain an embedding of the corresponding
derived categories of coherent sheaves,
\[  D(X) \arrow D(\Def(TX))\arrow D(\Def(B)).
\]
Using this approach, it is easy to prove that
\[ \dim X\leq \dim \Def(TX)\leq \dim  \Def(B).
\]
 Let $x\in X$ be an arbitrary point.
Consider the complex $C_x \in D(\Def(B))$ 
of coherent sheaves on $\Def(B)$, obtained as a composition
of the Fourier-Mukai transform maps. It is easy to check that
the lowest non-trivial cohomology sheaf of $C_x$ is 
a skyscraper sheaf in a point $F(x)\in \Def(B)$.
This gives an embedding
\[ F:\; X\arrow \Def(B).
\]
The map $F$ is complex analytic for all induced complex structure.
We obtained the following result.

\hfill

\lemma \label{_double_Fou_embedding_Lemma_}
In the above assumptions,
the embedding \[ F:\; X\arrow \Def(B) \]
is compatible with the hyperk\"ahler structure.

\endproof

\hfill

By \ref{_double_Fou_embedding_Lemma_},
the manifold $\Def(B)$ has a \index{terms}{trianalytic subvarieties} trianalytic subvariety
$F(X)$, of dimension 
$0<\dim F(X)< 40< \dim \Def(B)$.
On the other hand, for a hyperk\"ahler 
\index{terms}{hyperk\"ahler structures!generic} 
structure on $X$ generic,
the corresponding hyperk\"ahler structure on
$\Def(B)$ is also generic, so this manifold
has no \index{terms}{trianalytic subvarieties} trianalytic subvarieties. We obtained a contradiction.
Therefore, either $\Def(TX)$ or $\Def(B)$ is a new example of a 
simple hyperk\"ahler manifold. This proves 
\ref{_anti_Beauville_Conjecture}.

\hfill

\hfill

{\bf Acknowledegments:} 
I am grateful to V. Batyrev, A. Beilinson, P. \index{names}{Deligne, P.} Deligne, D. Gaitsgory,
D. Kaledin, D. \index{names}{Kazhdan, D.} Kazhdan, M. Koncevich and T. \index{names}{Pantev, T.} Pantev for valuable 
discussions. My gratitude to
D. Kaledin, who explained to me the results of \cite{_Swann_}.
This paper uses many ideas of our joint work on direct
and inverse twistor \index{terms}{twistor transform} 
transform (\cite{_NHYM_}).

\hfill

{\bf Year 2012.} Eyal Markman found an error in 
the statement of \ref{_hyperho_blow-up_stable_Theorem_}:
as stated, this theorem claimed that a blow-up of a reflexive
sheaf is polystable, while it is in fact a direct sum of stable
bundles of different degrees. This error is based on trivially
wrong (but heuristically obvious) assumption that an $SU(2)$-equivariant
line bundle over $\C P^n$ is flat. The error is corrected now; it is not
required to change anything else in the proofs. I am grateful to 
Eyal for  finding this error.


\end{document}